\pdfoutput=1
\documentclass[a4paper]{amsart}
\usepackage[utf8]{inputenc}
\usepackage[english,french]{babel}
\usepackage[T1]{fontenc}

\usepackage{listings}
\usepackage{color}      % support the colors
\usepackage{indentfirst}        %for the first line format
\usepackage{latexsym,bm}        %for mathematical equations and italic font
\usepackage{array}
\usepackage{extarrows}
\usepackage{subfiles}
\usepackage{amsxtra}
\usepackage{amstext}
\usepackage{xypic}
\usepackage{dsfont} % for \mathds{N}
\usepackage{amsmath,amsthm,amssymb,amscd,amsfonts}
\usepackage{calligra}
\usepackage{mathrsfs}
\usepackage{mathtools}
\usepackage{enumitem}

\usepackage[colorlinks,  linkcolor=blue,
anchorcolor=blue,
citecolor=red
]{hyperref}
\usepackage{cleveref}

\usepackage[mathscr]{eucal}
\usepackage{mathrsfs}
\usepackage{pdfsync}

\newcommand{\pp}{\mathfrak{p}}
\newcommand{\kk}{\mathfrak{k}}
\newcommand{\g}{\mathfrak{g}}
\newcommand{\z}{\mathfrak{z}}
\newcommand{\y}{\mathfrak{y}}
\newcommand{\kt}{\mathfrak{t}}

\newcommand{\kb}{\mathfrak{b}}
\newcommand{\ku}{\mathfrak{u}}
\newcommand{\kn}{\mathfrak{n}}

\newcommand{\kq}{\mathfrak{q}}
\newcommand{\kh}{\mathfrak{h}}

\newcommand{\km}{\mathfrak{m}}
\newcommand{\rd}{d}

\newcommand{\ii}{\sqrt{-1}}

\newcommand{\cE}{\mathcal{E}}
\newcommand{\cI}{\mathcal{I}}

\newcommand{\cP}{\mathcal{P}}

\newcommand{\Z}{\mathbb{Z}}
\newcommand{\C}{\mathbb{C}}
\newcommand{\R}{\mathbb{R}}
\newcommand{\bN}{\mathbb{N}}

\newcommand{\ad}{\mathrm{ad}}

\newtheorem{innercustomthm}{Theorem}
\newenvironment{customthm}[1]
{\renewcommand\theinnercustomthm{#1}\innercustomthm}
{\endinnercustomthm}

\begin{document}

%\linespread{1.6}

\title[full asymptotics of real analytic torsions]{On full asymptotics of 
real analytic torsions for compact 
locally symmetric orbifolds}
\author{Bingxiao LIU}
%\today

\selectlanguage{english}
\begin{abstract}
	We consider a certain sequence of flat vector bundles on a compact locally symmetric 
	orbifold, and we evaluate explicitly the associated asymptotic Ray-Singer 
	real analytic torsion. The basic idea is to computing the heat trace via Selberg's 
	trace formula, so that a key point in this paper is to evaluate the orbital 
	integrals associated with nontrivial elliptic elements. For that 
	purpose, we deduce a geometric localization formula, so that we can rewrite an elliptic 
	orbital integral as a sum of certain identity orbital integrals associated 
	with the centralizer of that elliptic element. The explicit geometric formula of 
	Bismut for semisimple orbital integrals plays an essential role in these 
	computations.
\end{abstract}

\maketitle

\theoremstyle{plain}
\newtheorem{theorem}{Theorem}[subsection]
\newtheorem{lemma}[theorem]{Lemma}
\newtheorem{proposition}[theorem]{Proposition}
\newtheorem{corollary}[theorem]{Corollary}

\theoremstyle{remark}
\newtheorem{definition}[theorem]{Definition}
\newtheorem{remark}[theorem]{Remark}
\newtheorem{example}[theorem]{Example}

%%%%%%%%%%%%%%%%%%%%%%%%%%%%%%%%%%%%%%%%%%%%%%%%%%%%%%%%%%

%%%%%%%%%%%%%%%%%%%%%%%%%%%%%%%%%%%%%%%%%%%%%%%%%%%%%%%%%%

\tableofcontents

\setcounter{page}{1}
\numberwithin{equation}{subsection}
\renewcommand\thesection{\arabic{section}}
\renewcommand\thesubsection{\arabic{section}.\arabic{subsection}}

%%%%%%%%%%%%%%%%%%%%%%%%%%%%%%%%%%%%%%%%%%%%%%%%%%%%%%%%%%%%%%%
\section{Introduction}

Let $(Z,g^{TZ})$ be a closed Riemannian manifold of dimension $m$, and 
let $F\rightarrow Z$ 
be a complex vector bundle equipped with a Hermitian metric $h^{F}$ 
and a flat connection $\nabla^{F,f}$. Let $(\Omega^\cdot(Z,F), 
d^{Z,F})$ be the associated de Rham complex valued in $F$. It is equipped with an 
$L_{2}$-metric induced by $g^{TZ}$, $h^{F}$. Let 
$\mathbf{D}^{Z,F,2}$ be the corresponding de Rham - Hodge Laplacian. The real 
analytic torsion $\mathcal{T}(Z,F)$ is a real valued (graded) spectral invariant 
of $\mathbf{D}^{Z,F,2}$ introduced by Ray and Singer 
\cite{RaySinger1971a,RaySinger1973}. When $Z$ is odd-dimensional and 
$(F,\nabla^{F,f})$ is acyclic, this invariant 
does not depend on the metric data $g^{TZ}$, $h^{F}$. Ray and Singer 
also conjectured that, for unitarily flat vector bundle $F$ 
(i.e., $\nabla^{F,f}h^{F}=0$), this invariant coincides 
with the Reidemeister torsion, a topological invariant associated with 
$(F,\nabla^{F,f})\rightarrow Z$. This conjecture was later proved by Cheeger \cite{MR528965} and M\"{u}ller 
\cite{MR498252}. 
Using the Witten deformation, Bismut and Zhang 
\cite{MR1139837, MR1185803} gave an extension of the 
Cheeger-M\"{u}ller theorem for arbitrary flat vector bundles.

If $Z$ is a compact orbifold, and if $F$ is a flat orbifold vector 
bundle on $Z$, the Ray-Singer analytic torsion $\mathcal{T}(Z,F)$ 
extends naturally to this case (see Definition \ref{def:2.2.5ss20}). In particular, if $F$ is 
acyclic, and if $Z$ as well as all the singular strata have odd 
dimensions, then $\mathcal{T}(Z,F)$ is 
independent of the metric data (see \cite[Corollary 
4.9]{2017arXiv170408369S}). We refer to \cite{MR2140438, 2017arXiv170408369S}, 
etc for more details.

In this paper, we consider a certain sequence of (acyclic) flat vector bundles 
$\{F_{d}\}_{d\in 
\bN}$ on a compact locally symmetric space $Z$, and we study the asymptotic behavior of 
$\mathcal{T}(Z,F_{d})$ as $d\rightarrow +\infty$. When $Z$ is manifold, such question was already studied by M\"{u}ller 
\cite{Muller2012torsion}, by Bismut-Ma-Zhang \cite{MR2838248,BMZ2015toeplitz} and by 
M\"{u}ller-Pfaff \cite{MP2013raysinger, MR3128980}. In particular, 
Bismut-Ma-Zhang \cite{MR2838248,BMZ2015toeplitz} worked on the manifolds 
which are more general than locally symmetric manifolds.  When $Z$ is a 
compact hyperbolic orbifold, such question was studied by Fedosova in \cite{Fedosova2015compact} using the method of harmonic 
analysis. Here, we consider this question for an arbitrary compact 
locally symmetric orbifold (of noncompact type).

Let $G$ be a connected linear reductive Lie group equipped with a Cartan 
involution $\theta\in \mathrm{Aut}(G)$ and an invariant nondegenerate 
symmetric
bilinear form $B$. 
Let $K\subset G$ be the fixed point set of $\theta$, which is a 
maximal compact subgroup of $G$. Put
\begin{equation}
	X=G/K.
	\label{eq:0.1.1fev}
\end{equation}
Then $X$ is a Riemannian symmetric space with the Riemannian metric 
induced from $B$. For convenience, we also assume that $G$ has a compact 
center, then $X$ is of noncompact type.

Now let
$\Gamma\subset G$ be a cocompact discrete subgroup. Set
\begin{equation}
	Z=\Gamma\backslash X.
\end{equation}
Then $Z$ is a compact locally symmetric space. In general, $Z$ is an 
orbifold. Let $\Sigma Z$ denote the orbifold resolution of the singular points 
in $Z$, whose connected components correspond exactly to the nontrivial 
elliptic conjugacy classes of $\Gamma$.

Since $G$ has compact center, the compact form $U$ of $G$ exists 
and is a connected compact linear Lie 
group. If $(E,\rho^{E}, h^{E})$ is a unitary (analytic) representation of $U$, then 
it extends uniquely to a representation of $G$ by unitary trick. This way, 
$F=G\times_K E$ is a vector bundle on $X$ equipped with an invariant 
flat connection $\nabla^{F,f}$ (see Subsection \ref{section3.4} and \eqref{eq:5.5.9bs}) and a unimodular Hermitian metric $h^{F}$ induced 
by $h^{E}$. Moreover, $(F,\nabla^{F,f},h^{F})$ descends to a flat 
Hermitian orbifold 
vector bundle on $Z$, which is still denoted by $(F,\nabla^{F,f},h^{F})$. Let 
$\mathbf{D}^{Z,F,2}$ denote the corresponding de Rham - Hodge Laplacian.

The fundamental rank $\delta(G)$ (or $\delta(X)$) of $G$ (or $X$) is the 
difference of the complex ranks of $G$ and of $K$. As we will see in Theorem \ref{thm:4.1.4paris20}, if 
$\delta(G)\neq 1$, we always have 
\begin{equation}
	\mathcal{T}(Z,F)=0.
	\label{eq:1.0.4fp2}
\end{equation}
If $F$ is defined instead 
by a unitary representation of $\Gamma$, this result is obtained by Moscovici and Stanton
\cite[Corollary 2.2]{MS1991}. If $\Gamma$ is torsion-free,  with $F$ defined via a representation 
of $G$ as above, \eqref{eq:1.0.4fp2} was proved in \cite[Remark 
8.7]{BMZ2015toeplitz} by using Bismut's formula for orbital integrals 
\cite[Theorem 6.1.1]{bismut2011hypoelliptic} (also see \cite[Theorems 5.4 \& 
5.5]{Ma2017bourbaki}). A new proof was given in 
\cite[Proposition 4.2]{MR3128980}(with a correction given in 
\cite[p44]{Matz-Mueller20}). Note that in \cite[Remark 
5.6]{Ma2017bourbaki}, Ma has indicated that, using essentially 
\cite[Theorem 5.4]{Ma2017bourbaki},
the identity \eqref{eq:1.0.4fp2} still holds if $\Gamma$ is not 
torsion-free (i.e., $Z$ is an orbifold), 
which gives us exactly Theorem \ref{thm:4.1.4paris20} in this paper. Due to 
this vanishing result, 
we only need to deal with the case $\delta(G)=1$.

We now describe the sequence of flat vector 
bundles $\{F_{d}\}_{d\in 
\bN}$ which is concerned here. Note that $U$ contains $K$ as a Lie subgroup. Let $T$ be a maximal torus of $K$, and let 
$T_{U}$ be the maximal torus of $U$ containing $T$. Let $\ku$ be the Lie algebra of $U$, and 
let $\kt_{U}\subset \ku$ be the Lie algebra of $T_{U}$. Let $R(\ku,\kt_{U})$ be the associated 
real root system with a system of positive roots
$R^{+}(\ku,\kt_{U})$. Then let $P_{++}(U)\subset \kt_{U}^{\ast}$ denote the set of (real) dominant weights of 
$U$ with respect to the above root system. If $\lambda\in P_{++}(U)$, let 
$(E_{\lambda},\rho^{E_{\lambda}})$ be the irreducible unitary 
representation of $U$ with the highest weight $\lambda$. We extend it to 
a representation of $G$. We require $\lambda$ to be 
nondegenerate, i.e., as $G$-representations, 
$(E_{\lambda},\rho^{E_\lambda})$ is not isomorphic to 
$(E_{\lambda},\rho^{E_{\lambda}}\circ\theta)$. We also take an 
arbitrary $\lambda_{0}\in 
P_{++}(U)$. If $d\in\mathbb{N}$, let 
$(E_{d},\rho^{E_d}, h^{E_{d}})$ be the unitary representation 
of $U$ with highest weight $d\lambda+\lambda_{0}$. By Weyl's 
dimension formula, $\dim E_{d}$ is a polynomial in $d$. This way, we 
get a sequence of (unimodular) flat vector bundles $\{(F_{d},\nabla^{F_{d}},h^{F_{d}})\}_{d\in 
\bN}$ on $X$ 
or on $Z$.

Note that in Subsection \ref{section:lower} (see also \cite[Lemma 4.1]{BV2013torsion}), the nondegeneracy of 
$\lambda$ implies that, for $d$ large enough,
\begin{equation}
	H^{\cdot}(Z,F_{d})=0.
	\label{eq:1.00cohom}
\end{equation}
Furthermore, $\dim Z$ is odd when $\delta(G)=1$. Then for any 
sufficiently large $d$,
$\mathcal{T}(Z,F_{d})$ is independent of the different 
choices of $h^{E_{d}}$ (or $h^{F_{d}}$).

Let $E[\Gamma]$ be the finite set of elliptic 
classes in $\Gamma$. Set $E^{+}[\Gamma]=E[\Gamma]\backslash\{1\}$. The first main 
result in this paper is the following theorem.
\begin{theorem}\label{thm:0.00001s}
	Assume that $\delta(G)=1$. There exists a (real) polynomial $P(d)$ in 
	$d$, and for each $[\gamma]\in 
	E^{+}[\Gamma]$, there exists a nice exponential polynomial 
	$PE^{[\gamma]}(d)$ in $d$ 
	(i.e., a finite sum of the terms of the form $\alpha 
	d^{j}e^{2\pi\sqrt{-1}\beta 
	d}$ with $\alpha\in\C, j\in\bN, \beta\in \mathbb{Q}$, see Definition 
	\ref{def:7.5.1kk20}), such that there 
	exists a constant $c>0$, for $d$ large, we have
	\begin{equation}
		\mathcal{T}(Z,F_{d})=P(d)+\sum_{[\gamma]\in E^{+}[\Gamma]} 
		PE^{[\gamma]}(d)+\mathcal{O}(e^{-cd}).
		\label{eq:1.0.6kkss}
	\end{equation}
	Moreover, the degrees of $P(d)$, $PE^{[\gamma]}(d)$ can be 
	determined in terms of $\lambda$, $\lambda_{0}$.
\end{theorem}

In \cite[Theorem 1.1]{Muller2012torsion}, for a hyperbolic $3$-manifold $Z$, 
M\"{u}ller computed explicitly the leading term of 
$\mathcal{T}(Z,F_{d})$ as $d\rightarrow +\infty$. In \cite{MR2838248,BMZ2015toeplitz}, under a more general setting for a closed 
manifold $Z$, Bismut, Ma 
and Zhang proved that there 
exists a constant $c>0$ such that \cite[Remark 7.8]{BMZ2015toeplitz}
\begin{equation}
	\mathcal{T}(Z,F_{d})=\mathcal{T}_{L_{2}}(Z,F_{d})+\mathcal{O}(e^{-cd}),
	\label{eq:L2torsion}
\end{equation}
where $\mathcal{T}_{L_{2}}(Z,F_{d})$ denotes the 
$L_{2}$-torsion \cite{MR1158345,MATHAI1992369} associated with 
$F_{d}\rightarrow Z$. Moreover, they constructed universally an 
element 
$W\in\Omega^{\bullet}(Z, o(TZ))$ (where $o(TZ)$ denotes the 
orientation bundle of $TZ$) such that if 
$n_{0}=\deg E_{d}$, then
\begin{equation}
	\mathcal{T}_{L_{2}}(Z,F_{d})=d^{n_{0}+1}\int_{Z} W 
	+\mathcal{O}(d^{n_{0}}).
	\label{eq:Wleading}
\end{equation}
The integral of $W$ in the right-hand side of \eqref{eq:Wleading} is called a $W$-invariant. If we 
specialize \eqref{eq:Wleading} for a compact 
locally symmetric manifold $Z$, we get 
\begin{equation}
	\mathcal{T}_{L_{2}}(Z,F_{d})=d^{n_{0}+1}\mathrm{Vol}(Z)[W]^{\mathrm{max}} 
	+\mathcal{O}(d^{n_{0}}).
	\label{eq:Wleading2}
\end{equation}
In 
\cite[Subsection 8.7]{BMZ2015toeplitz}, the explicit computation on $[W]^{\mathrm{max}}$ was carried out  for $G=\mathrm{SL}_{2}(\C)$ to 
recover the result of M\"{u}ller \cite[Theorem 1.1]{Muller2012torsion}.

We now compare \eqref{eq:1.0.6kkss} with \eqref{eq:L2torsion}. If ignoring that $\Gamma$ may act on $X$ non-effectively, we can extend 
the notion of $L_{2}$-torsion to the orbifold $Z$, so that $\mathcal{T}_{L_{2}}(Z,F_{d})$ is still defined in terms of the 
$\Gamma$-trace of the heat operators on $X$. Then $P(d)$ in \eqref{eq:1.0.6kkss} 
is exactly $\mathcal{T}_{L_{2}}(Z,F_{d})$. 
But different from \eqref{eq:L2torsion}, we still have the nontrivial terms 
$PE^{[\gamma]}(d)$, $[\gamma]\in E^{+}[\Gamma]$ in \eqref{eq:1.0.6kkss}. We 
will see, in a refined version of \eqref{eq:1.0.6kkss} stated in 
Theorem \ref{thm:maintheorem}, that $PE^{[\gamma]}(d)$ is essentially 
an linear combination of certain $L_{2}$-torsions for $\Sigma Z$ 
associated with $[\gamma]$ and $\lambda$, $\lambda_{0}$. 
Therefore, we can define an $L_{2}$-torsion for $\Sigma Z$ as 
follows,
\begin{equation}
	\widetilde{\mathcal{T}}_{L_{2}}(\Sigma Z, F_{d})=\sum_{[\gamma]\in E^{+}[\Gamma]} 
	PE^{[\gamma]}(d).
	\label{eq:1.0.11v2}
\end{equation}
Then, as an analogue to \eqref{eq:L2torsion},  we restate our Theorem \ref{thm:0.00001s} as follows.

\begin{customthm}{1.0.1'}\label{thm:1.0.1bis}
	Assume that $\Gamma$ acts on $X$ effectively. For 
	$Z=\Gamma\backslash X$, as $d\rightarrow +\infty$, we have
	\begin{equation}
		\mathcal{T}(Z,F_{d})=\mathcal{T}_{L_{2}}(Z,F_{d})+\widetilde{\mathcal{T}}_{L_{2}}(\Sigma Z, F_{d})+\mathcal{O}(e^{-cd}).
		\label{eq:L2orbifolds}
	\end{equation}
	Moreover, $\mathcal{T}_{L_{2}}(Z,F_{d})$ is a polynomial in $d$, 
	and $\widetilde{\mathcal{T}}_{L_{2}}(\Sigma Z, F_{d})$ is a nice 
	exponential polynomial in $d$. Their leading terms can be determined in terms of $W$-invariants as in \eqref{eq:Wleading2} .
\end{customthm}

To understand better on $\widetilde{\mathcal{T}}_{L_{2}}(\Sigma Z, 
F_{d})$, we need to recall the results of M\"{u}ller and 
Pfaff in \cite{MR3128980} (also in \cite{MP2013raysinger} for 
hyperbolic case) for a compact 
locally symmetric manifold $Z$. They gave a proof to 
\eqref{eq:L2torsion} using Selberg's trace formula, and then showed 
that $\mathcal{T}_{L_{2}}(Z,F_{d})$ is a polynomial in $d$. The 
Theorem \ref{thm:1.0.1bis} here is an extension of their results, 
which shows a nontrivial contribution from $\Sigma Z$. 

Let us give more detail on the results in \cite{MR3128980}. Let $\mathbf{D}^{X,F_{d},2}$ be the 
$G$-invariant Laplacian operator on $X$ which is the lift of $\mathbf{D}^{Z,F_{d},2}$. For $t>0$, let $p_{t}^{X,F_{d}}(x,x')$ denote the heat kernel of $\frac{1}{2}\mathbf{D}^{X,F_{d},2}$ with respect to the Riemannian volume element on $X$. For $t>0$, the identity orbital integral $\mathcal{I}_{X}(E_d,t)$ of $p^{X,F_{d}}_{t}$ is defined as 
\begin{equation}
	\mathcal{I}_{X}(F_d,t)=\mathrm{Tr_{s}}^{\Lambda^{\bullet}(T^{*}_{x}X)\otimes 
	F_{d,x}}[(N^{\Lambda^{\bullet}(T^{*}_{x}X)}-\frac{m}{2})p^{X,F_{d}}_{t}(x,x)],
	\label{eq:04.2.4ppaa}
\end{equation}
where $N^{\Lambda^{\bullet}(T^{*}_{x}X)}$ is the number operator on $\Lambda^{\bullet}(T^{*}_{x}X)$, and the right-hand side of \eqref{eq:04.2.4ppaa} is independent of the choice of $x\in X$. 
Let $\mathcal{MI}_{X}(F_d,s)$, $s\in \C$ denote the Mellin transform 
(see \eqref{eq:7.3.53conf}) of $\mathcal{I}_{X}(F_d,t)$, which is 
holomorphic at $0$. Set
\begin{equation}
	\mathcal{PI}_{X}(F_d)=\frac{\partial}{\partial s}|_{s=0} 
	\mathcal{MI}_{X}(F_d,s).
	\label{eq:04.2.6kk20}
\end{equation}
The $L_{2}$-torsion is defined as
\begin{equation}
	\mathcal{T}_{L^{2}}(Z,F_{d})=\mathrm{Vol}(Z)\mathcal{PI}_{X}(F_d).
	\label{eq:1.0.7mar20}
\end{equation}

Using essentially the Harish-Chandra's Plancherel theorem for 
$\mathcal{I}_{X}(F_d,t)$, M\"{u}ller-Pfaff \cite{MR3128980} managed to show that 
$\mathcal{PI}_{X}(F_d)$ is a polynomial in $d$ (for $d$ large 
enough). Moreover, if $\lambda_{0}=0$, there exists a constant $C_{\lambda}\neq 0$ such that 
\begin{equation}
	\mathcal{PI}_{X}(F_d)=C_{\lambda}d\dim E_{d}+R(d),
	\label{eq:1.0.9parisconf}
\end{equation}
where $R(d)$ is a polynomial in $d$ of degree no greater than $\deg \dim E_{d}$. They also gave concrete formulae for
$C_{\lambda}$ in some model cases \cite[Corollaries 1.4 
\& 1.5]{MR3128980}. 

In Subsection \ref{section7.3}, we use instead an explicit geometric formula of Bismut 
\cite[Theorem 6.1.1]{bismut2011hypoelliptic} for semisimple orbital 
integrals to give a different computation on $\mathcal{PI}_{X}(F_d)$. In 
Subsection \ref{section7.4kk}, we verify that our 
computational results coincide with the ones of M\"{u}ller-Pfaff 
\cite{MR3128980}.

For the orbifold case, i.e., $\Gamma$ contains 
nontrivial elliptic elements, a
key ingredient to Theorem \ref{thm:0.00001s} is to evaluate 
explicitly the elliptic orbital integrals associated with 
$[\gamma]\in E^{+}[\Gamma]$. For that purpose, we make use of the 
full power of Bismut's formula \cite[Theorem 
6.1.1]{bismut2011hypoelliptic}. Note that if $Z$ is a hyperbolic orbifold, i.e. $G=\mathrm{Spin}(1,2n+1)$, the 
result in Theorem \ref{thm:0.00001s} (or Theorem \ref{thm:1.0.1bis}) 
was obtained by Fedosova in \cite[Theorem 
1.1]{Fedosova2015compact}, where she evaluated the elliptic orbital 
integrals using Plancherel's theorem of Harish-Chandra.

In fact, we obtain in this paper a refined 
version of Theorem \ref{thm:0.00001s}, where we give more explicit descriptions 
on the exponential polynomials $PE^{[\gamma]}(d)$ or 
$\widetilde{\mathcal{T}}_{L_{2}}(\Sigma Z, F_{d})$. Before stating 
this refined result, we need to introduce some 
notations and facts.

Fix $k\in T$, and let $X(k)$ denote the fixed point set of $k$ acting on $X$. Then $X(k)$ is a connected symmetric space with $\delta(X(k))=1$. Let $Z(k)^{0}$ be the 
identity component of the centralizer $Z(k)$ of $k$ in $G$. Then 
$X(k)=Z(k)^{0}/K(k)^{0}$ with $K(k)^{0}=Z(k)^{0}\cap K$. Let $U(k)$ 
denote the centralizer of $k$ in $U$ with Lie algebra $\ku(k)\subset 
\ku$ . Then $U(k)^{0}$ is naturally 
a compact form of $Z(k)^{0}$, the triplet
$(X(k),Z(k)^{0},U(k)^{0})$ becomes a smaller version of $(X,G,U)$, 
except that $Z(k)^{0}$ may have noncompact center. Note that $T_{U}$ is also a 
maximal torus of $U(k)^{0}$. We get the following splitting of 
roots
\begin{equation}
	R(\ku,\kt_{U})=R(\ku(k),\kt_{U})\cup 
	R(\ku^{\perp}(k),\kt_{U}),
	\label{eq:mmmmtttt}
\end{equation}
where $\ku^{\perp}(k)$ is the 
orthogonal space of $\ku(k)$ in $\ku$ with respect to $B$. Let 
$R^{+}(\ku(k),\kt_{U})$, $R^{+}(\ku^{\perp}(k),\kt_{U})$ be 
the induced positive roots, and let $\rho_{\ku}$, 
$\rho_{\ku(k)}$ denote the half of the sum of the roots in 
$R^{+}(\ku,\kt_{U})$, $R^{+}(\ku(k),\kt_{U})$ respectively.

Let $W(\ku_{\C},\kt_{U,\C})$ be the Weyl group associated 
with the pair $(\ku,\kt_{U})$. Put
\begin{equation}
	W^{1}_{U}(k)=\{\omega\in W(\ku_{\C},\kt_{U,\C})\;|\; 
	\omega^{-1}(R^{+}(\ku(k),\kt_{U}))\subset R^{+}(\ku,\kt_{U})\}.
	\label{eq:05.4.13ss20}
\end{equation}
If $\sigma\in W^{1}_{U}(k)$, let $\varepsilon(\sigma)$ 
denote its sign. For $\mu\in P_{++}(U)$, set
\begin{equation}
	\varphi^{U}_{k}(\sigma,\mu)=\varepsilon(\sigma)\frac{\xi_{\sigma(\mu+\rho_{\ku})+\rho_{\ku}}(k)}{\Pi_{\alpha\in R^{+}(\ku^{\perp}(k),\kt_{U})}(\xi_{\alpha}(k)-1)}\in\C^{\ast},
	\label{eq:06.2.7kk20}
\end{equation}
where $\xi_{\alpha}$ is the character of $T_{U}$ with 
(dominant) weight 
$2\pi\sqrt{-1}\alpha$. It is clear that
$\varphi^{U}_{k}(\sigma,d\lambda+\lambda_{0})$ is 
an oscillating term of the form 
$c_{1}e^{2\pi\sqrt{-1}c_{2}d}$ 
with $c_{1}\in \C^{\ast}$, $c_{2}\in\R$. If $k$ is of finite 
order, then $c_{2}\in\mathbb{Q}$.

By an equivalent definition of nondegeneracy in Definition 
\ref{def:nondegdef}, for $\sigma\in W^{1}_{U}(k)$, $\sigma\lambda$ is a nondegenerate dominant weight of 
$U(k)^{0}$ with respect to $\theta|_{Z(k)^{0}}$. Let 
$E^{k}_{\sigma,d}$ denote the unitary representations 
of $U(k)^{0}$ (up to a finite central extension) with highest weight 
$d\sigma\lambda+\sigma(\lambda_{0}+\rho_{\ku})-\rho_{\ku(k)}$, 
$d\in\bN$, and let 
$\{F^{k}_{\sigma,d}\}_{d\in\bN}$ be the corresponding sequence of flat vector bundles on $X(k)$. 	

Now we state our second main theorem, which refines our Theorem 
\ref{thm:0.00001s}.
\begin{theorem}\label{thm:maintheorem}
	Assume that $\delta(G)=1$. 
	\begin{enumerate}
		\item If $\Gamma\subset G$ is a cocompact discrete subgroup, 
		if $\gamma\in \Gamma$ is elliptic, let $S(\gamma)$ denote the finite 
		subgroup of $\Gamma\cap Z(\gamma)$ which acts on $X(\gamma)$ 
		trivially. Then there exists a constant $c>0$, and for each $[\gamma]\in 
		E^{+}[\Gamma]$, there exists a nice exponential polynomial in 
		$d$ denoted by $\cP\cE_{X,\gamma}(F_d)$, such that for $Z=\Gamma\backslash X$, as $d\rightarrow +\infty$, we have
		\begin{equation}
			\begin{split}
				\mathcal{T}(Z, 
				F_{d})=&\frac{\mathrm{Vol}(Z)}{|S(1)|}\cP\cI_{X}(F_d)\\
				&+\sum_{[\gamma]\in 
				E^{+}[\Gamma]}\frac{\mathrm{Vol}(\Gamma\cap Z(\gamma)\backslash X(\gamma))}{|S(\gamma)|}  \cP\cE_{X,\gamma}(F_d)+\mathcal{O}(e^{-cd}).
			\end{split}
			\label{eq:04.2.5kk20}
		\end{equation}
		
		\item Fix an elliptic $[\gamma]\in E^{+}[\Gamma]$, then $\cP\cE_{X,\gamma}(F_d)$ depends only on the conjugacy class of 
		$\gamma$ in $G$ and is independent of the lattice $\Gamma$. 
		If $\gamma$ is conjugate to $k\in T$ by an element in $G$, 
		then we have the following identity
		\begin{equation}
			\mathcal{PE}_{X,\gamma}(F_d)=\sum_{\sigma\in W^{1}_{U}(k)} 
			\varphi^{U}_{k}(\sigma,d\lambda+\lambda_{0}) 
			\mathcal{PI}_{X(k)}(F^{k}_{\sigma,d}),
			\label{eq:1.0.22parisconf}
		\end{equation}		
	\end{enumerate}
\end{theorem}
%%%%%%%%%%%

Theorem \ref{thm:0.00001s} now is just a consequence 
of \eqref{eq:04.2.5kk20}. Note that for 
$[\gamma]\in E^{+}[\Gamma]$, the (compact) orbifold $\Gamma\cap 
Z(\gamma)\backslash X(\gamma)$ represents an orbifold stratum in 
$\Sigma Z$ (see \eqref{eq:3.4.22bbs}, Remark \ref{eq:3.4.3vogel}). An important observation 
on \eqref{eq:04.2.5kk20} is that the sequence $\{\mathcal{T}(Z, 
F_{d})\}_{d\in\bN}$ encodes the volume 
information on $Z$ as well as on $\Sigma Z$. Moreover, combining 
\eqref{eq:1.0.7mar20}, \eqref{eq:04.2.5kk20} with  
\eqref{eq:1.0.22parisconf}, we justify that the quantity 
$\widetilde{\mathcal{T}}_{L_{2}}(\Sigma Z, F_{d})$ defined by 
\eqref{eq:1.0.11v2} is indeed a linear combination of 
$L_{2}$-torsions such as $\mathcal{T}_{L^{2}}(\Gamma\cap 
Z(\gamma)\backslash X(\gamma),F^{\gamma}_{\sigma,d})$ for $\Sigma 
Z$.

Now we explain our 
approach to Theorem \ref{thm:maintheorem}.  Let us start with defining 
$\mathcal{PE}_{X,\gamma}(F_d)$ and \eqref{eq:04.2.5kk20}. In fact, $\mathcal{T}(Z,F_d)$ can be rewritten as 
the derivative at $0$ of the Mellin transform of
\begin{equation}
	\mathrm{Tr_{s}}[(N^{\Lambda^{\bullet}(T^{\ast}Z)}-\frac{m}{2})\exp(-t\mathbf{D}^{Z,F_{d},2}/2)],\; t>0, 
	\label{eq:1.0.24parisconf}
\end{equation}
where $\mathrm{Tr_{s}}[\cdot]$ denotes the supertrace with respect to the 
$\Z_{2}$-grading on $\Lambda^{\bullet}(T^{\ast}Z)$. 

If 
$\gamma\in G$ is semisimple, let $\mathcal{E}_{X,\gamma}(F_{d},t)$ 
denote the orbital integral (see Subsection \ref{section3.3}) of the Schwartz kernel of 
$(N^{\Lambda^{\bullet}(T^{*}X)}-\frac{m}{2})\exp(-t\mathbf{D}^{X,F_{d},2}/2)$ associated with $\gamma$. Note that in $\mathcal{E}_{X,\gamma}(F_{d},t)$, we take the supertrace of the endomorphism on $\Lambda^{\bullet}(T^{*}X)\otimes F$ (see \eqref{eq:6.2.8pl}). Moreover, $\mathcal{E}_{X,\gamma}(F_{d},t)$ depends only on the conjugacy class of $\gamma$ in $G$. Let $\mathcal{ME}_{X,\gamma}(F_d,s)$ denote the Mellin 
transform of $\mathcal{E}_{X,\gamma}(F_d,t)$, $t>0$ with appropriate 
$s\in \C$. If $\gamma=1$, they are just $\mathcal{I}_{X}(F_d,t)$, 
$\mathcal{MI}_{X}(F_d,s)$ introduced in \eqref{eq:04.2.4ppaa} - 
\eqref{eq:04.2.6kk20}.

We use the notation in Subsection \ref{section3.5}. Let $[\Gamma]$ denote the set of the conjugacy classes in $\Gamma$. By applying the Selberg's trace formula to $Z=\Gamma\backslash X$, we get
\begin{equation}
	\begin{split}
		\mathrm{Tr_{s}}[(N^{\Lambda^{\bullet}(T^{\ast}Z)}-\frac{m}{2})\exp(-t\mathbf{D}^{Z,F_{d},2}/2)]=\sum_{[\gamma]\in[\Gamma]}\frac{\mathrm{Vol}(\Gamma\cap 
		Z(\gamma)\backslash 
		X(\gamma))}{|S(\gamma)|}\mathcal{E}_{X,\gamma}(F_{d},t).
	\end{split}
	\label{eq:03.5.14ksd} 
\end{equation}
Now we compare \eqref{eq:04.2.5kk20} with \eqref{eq:03.5.14ksd}. Then 
a proof to \eqref{eq:04.2.5kk20} mainly includes the 
following three parts:

\begin{enumerate}[label=\arabic*.]
	\item We show that if $[\gamma]\in E[\Gamma]$, then 
	$\mathcal{ME}_{X,\gamma}(F_d,s)$ admits a meromorphic extension 
	to $s\in \C$ which is holomorphic at $s=0$. Thus we define
	\begin{equation}
		\mathcal{PE}_{X,\gamma}(F_d)=\frac{\partial}{\partial s}|_{s=0} 
		\mathcal{ME}_{X,\gamma}(F_d,s).
		\label{eq:01.0.25kk20}
	\end{equation}
	Such consideration also holds for arbitrary 
	elliptic element $\gamma\in G$.
	\item If $\gamma\in \Gamma$ is elliptic, then it is of finite 
	order, from \eqref{eq:1.0.22parisconf}, we get that
	$\mathcal{PE}_{X,\gamma}(F_d)$ is a nice exponential polynomial in $d$ for 
	$d$ large enough. 
	
	\item We prove that all the terms in the sum of 
	\eqref{eq:03.5.14ksd} associated with nonelliptic $[\gamma]\in 
	[\Gamma]$ contribute as $\mathcal{O}(e^{-cd})$ in $\mathcal{T}(Z,F_d)$. 
\end{enumerate}

Indeed, to handle the contribution of the nonelliptic $[\gamma]\in 
[\Gamma]$, we use a spectral gap of $\mathbf{D}^{Z,F_{d},2}$ due to 
the nondegeneracy of $\lambda$. By \cite[Th\'{e}or\`{e}me 
3.2]{MR2838248}, \cite[Theorem 4.4]{BMZ2015toeplitz} which holds for 
a more general setting (see also  
\cite[Proposition 7.5, Corollary 7.6]{MR3128980} for a proof by using representation 
theory for symmetric spaces), there exist constants $C>0$, 
$c>0$ such that for $d\in\bN$, 
\begin{equation}
	\mathbf{D}^{Z,F_{d},2}\geq cd^{2}-C.
	\label{eq:1.0.27kkss20}
\end{equation}
That also explains \eqref{eq:1.00cohom} for large $d$. The Part 3 follows essentially from the same arguments as in 
\cite[Section 8]{MR3128980} and  \cite[Subsections 6.6, 7.2, Remarks 
7.8, 8.15]{BMZ2015toeplitz} which 
makes good use of \eqref{eq:1.0.27kkss20} and the fact that nonelliptic elements in $\Gamma$ 
admit a uniform strictly positive lower bound for their displacement distances on $X$.

For elliptic $\gamma\in \Gamma$, we apply Bismut's formula 
\cite[Theorem 6.1.1]{bismut2011hypoelliptic} to evaluate
$\mathcal{E}_{X,\gamma}(F_d,t)$. Then we can write $\mathcal{E}_{X,\gamma}(F_d,t)$ as an 
Gaussian-like integral with the integrand given 
as a product of an analytic function determined by the 
adjoint action of $\gamma$ on Lie algebras and the character $\chi_{E_{d}}$ of the representation $E_{d}$. By 
coordinating these two factors, especially using 
all sorts of character formulae for $\chi_{E_{d}}$, we can integrate 
it out. We show that 
$\mathcal{E}_{X,\gamma}(F_d,t)$ is a finite sum of the terms as 
follows,
\begin{equation}
	t^{-j-\frac{1}{2}}e^{-t(cd+b)^{2}}Q(d),
	\label{eq:1.0.28kkss}
\end{equation}
where $j\in\bN$, $c\neq 0$, $b$ are real constants, and $Q(d)$ is a 
nice exponential polynomial in $d$. It is crucial that $c\neq 0$. Indeed, we will see in Subsection 
\ref{section4.3} that this quantity $c$ measures the difference between 
the representations $(E_{\lambda},\rho^{E_{\lambda}})$ and 
$(E_{\lambda},\rho^{E_{\lambda}}\circ \theta)$. 

As a consequence of \eqref{eq:1.0.28kkss}, $\mathcal{PE}_{X,\gamma}(F_d)$ in \eqref{eq:01.0.25kk20}
is well-defined, which is clearly a nice 
exponential polynomial in $d$ (for $d$ large enough). The details 
on these computations are carried out in Subsection \ref{sec7.2bath}, where we apply the techniques inspired by the 
computations in Shen's approach to the Fried conjecture \cite[Section 
7]{Shen_2016} and also in its extension to orbifold case by Shen and 
Yu \cite{2017arXiv170408369S}.

The formula \eqref{eq:1.0.22parisconf} gives a new and geometric approach 
to the above results on $\mathcal{PE}_{X,\gamma}(F_d)$. It is nicer 
in the sense that each $\mathcal{PI}_{X(k)}(F^{k}_{\sigma,d})$ is 
already well 
understood and related to the $L_{2}$-torsions for the singular 
stratum of $Z$. For proving it, we apply a 
geometric localization formula for $\mathcal{E}_{X,\gamma}(F_d,t)$ as 
follows.
\begin{theorem}\label{thm:local0000}
	Assume that $\delta(G)=1$. We use the same notation as in Theorem 
	\ref{thm:maintheorem}. Let $\gamma=k\in T$. Then for $t>0$, $d\in\bN$, 
	\begin{equation}
		\mathcal{E}_{X,\gamma}(F_d,t)=\sum_{\sigma\in W^{1}_{U}(k)} 
		\varphi^{U}_{k}(\sigma,d\lambda+\lambda_{0}) 
		\mathcal{I}_{X(k)}(F^{k}_{\sigma,d},t).
		\label{eq:1.0.29sskk}
	\end{equation}
\end{theorem}

After taking the Mellin transform on both sides of 
\eqref{eq:1.0.29sskk}, we get exactly \eqref{eq:1.0.22parisconf}.
In Theorem \ref{thm:6.2.1ss}, we will show a general version of the 
above geometric 
localization formula for $\mathcal{E}_{X,\gamma}(F_d,t)$ associated with any 
semisimple element $\gamma\in G$.

Our approach to Theorem \ref{thm:local0000} is a more delicate application of Bismut's 
formula \cite[Theorem 
6.1.1]{bismut2011hypoelliptic}.
As we said,
$\mathcal{E}_{X,\gamma}(F_d,t)$, 
$\mathcal{I}_{X(k)}(F^{k}_{\sigma,d},t)$ are equal to
integrals of some integrands involving $\chi_{E_{d}}$, 
$\chi_{E^{k}_{\sigma,d}}$ respectively.  To relate the both sides of 
\eqref{eq:1.0.29sskk}, we employ a generalized version of Kirillov 
character formula (see Theorem \ref{thm:5.4.4ss20}) which gives an explicit way of decomposing $\chi_{E_{d}}|_{U(k)^{0}}$ into a sum of 
$\chi_{E^{k}_{\sigma,d}}$, $\sigma\in W^{1}_{U}(k)$. This character 
formula was proved by Duflo, 
Heckman and Vergne in \cite[II. 3, Theorem (7)]{DufloHeckmanVergne1984} under 
a general setting, and 
we will recall its special case for our need in Subsection 
\ref{section5.4}. Then we expand the integral formula for 
$\mathcal{E}_{X,\gamma}(F_d,t)$ carefully into a sum of certain 
integrals involving $\chi_{E^{k}_{\sigma,d}}$, $\sigma\in 
W^{1}_{U}(k)$, which correspond respectively
$\mathcal{I}_{X(k)}(F^{k}_{\sigma,d},t)$ via Bismut's 
formula. This way, we prove \eqref{eq:1.0.29sskk}.

Theorem \ref{thm:local0000} can be interpreted as follows, the action 
of
elliptic element $\gamma$ on $X$ could lead to a geometric localization 
onto its 
fixed point set $X(k)$ when we evaluate the orbital integrals. Even though we only prove it for a very restrictive situation, we 
still expect such phenomenon in general due to a 
geometric formulation for the semisimple orbital integrals (see \cite[Chapter 
4]{bismut2011hypoelliptic}).

Finally, we note that in \cite[Section 8]{BMZ2015toeplitz}, the authors 
explained well how to use 
Bismut's formula for semisimple orbital integrals to study the 
asymptotic analytic torsion. Here, we go one step further in that 
direction to get a refined evaluation on it. Bergeron and 
Venkatesh \cite{BV2013torsion} also studied the asymptotic analytic 
torsion but under a totally different 
setting. In \cite{liu:tel-01841334,LIU2021109117}, the asymptotic equivariant analytic torsion 
for a locally symmetric space was studied, and the oscillating terms 
also appeared naturally in that case. Moreover, Finski \cite[Theorem 1.5]{FINSKI20183457} obtained the full asymptotic expansion of the holomorphic torsions for the tensor powers of a given positive line bundle over a compact 
Hermitian orbifold.

This paper is organized as follows.
In Section \ref{section2paris}, we recall the definition of 
Ray-Singer analytic torsion for compact orbifolds. We also include a 
brief introduction to the orbifolds at beginning.

In Section \ref{section2}, we introduce the explicit geometric formula of Bismut 
for semisimple orbital integrals and the Selberg's trace formula for 
compact locally symmetric orbifolds. They are the main tools to study 
the analytic torsions in this paper.

In Section \ref{section3.6}, we give a vanishing theorem for 
$\mathcal{T}(Z,F)$, so that we only need to focus on the case 
$\delta(G)=1$.

In Section \ref{section4}, we study the Lie algebra of $G$ provided $
\delta(G)=1$. Furthermore, we introduce a generalized Kirillov 
formula for 
compact Lie groups.

In Section \ref{section6paris}, we prove a general version of Theorem  \ref{thm:local0000}.

In Section \ref{section7paris}, given the sequence 
$\{F_{d}\}_{d\in\bN}$, we compute explicitly 
$\cE_{X,\gamma}(F_{d},t)$ in terms of root systems for elliptic 
$\gamma$, in particular, we prove \eqref{eq:1.0.28kkss}. Then we 
give the formulae for $\cP\cI_{X}(F_{d})$, $\cP\cE_{X,\gamma}(F_{d})$. 

Finally, in Section \ref{section8paris}, we introduce the spectral gap 
\eqref{eq:1.0.27kkss20} and we give a proof to Theorem 
\ref{thm:maintheorem}.

In this paper, if $V$ is a real vector spaces and if $E$ is a complex 
vector space, we will use the symbol $V\otimes E$ to denote the 
complex vector space $V\otimes_\R E$. If both $V$ and $E$ are complex 
vector spaces, then $V\otimes E$ is just the usual tensor over $\C$.

\noindent\textbf{Acknowledgments.}
I would like to thank Prof. Jean-Michel Bismut, and Prof. Werner M\"{u}ller for encouraging me to work on this subject, and for many useful 
discussions. I also thank Dr. Taiwang Deng for educating me about the cohomology of
arithmetic groups, and Dr. Ksenia Fedosova for explanations of her 
results on the hyperbolic case. 

This work is carried out during my stay 
in Max Planck Institute for Mathematics (MPIM) in Bonn. I also want to 
express my sincere gratitude to MPIM for providing so nice research 
environment. Last but not least, I also would like to thank the 
anonymous referee(s) for the suggestions and comments that greatly 
improved the paper.

%%%%%%%%%%%%%%%%%%%%%%%%%%%%%%%%%%%%%%%%%%%%%%%%%%%%%%%%%%%%%%%%%%%%%%%%
\section{Ray-Singer analytic torsion}\label{section2paris}
In this section, we recall the definitions of the orbifold and the orbifold 
vector bundle. We also refer to \cite{MR0079769,SATAKE_1957} and \cite[Chapter 
1]{Adem_2007} for more details. Then we recall the definition of Ray-Singer analytic 
torsion for compact orbifolds, where we refer to \cite{MR2140438,2017arXiv170408369S} for more details. In particular, Shen and 
Yu in \cite{2017arXiv170408369S} extended many important results on 
real analytic torsion from manifold setting to orbifold setting.

%%%%%%%%%%%%%%%%%%%%%%%%%%%%%%%%%%%%%%%%%%%%%%%%%%%%%%%%%%%%%%%%%%%%%%%%
\subsection{Orbifolds and orbifold vector bundles}\label{subs:2.1}
Let $Z$ be a topological space.
\begin{definition}\label{def:2.1.1vogel}
	If $U$ is a connected open subset of $Z$, an orbifold chart for $U$ is a triple $(\widetilde{U},\pi_U,G_U)$ such that 
	\begin{itemize}
		\item $\widetilde{U}$ is a connected open set of some $\R^m$, $G_U$ is a finite group acting smoothly and effectively on $\widetilde{U}$ on the left;
		\item $\pi_U$ is a continuous surjective 
		$\widetilde{U}\rightarrow U$, which is invariant by 
		$G_U$-action;
		\item $\pi_U$ induces a homeomorphism between $G_U\backslash \widetilde{U}$ and $U$.
	\end{itemize}
	
	If $V\subset U$ is a connected open subset, an embedding of orbifold chart for the inclusion $i: V\rightarrow U$ is an orbifold chart $(\widetilde{V}, \pi_V, G_V)$ for $V$ and an orbifold chart $(\widetilde{U},\pi_U,G_U)$ for $U$ together with a smooth embedding $\phi_{UV}: \widetilde{V}\rightarrow \widetilde{U}$ such that the following diagram commutes,
	\begin{equation}
		\begin{gathered}  \xymatrix{
			\widetilde{V}\ar[r]^{\phi_{UV}}\ar[d]_{\pi_V}& \widetilde{U}\ar[d]^{\pi_U} \\
			V\ar[r]^{i}& U
			}\end{gathered}.
		\end{equation}
		
		If $U_1$, $U_2$ are two connected open subsets of $Z$ with 
		the charts $(\widetilde{U}_1,\pi_{U_1},G_{U_1})$, $(\widetilde{U}_2,\pi_{U_2},G_{U_2})$ respectively. We say that these two orbifold charts are compatible if for any point $z\in U_1\cap U_2$, there exists an open connected neighborhood $V\subset U_1\cap U_2$ of $z$ with an orbifold chart $(\widetilde{V}, \pi_V, G_V)$ such that there exist two embeddings of orbifold charts $\phi_{U_1V}:(\widetilde{V}, \pi_V, G_V)\rightarrow (\widetilde{U}_{1}, \pi_{U_1}, G_{U_1})$, $\phi_{U_2V}:(\widetilde{V}, \pi_V, G_V)\rightarrow (\widetilde{U}_{2},\pi_{U_2}, G_{U_2})$. In this case, the diffeomorphism $\phi_{U_2 V}\circ \phi_{U_1V}^{-1}: \phi_{U_1 V}(\widetilde{V})\rightarrow \phi_{U_2 V}(\widetilde{V})$ is called a coordinate transformation.
	\end{definition}
	
	\begin{definition}
		An orbifold atlas on $Z$ is couple $(\mathcal{U},\widetilde{\mathcal{U}} )$ consisting of a cover $\mathcal{U}$ of open connected subsets of $Z$ and a family of compatible orbifold charts $\widetilde{\mathcal{U}}=\{(\widetilde{U},\pi_U, G_U)\}_{U\in\mathcal{U}}$.
		
		An orbifold atlas $(\mathcal{V},\widetilde{\mathcal{V}})$ is called a refinement of $(\mathcal{U},\widetilde{\mathcal{U}})$ if $\mathcal{V}$ is a refinement of $\mathcal{U}$ and if every orbifold chart in $\widetilde{\mathcal{V}}$ has an embedding into some orbifold chart in $\widetilde{\mathcal{U}}$. Two orbifold atlas are said to be equivalent if they have a common refinement, and the equivalent class of an orbifold atlas is called an orbifold structure on $Z$.
		
		An orbifold is a second countable Hausdorff space equipped with an orbifold structure. It is said to have dimension $m$ if all the orbifold charts which defines the orbifold structure are of dimension $m$.
	\end{definition}

	If $Z, Y$ are two orbifolds, a smooth map $f: Z\rightarrow Y$ is 
	a continuous map from $Z$ to $Y$ such that it lifts locally to an 
	equivariant smooth map from an orbifold chart of $Z$ to any orbifold chart of $Y$. In this way, we can define the notion of smooth functions and the smooth action of Lie groups.
	
	By \cite[Proposition 2.12]{2017arXiv170408369S}, if
	$\Gamma$ is discrete group acting smoothly and properly 
	discontinuously on the left on an orbifold $X$, then $Z=\Gamma\backslash X$ has a canonical orbifold structure induced from $X$.

	In the sequel, let $Z$ be an orbifold with an orbifold structure given by $(\mathcal{U},\widetilde{\mathcal{U}})$. If $z\in Z$, there exists an open connected neighborhood $U_z$ of $z$ with a compatible orbifold chart $(\widetilde{U}_z, G_z, \pi_z)$ such that $\pi_z^{-1}(z)$ contains only one point $x\in \widetilde{U}_z$. Then $G_z$ does not depend on the choice of such open connected neighborhood (up to canonical isomorphisms compatible with the orbifold structure), then $G_z$ is called the local group at $z$.
	
	Put
	\begin{equation}
		Z_{\mathrm{reg}}=\{z\in Z\,:\, G_z=\{1\}\},\;\; Z_{\mathrm{sing}}=\{z\in Z\,:\, G_z\neq\{1\}\}.
	\end{equation}
	Then $Z_{\mathrm{reg}}$ is naturally a smooth manifold. But 
	$Z_{\mathrm{sing}}$ is not necessarily an orbifold. In \cite[Section 
	2]{MR0474432}, Kawasaki explained two different methods to view 
	$Z_{\mathrm{sing}}$ as an immersed image of a disjoint union of 
	orbifolds. We just recall that method which appears naturally in Kawasaki's 
	local index theorems for orbifolds \cite{MR0474432,MR527023}.
	
	If $z\in Z_{\mathrm{sing}}$, let $1=(h_{z}^{0}), (h_{z}^{1}), \cdots, 
	(h^{l_{z}}_{z})$ be the conjugacy classes in $G_{z}$. Put
	\begin{equation}
		\Sigma Z=\{(z,(h^{j}_{z}))\;|\; z\in Z_{\mathrm{sing}}, 
		j=1,\cdots, l_{z}\}.
		\label{eq:1.1.3}
	\end{equation}
	Let $(\widetilde{U}_z, G_z, \pi_z)$ be the local orbifold chart for 
	$z\in Z_{\mathrm{sing}}$ such that $\pi_z^{-1}(z)$ contains only one 
	point. For $j=1,\cdots, l_{z}$, let 
	$\widetilde{U}_z^{h^{j}_{z}}\subset \widetilde{U}_z$ be the fixed 
	point set of $h^{j}_{z}$, which is a submanifold of 
	$\widetilde{U}_z$. Note that $\widetilde{U}_z^{h^{j}_{z}}\subset 
	Z_{\mathrm{sing}}$. Let $Z_{G_{z}}(h^{j}_{z})$ be the centralizer of 
	$h^{j}_{z}$ in $G_z$. Then $Z_{G_{z}}(h^{j}_{z})$ acts smoothly on 
	$\widetilde{U}_z^{h^{j}_{z}}$. Put
	\begin{equation}
		K^{j}_{z}=\ker (Z_{G_{z}}(h^{j}_{z})\rightarrow 
		\mathrm{Aut}(\widetilde{U}_z^{h^{j}_{z}})).
		\label{eq:1.1.4}
	\end{equation}
	Then $(\widetilde{U}_z^{h^{j}_{z}},Z_{G_{z}}(h^{j}_{z})/K^{j}_{z}, 
	\pi^{j}_{z}: \widetilde{U}_z^{h^{j}_{z}}\rightarrow 
	\widetilde{U}_z^{h^{j}_{z}}/Z_{G_{z}}(h^{j}_{z}))$ defines an 
	orbifold chart near $(z, (h^{j}_{z}))\in \Sigma Z$. They form an 
	orbifold structure for $\Sigma Z$. Let $Z^{i}, i=1,\cdots, l$ denote 
	the connected component of the orbifold $\Sigma Z$.
	
	The integer $m^{j}_{z}=|K^{j}_{z}|$ is called the 
	multiplicity of $\Sigma Z$ in $Z$ at $(z, (h^{j}_{z}))$. This defines 
	a function $m: \Sigma Z\rightarrow \mathbb{Z}_{+}$. As explained in 
	\cite[Section 1]{MR0474432}, $m$ is locally constant on $\Sigma Z$, 
	and let $m_{i}\in \mathbb{Z}_{+}$ be the value of $m$ on $Z^{i}$ for $i=1,\cdots,l$. We call $m_{i}$ 
	the multiplicity of $Z^{i}$ in $Z$. We will put
	\begin{equation}
		Z^{0}=Z,\; m_{0}=1.
		\label{eq:1.1.5}
	\end{equation}
	
	\begin{remark}\label{rk:2.1.3vogel}
		In Definition \ref{def:2.1.1vogel}, for an orbifold chart, we 
		require the action $G_{U}$ on $\widetilde{U}$ to be 
		effective. To emphasize this condition, the orbifold defined 
		above is often
		called an effective orbifold. In fact, we can drop this 
		effectiveness, then we get a general version of 
		the (possibly ineffective) orbifold, for example, using the 
		orbifold groupoid (see \cite[Definition 1.38]{Adem_2007}). 
		The point-view of orbifold groupoid provides a 
		unified way to deal with effective and ineffective orbifolds.
		
		As explained in \cite[Example 2.5]{Adem_2007}, for global quotient 
		groupoids (including all the effective orbifolds and certain 
		ineffective orbifolds), a natural stratification called 
		the inertia groupoid was
		introduced as an extension of the one $\cup_{i=0}^{l} Z^{i}$ 
		defined in \eqref{eq:1.1.3} - \eqref{eq:1.1.5}. It 
		plays a key role in the study of the geometry of orbifolds. 
		We will go back to this point in Subsections \ref{section3.4} 
		\& \ref{section3.5}. Through this paper, the terminology 
		orbifold will always refer to the effective one unless 
		otherwise stated.
	\end{remark}

	We say $E$ to be an orbifold vector bundle of rank $r$ on $Z$ if there exists a 
	smooth map of orbifolds $\pi: E\rightarrow Z$ such that for any $U\in 
	\mathcal{U}$ and $(\widetilde{U}, G_{U}, \pi_{U})\in 
	\widetilde{\mathcal{U}}$, there exists an orbifold chart 
	$(\widetilde{U}^{E}, G^{E}_{U}, \pi^{E}_{U})$ of $E$ such that 
	$\widetilde{U}^{E}$ is an vector bundle on $\widetilde{U}$ of rank 
	$r$ equipped an effective action of $G^{E}_{U}$ and 
	$\pi^{E}_{U}(\widetilde{U}^{E})=\pi^{-1}(U)$. Moreover, there exists a 
	surjective group morphism $\psi_{U}: G^{E}_{U}\rightarrow G_{U}$ such 
	that the action of $G^{E}_{U}$ on $\widetilde{U}$ is identified via $\psi_{U}$ with the 
	action of $G_{U}$ on $\widetilde{U}$. If we have an open embedding 
	$\phi_{UV}:(\widetilde{V}, \pi_V, G_V)\rightarrow (\widetilde{U}, 
	\pi_{U}, G_{U})$, we require that it lifts to 
	the open embedding
	$\phi^{E}_{UV}:(\widetilde{V}^{E}, \pi^{E}_V, G^{E}_V)\rightarrow 
	(\widetilde{U}^{E}, 
	\pi_{U}^{E}, G_{U}^{E})$ of the orbifold charts of $E$ such that 
	$\phi^{E}_{UV}:\widetilde{V}^{E} \rightarrow \widetilde{U}^{E}$ is a 
	morphism of vector bundles associated with the open embedding 
	$\phi_{UV}:\widetilde{V} \rightarrow \widetilde{U}$. If every 
	$\psi_{U}: G^{E}_{U}\rightarrow G_{U}$ is an isomorphism of groups, 
	we call $E$ a proper orbifold vector bundle on $Z$. 
	
	Note that if $E$ is proper, then the rank of $E$ can be extended to a locally constant 
	function $\rho$ on $\Sigma Z$. The orbifold chart of $Z^{i}$ is given 
	by the triples such as $(\widetilde{U}_z^{h^{j}_{z}},Z_{G_{z}}(h^{j}_{z})/K^{j}_{z}, 
	\pi^{j}_{z}: \widetilde{U}_z^{h^{j}_{z}}\rightarrow 
	\widetilde{U}_z^{h^{j}_{z}}/Z_{G_{z}}(h^{j}_{z}))$. By the above 
	definition of $E$, we have an orbifold chart $(\widetilde{U}^{E}, 
	G^{E}_{U}=G_{U}, \pi^{E}_{U})$ such that $\widetilde{U}^{E}$ is a 
	$G_{U}$-equivariant vector bundle on $\widetilde{U}$. Then for $x\in 
	\widetilde{U}_z^{h^{j}_{z}}$, $h^{j}_{z}$ acts on the fibres 
	$\widetilde{U}^{E}_{z}$ linearly, so that we can set
	$\rho(z,(h^{j}_{z}))=\mathrm{Tr}^{\widetilde{U}^{E}_{z}}[h^{j}_{z}]$. 
	Then $\rho$ is really a locally constant function 
	on $\Sigma Z$. For $i=1,\cdots, l$, let $\rho_{i}$ be the value of 
	$\rho$ on the component $Z^{i}$. We also put $\rho_{0}=r$.

	We call $s: Z\rightarrow E$ a smooth section of $E$ over $Z$ if 
	it is a smooth map between orbifolds such that $\pi\circ 
	s=\mathrm{Id}_{Z}$. We will use $C^{\infty}(Z,E)$ to denote the 
	vector space of smooth sections of $E$ over $Z$.
	
	Take an orbifold chart $(\widetilde{U}, G_{U}, \pi_{U})\in 
	\widetilde{\mathcal{U}}$ of $Z$, $G_{U}$ acts canonically on the 
	tangent vector bundle $T\widetilde{U}$ of $\widetilde{U}$. The open 
	embeddings of orbifold charts of $Z$ also lift to the open embeddings 
	of their tangent vector bundles. This way, we get a proper orbifold 
	vector bundle $TZ$ on $Z$, and the projection $\pi: TZ\rightarrow Z$ 
	is just given by the obvious projection $T\widetilde{U}\rightarrow 
	\widetilde{U}$. We call $TZ$ the tangent vector bundle of 
	$Z$. If we equipped $TZ$ with Euclidean metric $g^{TZ}$, we will 
	call $Z$ a Riemannian 
	orbifold and call
	$g^{TZ}$ a Riemannian metric of $Z$.

	Let $\Omega^{\bullet}(Z)$ denote the set of smooth differential forms 
	of $Z$, which has a $\Z$-graded structure by degrees. The 
	de Rham differential $d^{Z}: \Omega^{\bullet}(Z)\rightarrow 
	\Omega^{\bullet+1}(Z)$ is given by the family of de Rham 
	differential operators 
	$d^{\widetilde{U}}: \Omega^{\bullet}(\widetilde{U})\rightarrow 
	\Omega^{\bullet+1}(\widetilde{U})$. Then we can define the de Rham 
	complex $(\Omega^{\bullet}(Z), d^{Z})$ of $Z$ and the associated de 
	Rham cohomology $H^{\bullet}(Z,\R)$. By \cite[Section 
	1]{MR0474432}, there is a natural isomorphism between 
	$H^{\bullet}(Z,\R)$ and the singular cohomology of the underlying 
	topological space $Z$.

	Now let us recall the integrals on $Z$. Assume that $Z$ is 
	compact. We may take a finite open covering 
	$\{U_{i}\}_{i\in I}$ of the precompact orbifold charts for $Z$. Since $Z$ is Hausdorff, then there exists a partition of unity 
	subordinate to this open cover. We can find a family of smooth 
	functions $\{\phi_{i}\in C^{\infty}_{c}(Z)\}_{i\in I}$ with values in 
	$[0,1]$ such that $\mathrm{Supp}(\phi_{i})\subset U_{i}$, and that \begin{equation}
	\sum_{i\in I} \phi_{i}=1.
	\label{eq:1.1.6alpha}
\end{equation}
Take $\widetilde{\phi}_{i}=\pi^{*}_{U_{i}}(\phi_{i})\in 
C^{\infty}_{c}(\widetilde{U}_{i})^{G_{U_{i}}}$.

If $\alpha\in \Omega^{m}(Z,o(TZ))$, let 
$\widetilde{\alpha}_{U_{i}}$ be its lift on the chart 
$(\widetilde{U}_{i},\pi_{U_{i}},G_{U_{i}})$. We define
\begin{equation}
	\int_{Z} \alpha = \sum_{i}\frac{1}{|G_{U_{i}}|} 
	\int_{\widetilde{U}_{i}} 
	\widetilde{\phi}_{i}\widetilde{\alpha}_{U_{i}}.
	\label{eq:1.1.7alpha}
\end{equation}
By \cite[Section 3.2]{2017arXiv170408369S}, if $\alpha\in 
\Omega^{m}(Z,o(TZ))$, then $\alpha$ is also integrable on 
$Z_{\mathrm{reg}}$, so that
\begin{equation}
	\int_{Z} \alpha = \int_{Z_{\mathrm{reg}}}\alpha.
	\label{eq:1.1.8alpha}
\end{equation}
Also if $\alpha\in 
\Omega^{\bullet}(Z,o(TZ))$, we have
\begin{equation}
	\int_{Z} d^{Z}\alpha=0.
	\label{eq:1.19a}
\end{equation}

If $(Z,g^{TZ})$ is a Riemannian orbifold, we can define the 
integration of functions on $Z$ with respect to the Riemannian volume element. 
If we have a Hermitian orbifold vector bundle $(F,h^{F})\rightarrow 
(Z,g^{TZ})$, one can define the $L_{2}$ scalar product for the space 
of continuous sections of $F$ as usual. Then, after completion, we get the Hilbert space 
$L^{2}(Z,F)$. 

The Chern-Weil theory on the characteristic forms extends to 
orbifolds, where their constructions are parallel to 
the case of smooth manifolds. We 
refer to \cite[Subsection 3.4]{2017arXiv170408369S} for more details. Note 
that the 
characteristic forms are not only defined on $Z$ but also defined on 
$\Sigma Z$. The part $\Sigma Z$ has a nontrivial contribution in Kawasaki's 
local index theorems for orbifolds \cite{MR0474432,MR527023}.

Finally, we introduce the orbifold Euler characteristic number of $(Z, 
g^{TZ})$ \cite{SATAKE_1957}. Let 
$\nabla^{TZ}=\{\nabla^{T\widetilde{U}_{i}}\}_{U_{i}\in \mathcal{U}}$ be the Levi-Civita 
connection on $TZ$ associated with $g^{TZ}$. The Euler form $e(TZ,\nabla^{TZ})\in 
\Omega^{m}(Z,o(TZ))$ is given by the family of closed forms
\begin{equation}
	\big\{e(\widetilde{U}_{i},\nabla^{T\widetilde{U}_{i}})\in 
	\Omega^{m}(\widetilde{U}_{i},o(T\widetilde{U}_{i}))^{G_{U_i}}\big\}_{U_{i}\in \mathcal{U}}.
	\label{eq:2.1.10ssd}
\end{equation}
If $Z$ is oriented, then we can view $e(TZ,\nabla^{TZ})$ as a 
differential form on $Z$.

If $Z$ is compact, set
\begin{equation}
	\chi_{\mathrm{orb}}(Z)=\int_{Z}e(TZ,\nabla^{TZ}).
	\label{eq:2.1.11ssdd}
\end{equation}
By \cite[Section 3]{SATAKE_1957}, $\chi_{\mathrm{orb}}(Z)$ is a 
rational number, and it vanishes when $Z$ is odd dimensional.
%%%%%%%%%%%%%%%%%%%%%%%%%%%%%%%%%%%%%%%%%%%%%%%%%%%%%%%%%%%%%%%%%%%%%%%%%
\subsection{Flat vector bundles and analytic torsions of 
orbifolds}\label{subs:2.2}

If $(F,\nabla^{F})$ is an orbifold vector bundle over $Z$ 
with a connection $\nabla^{F}$, we call $(F,\nabla^{F})$ a flat 
vector bundle if the curvature $R^{F}=\nabla^{F,2}$ vanishes 
identically on $Z$. A detailed discussion for the flat vector bundles 
on $Z$ is given in \cite[Sections 2.3 - 2.5]{2017arXiv170408369S}.

Let $(Z,g^{TZ})$ be a compact Riemannian orbifold of dimension $m$. 
Let $(F,\nabla^F)$ be a flat complex orbifold vector bundle of rank 
$r$ on $Z$ with Hermitian metric $h^F$. Note that we do not assume 
that $F$ is proper.

Let $\Omega^\bullet(Z,F)$ be the set of smooth sections of 
$\Lambda^\bullet(T^*Z)\otimes F$ on $Z$. Let $d^Z$ be the exterior 
differential acting on $\Omega^\bullet(Z,\R)$.
\begin{definition} For $i=0,1,\cdots, m$, if $\alpha\in\Omega^i(Z,\R)$, $s\in C^\infty(Z,F)$, the operator $d^{Z,F}$ acting on $\Omega^i(Z,F)$ is defined by
	\begin{equation}
		d^{Z,F}(\alpha\otimes s)=(d^Z\alpha)\otimes s + (-1)^i \alpha\wedge \nabla^F s \in \Omega^{i+1}(Z,F).
		\label{eq:deRhamop}
	\end{equation}
\end{definition}

Since $\nabla^{F,2}=0$, then $(\Omega^\bullet(Z,F),d^{Z,F})$ is a 
complex, which is called the de Rham complex for the flat orbifold 
vector bundle $(F,\nabla^F)$ on $Z$. Let $H^\bullet(Z,F)$ denote the 
corresponding de Rham cohomology group of $Z$ valued in 
$F$, as in the case of closed manifolds, $H^\bullet(Z,F)$ is always 
finite dimensional.

Let $\langle\cdot,\cdot\rangle_{\Lambda^\bullet(T^*Z)\otimes F,z}$ be 
the Hermitian metric on $\Lambda^\bullet(T_z^*Z)\otimes F_z$, $z\in Z$ 
induced by $g^{TZ}_z$ and $h^F_z$. Let $dv$ be the Riemannian volume 
element on $Z$ induced by $g^{TZ}$. The $L_2$-scalar 
product on $\Omega^\bullet(Z,F)$ is given as follows, if $s, s'\in 
\Omega^\bullet(Z,F)$, then 
\begin{equation}\label{eq:L2metric}
	\langle s, s'\rangle_{L^2}=\int_{Z}\langle 
	s(z),s(z')\rangle_{\Lambda^\bullet(T^*Z)\otimes F,z} dv(z).
\end{equation}
By \eqref{eq:1.1.8alpha}, it will be the same if we take the 
integrals on $Z_{\mathrm{reg}}$.

Let $d^{Z,F,*}$ be the formal adjoint of $d^{Z,F}$ with respect to 
the above $L_{2}$-metric on $\Omega^\bullet(Z,F)$, i.e., for $s, s' \in 
\Omega^\bullet(Z,F)$, 
\begin{equation}\label{eq:formaladj}
	\langle d^{Z,F,*}s,s'\rangle_{L^2}=\langle s, d^{Z,F}s'\rangle_{L^2}.
\end{equation}
Then $d^{Z,F,*}$ is a first-order differential operator acting 
$\Omega^\bullet(Z,F)$ on which decreases the degree by $1$.

\begin{definition}
	The de Rham - Hodge operator $\mathbf{D}^{Z,F}$ of $\Omega^\bullet(Z,F)$ is defined as
	\begin{equation}\label{eq:DZFop}
		\mathbf{D}^{Z,F} = d^{Z,F}+d^{Z,F,*}.
	\end{equation}
	It is a first-order self-adjoint elliptic differential operator 
	acting on $\Omega^\bullet(Z,F)$.
\end{definition}

The Hodge Laplacian is
\begin{equation}\label{eq:HodgeLap}
	\mathbf{D}^{F,Z,2}=[d^{Z,F}, d^{Z,F,*}]=d^{Z,F}d^{Z,F,*}+d^{Z,F,*}d^{Z,F}.
\end{equation}
Here, $[\cdot,\cdot]$ denotes the supercommutator. Then $\mathbf{D}^{Z,F,2}$ 
is a second-order essentially self-adjoint non-negative elliptic 
operator, which preserves the degree.

The Hodge decomposition for $\Omega^{\bullet}(Z,F)$ still holds in this 
case (see \cite[Proposition 2.2]{MR2140438}, \cite[Proposition 
2.1]{MR3623749}),
\begin{equation}
	\Omega^{\bullet}(Z,F)= \ker (\mathbf{D}^{Z,F,2}|_{\Omega^{\bullet}(Z,F)}) \oplus 
	\mathrm{Im}(d^{Z,F}|_{\Omega^{\bullet-1}(Z,F)})\oplus 
	\mathrm{Im}(d^{Z,F,*}|_{\Omega^{\bullet+1}(Z,F)}).
	\label{eq:1.2.6ab}
\end{equation}
Then we have the canonical identification of vector spaces,
\begin{equation}
	\mathcal{H}^\bullet(Z,F):=\ker \mathbf{D}^{Z,F,2} \simeq H^{\bullet}(Z,F).
	\label{eq:1.2.7alpha}
\end{equation}

Put
\begin{equation}
	\chi(Z,F)=\sum_{j=0}^m (-1)^j 
	\dim H^j(Z,F).
	\label{eq:7.8.8s}
\end{equation}

If $F$ is proper, recall that the numbers $\rho_{i}$, $i=0,\cdots, l$ 
are defined in previous subsection as the extension of the rank of 
$F$. Then by \cite[Theorem 4.3]{2017arXiv170408369S}, we have
\begin{equation}
	\chi(Z,F)=\sum_{i=0}^{l}\rho_{i}\frac{\chi_{\mathrm{orb}}(Z_{i})}{m_{i}}.
	\label{eq:1.2.12alpha}
\end{equation}
The right-hand side of \eqref{eq:1.2.12alpha} contains the nontrivial 
contributions from $\Sigma Z$.

Let $P$ denote the orthogonal projection from 
$\Omega^\bullet(Z,F)$ to $\mathcal{H}^\bullet(Z,F)$. Let 
$\mathcal{H}^\perp$ denote the orthogonal subspace of 
$\mathcal{H}^\bullet(Z,F)$ in $\Omega^\bullet(Z,F)$, and let 
$(\mathbf{D}^{Z,F,2})^{-1}$ be the inverse of $\mathbf{D}^{Z,F,2}$ acting on 
$\mathcal{H}^\perp$. Let $N^{\Lambda^\bullet(T^*Z)}$ be the number 
operator on $\Lambda^\bullet(T^*Z)$ which acts on $\Lambda^j(T^*Z)$ 
by multiplication of $j$.

For $s\in\C$, 
$\Re(s)$ is large enough, set
\begin{equation}
	\label{eq:2.2.12vogel}
	\begin{split}
		\vartheta(F)(s) &= -\mathrm{Tr_s}[N^{\Lambda^{\bullet}(T^*Z)}(\mathbf{D}^{Z,F,2})^{-s}]\\
		&= - \frac{1}{\Gamma(s)}\int_0^{+\infty} 
		\mathrm{Tr_s}[N^{\Lambda^{\bullet}(T^*Z)}\exp(-t\mathbf{D}^{Z,F,2})(1-P)]t^{s-1}dt,		
	\end{split}
\end{equation}
where $\Gamma(s)$ is the Gamma function for $s\in \C$.
By the short time
asymptotic expansions of the heat trace (see \cite[Proposition 
2.1]{MR2140438}), $\vartheta(F)(s)$ admits a unique meromorphic extension to $s\in \C$ which is holomorphic at $s=0$. 

\begin{definition}\label{def:2.2.5ss20}
	Let $\mathcal{T}(g^{TZ},\nabla^F,h^F)\in \R$ be given by
	\begin{equation}\label{eq:1.1.12bonn}
		\mathcal{T}(g^{TZ},\nabla^F,h^F)=\frac{d}{ds}|_{s=0} \vartheta(F)(s).
	\end{equation}
	The quantity $\mathcal{T}(g^{TZ},\nabla^F,h^F)$ is called Ray-Singer analytic torsion associated with $(F,\nabla^F,h^F)$.
\end{definition}

By \cite[Proposition 4.6, Corollary 4.9]{2017arXiv170408369S}, for an 
orientable closed orbifold $Z$, if $m$ is even and $F$ is unitarily 
flat, then 
$\mathcal{T}(g^{TZ},\nabla^F,h^F)=0$; if $m$ is odd and $F$ is acyclic, then $\mathcal{T}(g^{TZ},\nabla^F,h^F)$ is independent of the metrics $g^{TZ}$ and $h^F$.

Now we explain how to 
evaluate $\mathcal{T}(g^{TZ},\nabla^F,h^F)$ in practice when $F$ is acyclic. Using the analog 
arguments in \cite[Theorem 
7.10, Section XI]{MR1185803}, as $t\rightarrow 0^{+}$, the 
heat supertrace
$\mathrm{Tr_s}\left[\left(N^{\Lambda^{\bullet}(T^*Z)}-\frac{m}{2}\right)\exp(-t\mathbf{D}^{Z,F,2}/2)\right]$ either has a leading term as a multiple of $\frac{1}{\sqrt{t}}$ or is a small quantity as $\mathcal{O}(\sqrt{t})$ (see \cite[Eq. (4.37)]{2017arXiv170408369S}). To deal with this possible divergent term $\frac{1}{\sqrt{t}}$ in the integral of \eqref{eq:2.2.12vogel}, we proceed as in the proof of \cite[Theorem 3.29]{BismutLott1995}, for $t>0$, put
\begin{equation}
	b_t(g^{TZ},F)=(1+2t\frac{\partial}{\partial 
	t})\mathrm{Tr_s}\bigg[\big(N^{\Lambda^{\bullet}(T^*Z)}-\frac{m}{2}\big)\exp(-t\mathbf{D}^{Z,F,2}/2)\bigg].
	\label{eq:7.8.7kkkk}
\end{equation}
By \cite[Theorem 7.10]{MR1185803}, \cite[Theorem 2.13]{BismutLott1995} and \cite[Subsection 
4.3]{2017arXiv170408369S} and that $F$ is acyclic, as $t\rightarrow 0$,
\begin{equation}
	b_t(g^{TZ},F)=\mathcal{O}(\sqrt{t})\,;
	\label{eq:7.8.8kkkk}
\end{equation}
as $t\rightarrow +\infty$,
\begin{equation}
	b_t(g^{TZ},F)= \mathcal{O}(1/\sqrt{t}).
	\label{eq:7.8.9kkkk}
\end{equation}

By \cite[Theorem 3.29]{BismutLott1995} and \cite[Corollary 4.14]{2017arXiv170408369S}, we have
\begin{equation}
	\mathcal{T}(g^{TZ}, \nabla^{F},h^{F})=-\int_0^{+\infty} 
	b_t(g^{TZ},F)\frac{dt}{t}.
	\label{eq:7.8.11bonn}
\end{equation}
One particular case is that if for $t>0$, we always have 
\begin{equation}
	\mathrm{Tr_s}\left[\left(N^{\Lambda^{\bullet}(T^*Z)}-\frac{m}{2}\right)\exp(-t\mathbf{D}^{Z,F,2}/2)\right]=0,
\end{equation}
then $\mathcal{T}(g^{TZ},\nabla^F,h^F)=0$. This holds even for 
non-acyclic $F$.

%%%%%%%%%%%%%%%%%%%%%%%%%%%%%%%%%%%%%%%%%%%%%%%%%%%%%%%%%%%%%%%%%%%%%%%%%
\section{Orbital integrals and locally symmetric spaces}\label{section2}
In this section, we recall the geometry of the symmetric 
space $X$, and we recall an explicit geometric formula for semisimple orbital 
integrals obtained by Bismut 
\cite[Chapter 6]{bismut2011hypoelliptic} . Then, given a cocompact discrete subgroup $\Gamma\subset G$, we 
describe the orbifold structure on $Z=\Gamma\backslash 
X$, and we give the Selberg's trace formula for $Z$.

In this section, $G$ is taken to be a connected linear real reductive Lie 
group, we do not require that it has a compact center. Then $X$ is a 
symmetric space which may have de Rham components of both noncompact 
type and Euclidean type.
%%%%%%%%%%%%%%%%%%%%%%%%%%%%%%%%%%%%%%%%%%%%%%%%%%%%%%%%%%%%%%%%%%%%%%%%%
\subsection{Real reductive Lie group}\label{section3.1}
Let $G$ be a connected linear real reductive Lie group
with Lie algebra $\g$, and let $\theta\in \mathrm{Aut}(G)$ be a 
Cartan involution. Let $K$ be the fixed point set of $\theta$ in $G$. 
Then $K$ is a maximal compact subgroup of $G$, and let $\kk$ be its Lie algebra. Let $\pp\subset \g$ be the eigenspace of $\theta$ associated with the eigenvalue $-1$. The Cartan decomposition of $\g$ is given by
\begin{equation}
	\mathfrak{g}=\mathfrak{p}\oplus \mathfrak{k}.
	\label{eq:0.2.1ugc}
\end{equation}
Put $m=\dim \pp$, $n=\dim \kk$. 

Let $B$ be a $G$- and $\theta$-invariant nondegenerate symmetric bilinear form 
on $\g$, which is positive on $\pp$ and negative on $\kk$. It induces 
a symmetric bilinear form $B^{*}$ on $\g^{*}$, which extends to a 
symmetric bilinear form on $\Lambda^{\bullet}(\g^{*})$. The 
$K$-invariant bilinear form 
$\langle\cdot,\cdot\rangle=-B(\cdot,\theta\cdot)$ is a scalar 
product on $\g$, which extends to a scalar product on 
$\Lambda^{\bullet}(\g^{*})$. We will use $|\cdot|$ to denote the norm 
under this scalar product.

Let $U\g$ be the universal enveloping algebra of $\g$. Let $C^\g\in U\g$ be the Casimir element associated with $B$, i.e., if $\{e_i\}_{i=1,\cdots, m+n}$ is a basis of $\g$, and if $\{e^*_i\}_{i=1,\cdots, m+n}$ is the dual basis of $\g$ with respect to $B$, then
\begin{equation}
	C^\g= -\sum e^*_i e_i.
	\label{eq:3.1.2th19}
\end{equation}
We can identify $U\g$ with the algebra of left-invariant differential 
operators over $G$, then $C^\g$ is a second-order differential 
operator, which is $\mathrm{Ad}(G)$-invariant. Similarly, let 
$C^{\kk}\in U\kk$ denote the Casimir operator associated with 
$(\kk,B|_{\kk})$.

Let $\z_\g\subset \g$ be the center of $\g$. Put
\begin{equation}
	\g_{\mathrm{ss}}=[\g,\g].
	\label{eq:3.1.3th19}
\end{equation}
Then
\begin{equation}
	\g=\z_\g\oplus\g_{\mathrm{ss}}.
	\label{eq:3.1.4th19}
\end{equation}
They are orthogonal with respect to $B$.

Let $Z_{G}$ be the center of $G$, let $G_{\mathrm{ss}}$ be the 
closed analytic subgroup of $G$ associated with $\g_{\mathrm{ss}}$ 
(see \cite[Corollary 7.11]{knapp2002liegroupe}). Then $G$ 
is the commutative product of $Z_{G}$ and $G_{\mathrm{ss}}$, in particular,
\begin{equation}
	G=Z_{G}^0G_{\mathrm{ss}}.
	\label{eq:3.1.5th19}
\end{equation}

Let $i=\sqrt{-1}$ denote one square root of $-1$. Put
\begin{equation}
	\ku=\sqrt{-1} \pp\oplus \kk.
	\label{eq:3.1.6th19}
\end{equation}
For saving notation, if $a\in\pp$, we use notation $ia$ or $\ii a\in\sqrt{-1}\pp\subset\ku$ to denote the 
corresponding vector.

Then $\ku$ is a (real) Lie algebra, which is called the compact form of $\g$. Then 
\begin{equation}
	\g_\C=\ku_\C.
	\label{eq:3.1.7th19}
\end{equation}

Let $G_\C$ be the complexification of $G$ with Lie algebra $\g_\C$, 
which is closed and linear reductive \cite[Proposition 5.6]{knapp1986representation}. 
Then $G$ is the analytic subgroup of $G_\C$ with Lie algebra $\g$. 
Let $U\subset G_\C$ be the analytic subgroup associated with $\ku$. 
If $G$ has 
compact center, i.e., $\z_{\g}\cap\pp=\{0\}$, then by \cite[Proposition 5.3]{knapp1986representation}, $U$ is compact; 
since $G_{\C}$ is closed, $U$ is a maximal compact subgroup of $G_\C$.

\begin{definition}
	An element $\gamma\in G$ is said to be semisimple if there exists $g\in G$ such that 
	\begin{equation}\label{eq:3.1.8vogel}
		\gamma=g(e^ak)g^{-1}, a\in\pp, k\in K, \mathrm{Ad}(k)a=a.
	\end{equation}
	We call $\gamma_h=ge^ag^{-1}$, $\gamma_e=gkg^{-1}$ the 
	hyperbolic, elliptic parts of $\gamma$. These two parts are uniquely 
	determined by $\gamma$. If $\gamma_h=1$, we say $\gamma$ to be 
	elliptic, and if $\gamma_e=1$ and $\gamma_{h}\neq 1$, we say $\gamma$ to be hyperbolic.
\end{definition}

Let $Z(\gamma)$ be the centralizer of $\gamma$ in $G$. If $ v\in \g$, 
let $Z(v)\subset G$ be the stabilizer of $v$ in $G$ via the adjoint 
action. Let $\z(\gamma)$, $\z(v)$ be the Lie algebras 
of $Z(\gamma)$, $Z(v)$ respectively. If $\gamma=\gamma_{h}\gamma_{e}$ is semisimple as above, by \cite[Theorem 
2.19.23]{eberlein1996geometry} and \cite[Lemma 7.36]{knapp2002liegroupe}, 
\begin{equation}
	Z(\gamma)=Z(\gamma_{h})\cap Z(\gamma_{e}), \; 
	Z(\gamma_{h})=Z(\mathrm{Ad}(g)a).
	\label{eq:3.1.3ff}
\end{equation}

By \cite[Proposition 7.25]{knapp2002liegroupe}, $Z(\gamma)$ is 
reductive (possibly with several connected components). Set 
\begin{equation}
	\theta_{g}=C(g)\theta C(g^{-1}).
	\label{eq:3.1.4ff}
\end{equation}
Then $\theta_{g}$ defines a Cartan involution on $Z(\gamma)$. Let 
$K(\gamma)$ be the fixed point set of $\theta_{g}$ in $Z(\gamma)$, 
then
\begin{equation}
	K(\gamma)=Z(\gamma)\cap gKg^{-1}.
	\label{eq:ap3.1.6}
\end{equation}

Let $Z(\gamma)^{0}$, $K(\gamma)^{0}$ be the connected components of the identity of 
$Z(\gamma)$, $K(\gamma)$ respectively. By \cite[Theorem 
3.3.1]{bismut2011hypoelliptic}, 
\begin{equation}
	Z(\gamma)/K(\gamma)=Z(\gamma)^{0}/K(\gamma)^{0}.
	\label{eq:3.3.3ff}
\end{equation}
Moreover, $K(\gamma)$, $K(\gamma)^{0}$ are maximal compact subgroups 
of $Z(\gamma)$, $Z(\gamma)^{0}$ respectively. 

Taking the corresponding Lie algebras in \eqref{eq:3.1.3ff}, we have
\begin{equation}
	\z(\gamma)=\z(\gamma_{h})\cap \z(\gamma_{e}),\; 
	\z(\gamma_{h})=\z(\mathrm{Ad}(g)a).
	\label{eq:3.1.4ffb}
\end{equation}
Let $\kk(\gamma)\subset\z(\gamma)$ be the Lie algebra of 
$K(\gamma)$.  Put
\begin{equation}
	\pp(\gamma)=\z(\gamma)\cap\mathrm{Ad}(g)\pp.
	\label{eq:3.1.11ff}
\end{equation}
Then the Cartan decomposition of $\z(\gamma)$ with respect to 
$\theta_{g}$ is given by
\begin{equation}
	\z(\gamma)=\kk(\gamma)\oplus\pp(\gamma).
	\label{eq:3.1.12ff}
\end{equation}
Let 
$B_{\z(\gamma)}$ denote the restriction of $B$ on $\z(\gamma)\times 
\z(\gamma)$, then $B_{\z(\gamma)}$ is invariant under the adjoint 
action of $\theta_{g}$ on $\z(\gamma)$.
Moreover, $B_{\z(\gamma)}$ is positive on $\pp(\gamma)$, and negative on 
$\kk(\gamma)$. The splitting in \eqref{eq:3.1.12ff} is orthogonal 
with respect to $B_{\z(\gamma)}$.

%%%%%%%%%%%%%%%%%%%%%%%%%%%%%%%%%%%%%%%%%%%%%%%%%%%%%%%%%%%%%%%%%%%%%%%
\subsection{Symmetric space}\label{section3.2}
Set
\begin{equation}
	X=G/K.
	\label{eq:3.2.1ee}
\end{equation}
Then $X$ is a smooth manifold with the smooth structure induced 
by $G$. By definition, $X$ is diffeomorphism to $\pp$. 

Let $\omega^\g\in \Omega^1(G,\g)$ be the canonical left-invariant $1$-form on $G$. Then by \eqref{eq:0.2.1ugc}, 
\begin{equation}
	\omega^\g=\omega^\pp+\omega^\kk.
\end{equation}

Let $p: G\rightarrow X$ denote the obvious projection. Then $p$ is a $K$-principal bundle over $X$. Then $\omega^\kk$ is a connection form of this principal bundle. The associated curvature form 
\begin{equation}
	\Omega^\kk=d\omega^\kk+\frac{1}{2}[\omega^\kk,\omega^\kk]=-\frac{1}{2}[\omega^\pp,\omega^\pp].
	\label{eq:3.2.2ee}
\end{equation}

If $(E,\rho^E, h^E)$ is a finite dimensional unitary or Euclidean 
representation of $K$, then $F=G\times_K E$ defines a vector bundle 
over $X$ equipped with a metric $h^{F}$ induced by $h^{E}$ and a 
unitary or an Euclidean 
connection $\nabla^F$ induced by $\omega^\kk$. Note that $G$ acts on 
$(F,h^{F},\nabla^{F})\rightarrow X$ equivariantly on the left, more precisely, for 
$\gamma\in 
G$, $(g,v)\in G\times_{K} E$, the action of $\gamma$ on $F$ is represented by
\begin{equation}
	\gamma(g,v)=(\gamma g,v)\in G\times_{K} E.
	\label{eq:3.2.4vogel}
\end{equation}

In particular, we have the 
identification
\begin{equation}
	TX=G\times_K \pp,
	\label{eq:3.2.3ee}
\end{equation}
where the right-hand side is defined from the adjoint action of $K$ 
on $\pp$.
The bilinear form $B$ restricting to $\pp$ gives a Riemannian metric 
$g^{TX}$, and $\omega^\kk$ induces the associated Levi-Civita 
connection $\nabla^{TX}$. Then $G$ acts on $(X,g^{TX})$ 
isometrically. Let $d(\cdot, \cdot)$ denote the Riemannian 
distance on $X$.

Let $C(G,E)$ denote the set of continuous map from $G$ into $E$. If $k\in K$, $s\in C(G,E)$, put
\begin{equation}
	(k.s)(g)=\rho^E(k)s(gk).
	\label{eq:3.2.5bbs}
\end{equation}
Let $C_K(G,E)$ be the set of $K$-invariant maps in $C(G,E)$. Let $C(X,F)$ denote the continuous sections of $F$ over $X$. Then 
\begin{equation}
	C_K(G,E)=C(X,F).
	\label{eq:3.2.6bbs}
\end{equation}
Also $C^\infty_K(G,E)=C^\infty(X,F)$.

The Casimir operator $C^\g$ acting on $C^\infty(G,E)$ preserves 
$C^\infty_K(G,E)$, so it induces an operator $C^{\g,X}$ acting on 
$C^\infty(X,F)$. Let $\Delta^{H,X}$ be the Bochner Laplacian acting 
on $C^\infty(X,F)$ given by $\nabla^F$, and let $C^{\kk, 
E}\in\mathrm{End}(E)$ be the action of the Casimir $C^\kk$ on $E$ via $\rho^E$. The element $C^{\kk,E}$ induces an self-adjoint section of $\mathrm{End}(F)$ over $X$. Then
\begin{equation}
	C^{\g,X}=-\Delta^{H,X}+C^{\kk,E}.
	\label{eq:3.2.7bbsd}
\end{equation}

Let $C^{\kk,\pp}\in\mathrm{End}(\pp)$, 
$C^{\kk,\kk}\in\mathrm{End}(\kk)$ be the actions of $C^{\kk}$ 
acting on $\pp$, $\kk$ via the adjoint actions. Moreover, we 
can also view $C^{\kk,\pp}$ as a parallel section of $\mathrm{End}(TX)$.

If $A\in\mathrm{End}(E)$ commutes with $K$, then it can be viewed a parallel section of $\mathrm{End}(F)$ over $X$. Let 
$dx$ be the Riemannian volume element of $(X, g^{TX})$.
\begin{definition}
	Let $\mathcal{L}^X_A$ be the Bochner-like Laplacian acting on $C^\infty(X,F)$ given by
	\begin{equation}
		\mathcal{L}^X_A=\frac{1}{2}C^{\g,X} + 
		\frac{1}{16}\mathrm{Tr}^{\pp}[C^{\kk,\pp}]+\frac{1}{48}\mathrm{Tr}^{\kk}[C^{\kk,\kk}]+A.
		\label{eq:3.2.9ssd}
	\end{equation}
	For $t>0$, $x,x'\in X$, let $p^X_t(x,x')$ denote its heat kernel with 
	respect to $dx'$.
\end{definition}

Since $\mathcal{L}^{X}_{A}$ is $G$-invariant, then $p^X_t(x,x')$ 
lifts to a function $p^{X}_{t}(g,g')$ on $G\times G$ valued in 
$\mathrm{End}(E)$ such that for $g''\in G$, $k,k'\in K$,
\begin{equation}
	p^{X}_{t}(g''g,g''g')=p^{X}_{t}(g,g'), 
	p^{X}_{t}(gk,g'k')=\rho^{E}(k^{-1})p^{X}_{t}(g,g')\rho^{E}(k').
	\label{eq:3.2.10ffe}
\end{equation}
We set
\begin{equation}
	p^{X}_{t}(g)=p^{X}_{t}(1,g).
	\label{eq:3.2.11ffe}
\end{equation}
Then $p^{X}_{t}$ is a $K\times K$-invariant smooth function on $G$ 
valued in $\mathrm{End}(E)$. We will not distinguish the heat kernel $p^{X}_{t}(x,x')$ 
and the function $p^{X}_{t}(g)$ in the sequel.

%%%%%%%%%%%%%%%%%%%%%%%%%%%%%%%%%%%%%%%%%%%%%%%%%%%%%%%%%%%%%%%%%%%
\subsection{Bismut's formula for semisimple orbital 
integrals}\label{section3.3}

Let $dg$ be the left-invariant Haar measure on $G$ induced by 
$(\g,\langle\cdot,\cdot\rangle)$. Since $G$ is 
unimodular, then $dg$ is also right-invariant.
Let $dk$ be the Haar measure on $K$ induced by $-B|_{\kk}$, then
\begin{equation}
	dg=dxdk.
	\label{eq:3.2.9cc}
\end{equation}

Now let $\gamma\in G$ be a semisimple element given as in 
\eqref{eq:3.1.8vogel}.

By 
\cite[Definition 2.19.21]{eberlein1996geometry} and \cite[Theorem 
3.1.2]{bismut2011hypoelliptic}, $\gamma\in G$ is 
semisimple if and only if the displacement function $X\ni x\mapsto 
d(x,\gamma x)$ on $X$ associated with $\gamma$ can reach its minimum 
$m_{\gamma}\geq 0$ in $X$. In this case, the minimizing set $X(\gamma)$ of this 
displacement function is a geodesically convex submanifold of 
$X$, and by \cite[Theorem 
3.3.1]{bismut2011hypoelliptic},
\begin{equation}
	X(\gamma)\simeq Z(\gamma)^{0}/K(\gamma)^{0}=Z(\gamma)/K(\gamma).
	\label{eq:3.3.3vogel}
\end{equation}
Moreover, we have
\begin{equation}
	m_{\gamma}=|a|.
	\label{eq:3.3.2ff}
\end{equation}

Let $dy$ be the Riemannian volume element of $X(\gamma)$, and let 
$dz$ be the bi-invariant Haar measure on 
$Z(\gamma)$ induced by $B_{\z(\gamma)}$. Let $dk(\gamma)$ be the Haar 
measure on $K(\gamma)$ such that
\begin{equation}
	dz=dydk(\gamma).
	\label{eq:3.3.7ffd}
\end{equation}
Let $\mathrm{Vol}(K(\gamma)\backslash K)$ be the volume of 
$K(\gamma)\backslash K$ with respect to $dk, dk(\gamma)$. Then we have
\begin{equation}
	\mathrm{Vol}(K(\gamma)\backslash 
	K)=\frac{\mathrm{Vol}(K)}{\mathrm{Vol}(K(\gamma))}.
	\label{eq:3.3.7bisf}
\end{equation}

Let $dv$ be the $G$-left invariant measure on 
$Z(\gamma)\backslash G$ such that
\begin{equation}
	dg=dzdv.
	\label{eq:3.3.6ff}
\end{equation}

By \cite[Definition 4.2.2, Proposition 
4.4.2]{bismut2011hypoelliptic}, for $t>0$, the orbital integral
\begin{equation}
	\mathrm{Tr}^{[\gamma]}[\exp(-t\mathcal{L}^{X}_{A})]=\frac{1}{\mathrm{Vol}(K(\gamma)\backslash K)}\int_{Z(\gamma)\backslash G}\mathrm{Tr}^{E}[p^{X}_{t}(v^{-1}\gamma v)]dv
	\label{eq:3.3.8ffe}
\end{equation}
is well-defined. As indicated by the notation, it only depends on the 
conjugacy class $[\gamma]$ of $\gamma$ in $G$.

Using the theory of hypoelliptic Laplacian and the techniques from local 
index theory, Bismut obtained an explicit geometric formula for 
$\mathrm{Tr}^{[\gamma]}[\exp(-t\mathcal{L}^{X}_{A})]$ in 
\cite[Theorem 6.1.1]{bismut2011hypoelliptic} as well as its extension to the 
wave operators of $\mathcal{L}^{X}_{A}$ \cite[Section 6.3]{bismut2011hypoelliptic}. Now we describe in detail 
this formula. We may and we will assume that
\begin{equation}
	\gamma=e^{a}k,\; a\in\pp,\; k\in K,\; \mathrm{Ad}(k)a=a.
	\label{eq:3.3.10ff}
\end{equation}

Put 
\begin{equation}
	\z_0=\z(a),\;\;\pp_0=\ker{\ad(a)}\cap \pp,\;\;\kk_0=\ker{\ad(a)}\cap \kk.
	\label{eq:1.6.4ugc}
\end{equation}
Let $\z^{\perp}_0$, $\pp^{\perp}_0$, $\kk^{\perp}_0$ be the 
orthogonal vector spaces to $\z_0$, $\pp_0$, $\kk_0$ in $\g, \pp, 
\kk$ with respect to $B$. Then 
\begin{equation}
	\z_0=\pp_0\oplus \kk_0,\;\; \z^{\perp}_0=\pp^{\perp}_0\oplus \kk^{\perp}_0.
	\label{eq:1.6.6ugcd}
\end{equation}

By \cite[Eq. (3.3.6)]{bismut2011hypoelliptic}, 
\begin{equation}
	\z(\gamma)=\z_{0}\cap \z(k).
\end{equation}
Also $\pp(\gamma)$, $\kk(\gamma)$ are subspaces of $\pp_0$, $\kk_0$ respectively. Let $\z^{\perp}_{0}(\gamma)$, $\pp^{\perp}_{0}(\gamma)$, $\kk^{\perp}_{0}(\gamma)$ be the orthogonal spaces to $\z(\gamma)$, $\pp(\gamma)$, $\kk(\gamma)$ in $\z_0$, $\pp_0$, $\kk_0$. Then
\begin{equation}
	\z^{\perp}_{0}(\gamma)=\pp^{\perp}_{0}(\gamma)\oplus \kk^{\perp}_{0}(\gamma).
	\label{eq:1.6.7ugcd}
\end{equation}
Also the action $\mathrm{ad}(a)$ gives an isomorphism between $\pp^{\perp}_0$ and $\kk^{\perp}_0$.

For $Y_0^{\kk}\in \kk(\gamma)$, $\mathrm{ad}(Y^{\kk}_0)$ 
preserves $\pp(\gamma), \kk(\gamma), 
\pp^{\perp}_{0}(\gamma), \kk^{\perp}_{0}(\gamma)$, and it is 
an antisymmetric endomorphism with respect to the scalar product.

Recall that the function $\widehat{A}$ is given by 
\begin{equation}\label{eq:Ahatfunction}
	\widehat{A}(x)=\frac{x/2}{\sinh(x/2)}.
\end{equation}
Let $H$ be a finite-dimensional Hermitian vector space. If $B\in 
\mathrm{End}(H)$ is self-adjoint, then $\dfrac{B/2}{\sinh(B/2)}$ is a 
self-adjoint positive endomorphism. Put
\begin{equation}\label{eq:AhatB}
	\widehat{A}(B)=\det{}^{1/2}\bigg[\frac{B/2}{\sinh(B/2)}\bigg].
\end{equation}
In \eqref{eq:AhatB}, the square root is taken to be the positive square root.

If $Y^\kk_0\in\kk(\gamma)$, as explained in 
\cite[pp105]{bismut2011hypoelliptic}, the following function $A(Y^{\kk}_{0})$ has a natural square 
root that is analytic in $Y^{\kk}_{0}\in \kk(\gamma)$,
\begin{equation}
	A(Y^{\kk}_{0})=\frac{1}{\det 
	(1-\mathrm{Ad}(k))|_{\z^{\perp}_{0}(\gamma)}} \cdot\frac{\det 
	(1-\exp(-i\mathrm{ad}(Y_0^\kk))\mathrm{Ad}(k))|_{\kk^{\perp}_{0}(\gamma)}}{\det (1-\exp(-i\mathrm{ad}(Y_0^\kk))\mathrm{Ad}(k))|_{\pp^{\perp}_{0}(\gamma)}}. 
	\label{eq:AinJ}
\end{equation}
Its square root is denoted by
\begin{equation}
	\bigg[ \frac{1}{\det 
	(1-\mathrm{Ad}(k))|_{\z^{\perp}_{0}(\gamma)}} \cdot\frac{\det 
	(1-\exp(-i\mathrm{ad}(Y_0^\kk))\mathrm{Ad}(k))|_{\kk^{\perp}_{0}(\gamma)}}{\det (1-\exp(-i\mathrm{ad}(Y_0^\kk))\mathrm{Ad}(k))|_{\pp^{\perp}_{0}(\gamma)}}  \bigg]^{1/2}.
	\label{eq:Aroot}
\end{equation}
The value of \eqref{eq:Aroot} at $Y^{\kk}_{0}=0$ is taken to be such that 
\begin{equation}
	\frac{1}{\det(1-\mathrm{Ad}(k))|_{\pp^{\perp}_{0}(\gamma)}}.
	\label{eq:3.3.16ffed}
\end{equation}

We recall an important function $J_\gamma$ defined in \cite[Eq. (5.5.5)]{bismut2011hypoelliptic}.
\begin{definition}\label{def:3.3.1ss20}
	Let $J_{\gamma}(Y^{\kk}_{0})$ be the analytic function of 
	$Y_0^{\kk}\in \kk(\gamma)$ given by 
	\begin{equation}\label{Jfunction}
		\begin{split}
			&J_{\gamma}(Y_0^{\kk})=\frac{1}{|\det (1-\mathrm{Ad}(\gamma))|_{\z^{\perp}_0}|^{1/2}} \frac{\widehat{A}(i\mathrm{ad}(Y_0^\kk)|_{\pp(\gamma)})}{\widehat{A}(i\mathrm{ad}(Y_0^\kk)|_{\kk(\gamma)})}\\
			&\bigg[ \frac{1}{\det (1-\mathrm{Ad}(k))|_{\z^{\perp}_{0}(\gamma)}} \frac{\det (1-\exp(-i\mathrm{ad}(Y_0^\kk))\mathrm{Ad}(k))|_{\kk^{\perp}_{0}(\gamma)}}{\det (1-\exp(-i\mathrm{ad}(Y_0^\kk))\mathrm{Ad}(k))|_{\pp^{\perp}_{0}(\gamma)}}     \bigg]^{1/2}.
		\end{split}
	\end{equation}
\end{definition}

By \cite[Eq. (6.1.1)]{bismut2011hypoelliptic}, there 
exist $C_{\gamma}>0$, $c_{\gamma}>0$ such that if $Y^{\kk}_{0}\in 
\kk(\gamma)$,
\begin{equation}
	|J_{\gamma}(Y^{\kk}_{0})|\leq 
	C_{\gamma}e^{c_{\gamma}|Y^{\kk}_{0}|}.
	\label{eq:3.3.17ffed}
\end{equation}

Put $p=\dim \pp(\gamma)$, $q=\dim \kk(\gamma)$. Then $r=\dim 
\z(\gamma)=p+q$. By \cite[Theorem 6.1.1]{bismut2011hypoelliptic}, for 
$t>0$, we have
\begin{equation}
	\begin{split}
		&\mathrm{Tr}^{[\gamma]}[\exp(-t\mathcal{L}^{X}_{A})]\\
		&=\frac{e^{-\frac{|a|^2}{2t}}}{(2\pi 
		t)^{p/2}}\int_{\kk(\gamma)}J_{\gamma}(Y^{\kk}_{0})\mathrm{Tr}^{E}\big[\rho^{E}(k)\exp(-i\rho^{E}(Y^{\kk}_{0})-tA)\big]e^{-|Y^{\kk}_{0}|^{2}/2t}\frac{dY^{\kk}_{0}}{(2\pi t)^{q/2}}.
	\end{split}
	\label{eq:3.3.18ffed}
\end{equation}

\begin{remark}
	A generalization of Bismut's formula \eqref{eq:3.3.18ffed} to the twisted case is 
	obtained in \cite{liu:tel-01841334, LIU201974}. An extension 
	of this formula for considering arbitrary elements in the center of 
	enveloping algebra instead of Casimir operator 
	\eqref{eq:3.2.7bbsd} was obtained in \cite{bismut2019geometric} 
	by Bismut and Shen.
\end{remark}

%%%%%%%%%%%%%%%%%%%%%%%%%%%%%%%%%%%%%%%%%%%%%%%%%%%%%%%%%%%%%%%%%%%
\subsection{Compact locally symmetric spaces}\label{section3.4}
Let $\Gamma$ be a cocompact discrete subgroup of $G$. Then $\Gamma$ 
acts on $X$ isometrically and properly discontinuously. Then 
$Z=\Gamma\backslash X$ is compact second countable Hausdorff space. 

If $x\in X$, put
\begin{equation}
	\Gamma_x=\{\gamma\in \Gamma : \gamma x=x\}.
	\label{eq:3.2.9bbc}
\end{equation}
Then $\Gamma_{x}$ is a finite subgroup of $\Gamma$. 
Put
\begin{equation}
	r_{x}=\inf_{\gamma\in\Gamma - \Gamma_{x}}d(x,\gamma x).
	\label{eq:3.2.10bbc}
\end{equation}
Then we always have $r_{x}>0$. Set
\begin{equation}
	U_{x}=B(x,\frac{r_{x}}{4})\subset X.
	\label{eq:3.2.11bbc}
\end{equation}
If $x\in X$, $\gamma\in\Gamma$, we have
\begin{equation}
	r_{\gamma x}=r_{x},\; U_{\gamma x}=\gamma U_{x}.
	\label{eq:3.2.12bbc}
\end{equation}
It is clear that $\Gamma_{x}\backslash U_{x}$ can identified with a 
connected open subset of $Z$. 

Set
\begin{equation}
	S=\ker (\Gamma\rightarrow \mathrm{Diffeo}(X))=\Gamma\cap 
	\ker(K\xrightarrow{\mathrm{Ad}}\mathrm{Aut}(\pp)).
	\label{eq:3.2.9aa}
\end{equation}
Then $S$ is a finite subgroup of $\Gamma\cap K$, and a normal 
subgroup of $\Gamma$. 

\begin{remark}
	Note that $G_{\mathrm{ss}}$ is a connected noncompact simple linear Lie group, then
	\begin{equation}
		S=Z_G\cap \Gamma\cap K.
		\label{eq:3.4.16th19}
	\end{equation}
\end{remark}

Put
\begin{equation}
	\Gamma'=\Gamma/S.
	\label{eq:3.2.10aa}
\end{equation}
Then $\Gamma'$ acts on $X$ effectively and we have 
$Z=\Gamma'\backslash X$. 

If $x\in X$, we have
\begin{equation}
	S\subset \Gamma_{x},\; \Gamma'_{x}=\Gamma_{x}/S.
	\label{eq:3.3.7aa}
\end{equation}
Then the orbifold charts $(U_{x}, \Gamma'_{x}, \pi_{x}: 
U_{x}\rightarrow \Gamma'_{x}\backslash U_{x})_{x\in X}$ together with 
the action of $\Gamma'$ on these charts give an (effective) orbifold structure 
for
$Z$, so that $Z=\Gamma\backslash X$ is a compact orbifold with a 
Riemannian metric $g^{TZ}$ induced by $g^{TX}$.

By \cite[Lemma 1]{Selberg1960}, if $\gamma\in \Gamma$, then $\gamma$ is semisimple. Let 
$[\Gamma]$ denote the set of the conjugacy classes of $\Gamma$. If 
$\gamma\in \Gamma$, we say $[\gamma]\in [\Gamma]$  to be an elliptic 
class if $\gamma$ is elliptic. Let $\mathrm{E}[\Gamma]\subset [\Gamma]$ be 
the set of elliptic classes, then $\mathrm{E}[\Gamma]$ is always a finite 
set. If $\mathrm{E}[\Gamma]$ only contains the trivial conjugacy 
class $[1]$, i.e. $\Gamma$ is torsion free, then $Z$ is 
compact smooth manifold.

Let $[\Gamma']$ be the set of conjugacy classes in $\Gamma'$, and let 
$\mathrm{E}[\Gamma']$ denote the set of elliptic classes in 
$[\Gamma']$. If $\gamma'\in \Gamma'$, let 
$Z_{\Gamma'}(\gamma')$ denote the centralizer of $\gamma'$ in 
$\Gamma'$, and let $[\gamma']'$ denote the conjugacy class of 
$\gamma'$ in $\Gamma'$. If $\gamma'\in \Gamma'$ is elliptic, let $X(\gamma')$ be 
its fixed point set in $X$ on which $Z_{\Gamma'}(\gamma')$ acts 
isometrically and properly discontinuously (see \cite[Lemma 
2]{Selberg1960}). Note that if $\gamma\in\Gamma$ is a lift of 
$\gamma'\in \Gamma'$, then $X(\gamma)=X(\gamma')$, and $\gamma$ is elliptic if and only if 
$\gamma'$ is elliptic.

\begin{proposition}\label{prop:3.4.1lead}
	We have
	\begin{equation}
		Z_{\mathrm{sing}}=\Gamma'\backslash 
		\Gamma'\big(\cup_{[\gamma']'\in 
		\mathrm{E}[\Gamma']\backslash \{1\}} X(\gamma')\big)\subset Z.
		\label{eq:3.3.8aa}
	\end{equation}
	Moreover, we have
	\begin{equation}
		\Sigma Z=\cup_{[\gamma']'\in \mathrm{E}[\Gamma']\backslash 
		\{1\}} Z_{\Gamma'}(\gamma')\backslash X(\gamma').
		\label{eq:3.2.11aa}
	\end{equation}
	Note that the right-hand side of \eqref{eq:3.2.11aa} is a 
	disjoint union of compact orbifolds.
	
	If $\gamma'\in 
	\Gamma'$, put
	\begin{equation}
		S'(\gamma')=\ker (Z_{\Gamma'}(\gamma')\rightarrow 
		\mathrm{Diffeo}(X(\gamma'))).
		\label{eq:3.5.15oct9s}
	\end{equation}
	Then $|S'(\gamma')|$ is the multiplicity of the connected component 
	$Z_{\Gamma'}(\gamma')\backslash X(\gamma')$ in $\Sigma Z$
\end{proposition}
\begin{proof}
	Note that $z\in Z$ with a lift $x\in X$ belongs to 
	$Z_{\mathrm{sing}}$ if and only if the stabilizer $\Gamma'_{x}$ 
	is nontrivial. Thus $x$ is a fixed point of some $\gamma'\in 
	\Gamma'$, from which \eqref{eq:3.3.8aa} follows. By definition in 
	Subsection \ref{subs:2.1}, we get the rest part of this 
	proposition. This completes the proof.
\end{proof}

Note that $\Gamma\backslash G$ is a compact smooth homogeneous space 
equipped with a right action of $K$. Moreover, the action of $K$ is 
almost free, i.e. for each $\bar{g}\in \Gamma\backslash G$, the 
stabilizer $K_{\bar{g}}$ is finite. Then the quotient space 
$(\Gamma\backslash G)/K$ also have a natural orbifold structure, 
which, after examining the local charts, is equivalent to $Z$. 

Let $d\bar{g}$ be the volume element on $\Gamma\backslash G$ induced 
by $dg$. By \eqref{eq:3.2.9cc}, we get
\begin{equation}
	\mathrm{Vol}(\Gamma\backslash G)=\frac{\mathrm{Vol}(K)}{|S|}\mathrm{Vol}(Z).
	\label{eq:3.3.13ddd}
\end{equation}

In the context of geometry, we 
have many interesting cases where $S=\{1\}$. For instance, given a Riemannian symmetric space $(X,g^{TX})$ of noncompact 
type, let $G=\mathrm{Isom}(X)^{0}$ be the connected component of 
identity of the Lie group of isometries of $X$. By \cite[Proposition 
2.1.1]{eberlein1996geometry}, $G$ is a semisimple Lie group with 
trivial center (which might not be 
linear, but we do not need that linearity for the geometry of $Z$). We refer to 
\cite[Chapter 2]{eberlein1996geometry} and \cite[Chapter 
3]{bismut2011hypoelliptic} for more details. This way, any subgroup 
of $G$ acts on $X$ effectively. In particular, if $\Gamma$ is 
a cocompact discrete subgroup of $G$, then $Z=\Gamma\backslash X$ is 
a compact good orbifold with the orbifold fundamental group $\Gamma$. By 
\eqref{eq:3.2.11aa}, we have
\begin{equation}
	\Sigma Z=\cup_{[\gamma]\in E[\Gamma]\backslash \{1\}} \Gamma\cap 
	Z(\gamma)\backslash X(\gamma).
	\label{eq:3.4.22bbs}
\end{equation}
In general, by \cite[Ch.V \S 4, Theorem 
4.1]{helgason1979differential}, $G=\mathrm{Isom}(X=G/K)^{0}$ if 
and only if $K$ acts on $\pp$ effectively. 

\begin{remark}\label{eq:3.4.3vogel}
	Note that, as mentioned in Remark \ref{rk:2.1.3vogel}, when 
	$S\neq\{1\}$, we can 
	also consider $Z=\Gamma\backslash X$ as an ineffective orbifold by taking the 
	action of $\Gamma$ instead of $\Gamma'$ on the local charts. This 
	way, the role of the above $Z\cup \Sigma Z$ is replaced by the 
	inertia groupoid defined in \cite[Example 2.5]{Adem_2007}, which 
	is exactly
	\begin{equation}
		\cup_{[\gamma]\in E[\Gamma]} \Gamma\cap 
		Z(\gamma)\backslash X(\gamma).
		\label{eq:3.4.14vogel}
	\end{equation}
	It is a much natural object to use in the context here, for instance, 
	for the Selberg's trace formula in the next subsection. 
	In the 
	problems concerned by us, these two 
	point-views on $Z$ are equivalent. 
\end{remark}

If $\rho:\Gamma'\rightarrow \mathrm{GL}(\C^{k})$ is a representation 
of $\Gamma'$, which can be viewed as a representation of $\Gamma$ via 
the projection $\Gamma\rightarrow \Gamma'=\Gamma/S$, then 
$F=\Gamma'\backslash (X\times \C^{k})$ is a proper flat 
orbifold vector bundle on $Z$ with the flat connection $\nabla^{F,f}$ 
induced from the exterior differential $d^{X}$ on 
$\C^{k}$-valued functions. By \cite[Theorem 
2.35]{2017arXiv170408369S}, all the proper orbifold vector bundle on 
$Z$ of rank $k$ comes from this way.

Now let $\rho:\Gamma\rightarrow \mathrm{GL}(\C^{k})$ be a 
representation of $\Gamma$, we do not assume that it comes from a 
representation of $\Gamma'$. We still have a flat orbifold vector 
bundle $(F=\Gamma\backslash (X\times \C^{k}),\nabla^{F,f})$ on $Z$, which may not be proper in general. 
Note that $\Gamma$ acts on $C^{\infty}(X,\C^{k})$ so that if
$\varphi\in C^{\infty}(X,\C^{k})$, $\gamma\in\Gamma$, then
\begin{equation}
	(\gamma\varphi)(x) =\rho(\gamma)\varphi(\gamma^{-1}x).
	\label{eq:3.4.16bdf}
\end{equation}
Let $C^{\infty}(X,\C^{k})^{\Gamma}$ denote the $\Gamma$-invariant 
sections in $C^{\infty}(X,\C^{k})$. Then
\begin{equation}
	C^{\infty}(Z,F)=C^{\infty}(X,\C^{k})^{\Gamma}.
	\label{eq:eq:3.4.15bk}
\end{equation}

\begin{definition}\label{def:properFvogel}
	Let $(V,\rho^{V})$ be the isotypic 
	component of $(\C^{k},\rho|_S)$ corresponding to the trivial representation of $S$ on 
	$\C$, i.e. the maximal $S$-invariant subspace of $\C^{k}$ via 
	$\rho$. Set 
	\begin{equation}
		F^{\mathrm{pr}}=\Gamma\backslash(X\times V).
	\end{equation}
	It is clear that $F^{\mathrm{pr}}$ is a proper flat orbifold 
	vector bundle on $Z$.
\end{definition}

We have the following results.
\begin{proposition}\label{prop:3.4.2octs}
	We have
	\begin{equation}
		C^{\infty}(Z,F)=C^{\infty}(Z,F^{\mathrm{pr}}).
		\label{eq:3.4.14bf}
	\end{equation}
	
	In particular, if $\rho|_{S}:S\rightarrow 
	\mathrm{GL}(\C^{k})$ does not have the isotypic component of the
	trivial representation of $S$ on $\C$, then
	\begin{equation}
		C^{\infty}(Z,F)=\{0\}.
		\label{eq:3.4.15bf}
	\end{equation}
\end{proposition}

Let $(E,\rho^{E})$ be a finite dimensional complex representation of 
$G$. When restricting to $\Gamma$, $K$, we get the corresponding 
representations of $\Gamma$, $K$ respectively, which are still 
denoted by $\rho^{E}$. As discussed in Subsection \ref{section3.2}, 
associated with $K$-representation $(E,\rho^{E})$
we define a homogeneous vector bundle $F=G\times_{K} E$ on $X$. 
Moreover, $G$ acts on $F$ equivariantly. By taking a 
$\Gamma$-quotient on the left, it descends to an orbifold vector 
bundle on $Z$, which we still denote by the same notation.

The map $(g,v)\in G\times_{K} E\rightarrow (pg, \rho^{E}(g)v)\in 
X\times E$ gives a canonical trivialization of $F$ over $X$. This 
trivialization provides a flat connection $\nabla^{X,F,f}$ for 
$F\rightarrow X$, which is $G$-invariant. Then it descends to a 
flat connection $\nabla^{Z,F,f}$ on the orbifold vector bundle $F$ 
over $Z$. Moreover, 
the above trivialization of $F\rightarrow X$ implies that the flat 
orbifold vector bundle $(F,\nabla^{Z,F,f})$ is exactly the one given by 
$\Gamma\backslash(X\times E)$ with the flat connection $\nabla^{F,f}$ induced by $d^{X}$. 
We will always use the notation $\nabla^{F,f}$ for the above flat 
connection. By 
\eqref{eq:3.2.6bbs}, \eqref{eq:eq:3.4.15bk}, we get
\begin{equation}
	C^{\infty}(Z,F)=C_{K}^{\infty}(G,E)^{\Gamma}.
	\label{eq:3.4.21bbs}
\end{equation}

%%%%%%%%%%%%%%%%%%%%%%%%%%%%%%%%%%%%%%%%%%%%%%%%%%%%%%%%%%%%%%%%%%%%%%%%
\subsection{Selberg's trace formula}\label{section3.5}
Let $Z$ be the compact locally symmetric space discussed in 
Subsection \ref{section3.4}, and let $(F,h^{F},\nabla^{F})$ be a Hermitian 
vector bundle on $X$ defined by a unitary representation 
$(E,\rho^{E})$ of $K$. As said before, $(F,h^{F},\nabla^{F})$ descends to 
a Hermitian orbifold vector bundle on $Z$. Recall the Bochner-like 
Laplacian $\mathcal{L}^{X}_{A}$ is defined by \eqref{eq:3.2.9ssd}. 
Since it commutes with $G$, then it descends to a Bochner-like 
Laplacian $\mathcal{L}^{Z}_{A}$ acting on $C^{\infty}(Z,F)$. 

Here the convergences of the integrals and infinite sums are already 
guaranteed by the results in \cite[Chapters 2 \& 
4]{bismut2011hypoelliptic} and in \cite[Section 4D]{Shen_2016}. 

For 
$t>0$, let $p^{Z}_{t}(z,z')$, $z,z'\in Z$ be the heat kernel of 
$\mathcal{L}^{Z}_{A}$ over $Z$ with respect to $dz'$. If $z,z'$ are identified with their 
lifts in $X$, then
\begin{equation}
	p^{Z}_{t}(z,z')=\frac{1}{|S|}\sum_{\gamma\in\Gamma} \gamma 
	p^{X}_{t}(\gamma^{-1}z,z')=\frac{1}{|S|}\sum_{\gamma\in\Gamma} 
	p^{X}_{t}(z,\gamma z')\gamma.
	\label{eq:3.5.1ksd}
\end{equation}
Note that the action of $\gamma$ on $F_{\gamma^{-1}z}$ or on
the metric dual of
$F_{z'}$ is given as in \eqref{eq:3.2.4vogel}.

Since $Z$ is compact, then for $t>0$, $\exp(-t\mathcal{L}^{Z}_{A})$ 
is trace class. We have
\begin{equation}
	\mathrm{Tr}[\exp(-t\mathcal{L}^{Z}_{A})]=\int_{Z}\mathrm{Tr}^{F}[p^{Z}_{t}(z,z)]dz.
	\label{eq:3.5.3ksd}
\end{equation}

Combining \eqref{eq:3.2.10ffe}, \eqref{eq:3.2.11ffe}, 
\eqref{eq:3.3.13ddd} and \eqref{eq:3.5.1ksd}, \eqref{eq:3.5.3ksd}, 
and proceeding as in \cite[Eqs. 
(4.8.8)-(4.8.12)]{bismut2011hypoelliptic}, we 
get
\begin{equation}
	\begin{split}
		\mathrm{Tr}[\exp(-t\mathcal{L}^{Z}_{A})]&=\frac{1}{\mathrm{Vol}(K)}\int_{\Gamma\backslash G}\sum_{\gamma\in\Gamma} \mathrm{Tr}^{E}[p^{X}_{t}(\bar{g}^{-1}\gamma\bar{g})]d\bar{g}\\
		&=\sum_{[\gamma]\in[\Gamma]}\frac{\mathrm{Vol}(\Gamma\cap 
		Z(\gamma)\backslash Z(\gamma))}{\mathrm{Vol}(K(\gamma))}\mathrm{Tr}^{[\gamma]}[\exp(-t\mathcal{L}^{X}_{A})].
	\end{split}
	\label{eq:3.5.4ksd}
\end{equation}

Take $\gamma\in \Gamma$. Recall that $X(\gamma)=Z(\gamma)/K(\gamma)$ 
defined in Subsection \ref{section3.3}. Then $K(\gamma)$ acts on $Z(\gamma)$ on the right, which induces an 
action on $\Gamma\cap Z(\gamma)\backslash Z(\gamma)$ on the right. Set
\begin{equation}
	S(\gamma)=\ker (\Gamma\cap Z(\gamma)\rightarrow 
	\mathrm{Diffeo}(X(\gamma))).
	\label{eq:3.3.12aa}
\end{equation}
Then $S(\gamma)$ represents the 
isotropy group of the principal orbit type for the right action of 
$K(\gamma)$ on $\Gamma\cap Z(\gamma)\backslash Z(\gamma)$.
As in \eqref{eq:3.3.13ddd}, we have
\begin{equation}
	\mathrm{Vol}(\Gamma\cap Z(\gamma)\backslash 
	Z(\gamma))=\frac{\mathrm{Vol}(K(\gamma))}{|S(\gamma)|}\mathrm{Vol}(\Gamma\cap Z(\gamma)\backslash X(\gamma)).
	\label{eq:3.5.13ksd}
\end{equation}

\begin{theorem}\label{thm:3.5.2kkss}
	For $t>0$, we have the following identity,
	\begin{equation}
		\mathrm{Tr}[\exp(-t\mathcal{L}^{Z}_{A})]=\sum_{[\gamma]\in[\Gamma]}\frac{\mathrm{Vol}(\Gamma\cap Z(\gamma)\backslash X(\gamma))}{|S(\gamma)|} \mathrm{Tr}^{[\gamma]}[\exp(-t\mathcal{L}^{X}_{A})].
		\label{eq:3.5.14ksd}
	\end{equation}
\end{theorem}
\begin{proof}
	This is a direct consequence of \eqref{eq:3.5.4ksd} and \eqref{eq:3.5.13ksd}.
\end{proof}

In the case where $S=1$, the trace formula 
\eqref{eq:3.5.14ksd} shows clearly the different contributions from $Z$ 
and from each components of $\Sigma Z$. 
Then combining 
\eqref{eq:3.2.11aa}, \eqref{eq:3.5.14ksd} with the results in \cite[Theorem 7.8.2]{bismut2011hypoelliptic}
\cite[Theorem 7.7.1]{liu:tel-01841334}, we can recover \eqref{eq:1.2.12alpha} 
for $Z$. If we use the same settings as in \cite[Sections 
7.1, 7.2]{bismut2011hypoelliptic} and we use instead the results in 
\cite[Theorem 7.7.1]{bismut2011hypoelliptic}, then we can 
recover the Kawasaki's local index theorem \cite{MR527023} for $Z$. 
By taking account of Remarks \ref{rk:2.1.3vogel} \& 
\ref{eq:3.4.3vogel}, the above considerations also hold even for 
$S\neq\{1\}$.

%%%%%%%%%%%%%%%%%%%%%%%%%%%%%%%%%%%%%%%%%%%%%%%%%%%%%%%%%%%%%%%%%%%%%%%%%%%%
%%%%%%%%%%%%%%%%%%%%%%%%%%%%%%%%%%%%%%%%%%%%%%%%%%%%%%%%%%%%%%%%%%%%%%%%%%%%%
\section{Analytic torsions for compact locally symmetric spaces}
\label{section3.6}
In this section, we explain how to make use of Bismut's formula 
\eqref{eq:3.3.18ffed} and Selberg's trace formula 
\eqref{eq:3.5.14ksd} to study the analytic torsions of $Z$. We continue using the same settings 
as in 
Section \ref{section2}. We will see that by a vanishing result on the analytic torsion, only the 
case $\delta(G)=1$ remains interesting. For studying this case, more 
tools will be carried out in 
Sections \ref{section4} \& \ref{section6paris}.

%%%%%%%%%%%%%%%%%%%%%%%%%%%%%%%%%%%%%%%%%%%%%%%%%%%%%%%%%%%%%%%%%%%%%%%%%
\subsection{A vanishing result on the analytic 
torsions}\label{subsection4.1paris}
Recall that $G$ is a connected linear real reductive Lie group. 
Recall that $\z_{\g}$ is the center of $\g$. Set
\begin{equation}
	\z_{\pp}=\z_{\g}\cap \pp,\; \z_{\kk}=\z_{\g}\cap\kk.
	\label{eq:3.6.1oct19}
\end{equation}
Then
\begin{equation}
	\z_{\g}=\z_{\pp}\oplus \z_{\kk},\; Z_{G}=\exp(\z_{\pp})(Z_{G}\cap 
	K).
	\label{eq:3.6.2}
\end{equation}

Let $T$ be a maximal torus of $K$ with Lie algebra $\kt$, put
\begin{equation}\label{eq:defofb}
	\kb=\{f\in\pp\,:\, [f,\kt]=0\}.
\end{equation}
It is clear that
\begin{equation}
	\z_{\pp}\subset \kb.
	\label{eq:3.6.2bis19}
\end{equation}

Put $\kh=\kb\oplus\kt$, then $\kh$ is a Cartan subalgebra of $\g$; 
let $H$ be analytic subgroup of $G$ associated with $\kh$, then it is 
also a Cartan subgroup of $G$ (see \cite[p.129 and Theorem 5.22 (b)]{knapp1986representation}).
Moreover, $\dim\kt$ is just the complex rank of $K$, and $\dim \kh$ 
is the complex rank of $G$. 

\begin{definition}
	Using the above notations, the deficiency of $G$, or the fundamental rank of $G$ is defined as
	\begin{equation}\label{eq:deltaG}
		\delta(G)=\mathrm{rk}_{\C}G-\mathrm{rk}_{\C}K=\dim_\R \kb.
	\end{equation}
	The number $m-\delta(G)$ is even.
\end{definition}

The following result was proved in \cite[Proposition 3.3]{Shen_2016}.
\begin{proposition}
	If $\gamma\in G$ is semisimple, then 
	\begin{equation}
		\delta(G)\leq \delta(Z(\gamma)^{0}).
		\label{eq:3.6.5oct19}
	\end{equation}
	The two sides of \eqref{eq:3.6.5oct19} are equal if and only if 
	$\gamma$ can be conjugated into $H$.
\end{proposition}

Recall that $\ku=\sqrt{-1}\pp\oplus\kk$ is the compact form of $G$, 
and that $U\subset G_{\C}$ is the analytic subgroup with Lie algebra 
$\ku$.
Let $U\ku$, $U\g_\C$ be the enveloping algebras of $\ku$, $\g_\C$ respectively.
Then $U\g_{\C}$ can 
be identified with the left-invariant holomorphic differential 
operators on $G_{\C}$. 
Let $C^\ku\in U\ku$ be the Casimir operator of $\ku$ associated with 
$B$, then
\begin{equation}
	C^\ku=C^\g\in U\g\cap U\ku\subset U\g_{\C}.
	\label{eq:7.7.1hh}
\end{equation}

In the sequel, we always assume that $U$ is compact, this is the case when $G$ has 
compact center. 

\begin{proposition}[Unitary trick]\label{prop:unitarytrick}
	Assume that $U$ is compact. Then any irreducible finite 
	dimensional (analytic) complex representation of $U$ extends uniquely to an 
	irreducible finite dimensional complex representation of $G$ such that their induced representations of Lie algebras are compatible. 
\end{proposition}

We now fix a unitary representation $(E,\rho^{E}, h^{E})$ of $U$, and we 
extend it to a representation of $G$, whose restriction to $K$ is 
still unitary.  Put $F=G\times_K E$ with the Hermitian metric $h^{F}$ induced by 
$h^{E}$. Let $\nabla^F$ be the Hermitian connection 
induced by the connection form $\omega^\kk$. 

Furthermore, as explained in the last part of Subsection 
\ref{section3.4}, $F$ is equipped with a canonical flat connection $\nabla^{F,f}$ as 
follows,
\begin{equation}
	\nabla^{F,f}=\nabla^{F}+ \rho^E(\omega^\pp).
	\label{eq:5.5.9bs}
\end{equation}
If $G$ has compact 
center, then $(F,h^{F},\nabla^{F,f})$ is a unimodular flat vector bundle.

Let $(\Omega_c^\bullet(X,F), d^{X,F})$ be the (compactly supported) de Rham 
complex twisted by $F$. Let $d^{X,F,*}$ be the adjoint operator of $d^{X,F}$ 
with respect to the $L_2$ metric on $\Omega_c^\bullet(X,F)$. 
The de Rham-Hodge operator $\mathbf{D}^{X,F}$ of this de Rham complex is given by
\begin{equation}\label{eq5.90}
	\mathbf{D}^{X,F}=d^{X,F}+d^{X,F,*}.
\end{equation}

The Clifford algebras $c(TX)$, 
$\widehat{c}(TX)$ act on $\Lambda^\bullet(T^*X)$. 
We still use $e_1$, $\cdots$, $e_m$ to denote an orthonormal basis of 
$\pp$ or $TX$, and let $e^1$, $\cdots$, $e^m$ be the corresponding 
dual basis of $\pp^{*}$ or $T^*X$.

Let $\nabla^{\Lambda^\bullet(T^*X)\otimes F,u}$ be the unitary 
connection on $\Lambda^\bullet(T^*X)\otimes F$ 
induced by $\nabla^{TX}$ and $\nabla^F$. 
Then the standard Dirac operator is given by
\begin{equation}
	D^{X,F}=\sum^m_{j=1} c(e_j)\nabla^{\Lambda^\bullet(T^*X)\otimes F,u}_{e_j}.
	\label{eq:5.5.14bs}
\end{equation}

By \cite[Eq.(8.42)]{BMZ2015toeplitz}, we have
\begin{equation}
	\mathbf{D}^{X,F}= D^{X,F}+ \sum^m_{j=1} \widehat{c}(e_j)\rho^E(e_j).
	\label{eq:5.5.15bs}
\end{equation}

In the same time, as explained in Subsection \ref{section3.2}, $C^{\g}$ descends to an elliptic differential operator $C^{\g,X}$ acting on 
$C^\infty(X,\Lambda^\bullet(T^*X)\otimes F)$. 
As in \eqref{eq:3.2.9ssd}, we put
\begin{equation}
	\mathcal{L}^{X,F}=\frac{1}{2}C^{\g,X}+\frac{1}{16}\mathrm{Tr}^{\pp}[C^{\kk,\pp}]+\frac{1}{48}\mathrm{Tr}^{\kk}[C^{\kk,\kk}].
	\label{eq:7.7.14sss}
\end{equation}
For simplicity, we will always put
\begin{equation}
	\beta_{\g}=\frac{1}{16}\mathrm{Tr}^{\pp}[C^{\kk,\pp}]+\frac{1}{48}\mathrm{Tr}^{\kk}[C^{\kk,\kk}]\in\R.
	\label{eq:4.1.13vogel}
\end{equation}

By \cite[Proposition 8.4]{BMZ2015toeplitz}, we have
\begin{equation}\label{eq:7.7.7}
	\frac{\mathbf{D}^{X,F, 2}}{2}= \mathcal{L}^{X,F}-
	\frac{1}{2}C^{\g,E}-\beta_{\g}=:\mathcal{L}^{X,F}_{A},
\end{equation}
where $A=-
\frac{1}{2}C^{\g,E}-\beta_{\g}$.

Let $\gamma\in G$ be a semisimple element. In the sequel, we may assume that 
\begin{equation}
	\gamma=e^ak, a\in \pp, k\in K, \mathrm{Ad}(k)a=a.
	\label{eq:s4formst}
\end{equation}
We also use the same notation as in Subsection \ref{section3.3}.

Recall that $p=\dim\pp(\gamma)$, $q=\dim \kk(\gamma)$. 
By \eqref{eq:3.3.18ffed} and \eqref{eq:7.7.7}, we have
\begin{equation}
	\begin{split}
		&\mathrm{Tr_s}^{[\gamma]}\big[(N^{\Lambda^\bullet(T^*X)}-\frac{m}{2})\exp(-t\mathbf{D}^{X,F,2}/2)\big]\\
		&=\frac{e^{-\frac{|a|^2}{2t}}}{(2\pi t)^{p/2}} 
		\exp\big(t\beta_{\g}) \int_{\kk(\gamma)} 
		J_\gamma(Y^\kk_0)\mathrm{Tr_s}^{\Lambda^\bullet(\pp^*)}\big[(N^{\Lambda^\bullet(\pp^*)}-\frac{m}{2})\mathrm{Ad}(k)\exp(-i\mathrm{ad}(Y^\kk_0))\big]\\
		&\qquad\qquad\qquad\qquad\qquad
		\cdot\mathrm{Tr}^{E}[\rho^{E}(k)\exp(-i\rho^{E}(Y^\kk_0)+\frac{t}{2}C^{\ku,E})] e^{-|Y^\kk_0|^2/2t}\frac{dY^\kk_0}{(2\pi t)^{q/2}}.
	\end{split}
	\label{eq:6.2.8pl}
\end{equation}

Now we take a cocompact discrete subgroup $\Gamma\subset G$. Then 
$Z=\Gamma\backslash X$ is a compact locally symmetric orbifold. We use the same notation as in 
Subsections \ref{section3.4} \& \ref{section3.5}. Then we get a flat orbifold vector bundle $(F,\nabla^{F,f},h^{F})$ on 
$Z$. Furthermore, $\mathbf{D}^{X,F}$ descends to the corresponding de 
Rham - Hodge operator 
$\mathbf{D}^{Z,F}$ acting on $\Omega^{\cdot}(Z,F)$. Let 
$\mathcal{T}(Z,F)$ denote the associated 
analytic torsion as in Definition \ref{def:2.2.5ss20}, i.e.,
\begin{equation}
	\mathcal{T}(Z,F)=\mathcal{T}(g^{TZ}, \nabla^{F,f},h^{F}).
	\label{eq:4.1.20paris2020}
\end{equation}

As explained in Subsection \ref{subs:2.2}, for 
computing $\mathcal{T}(Z,F)$, it is enough to evaluate 
\begin{equation}
	\mathrm{Tr_{s}}[(N^{\Lambda^{\bullet}(T^{\ast}Z)}-\frac{m}{2})\exp(-t\mathbf{D}^{Z,F,2}/2)], \; t>0.
	\label{eq:4.0.23}
\end{equation}
Then we apply the Selberg's trace formula in Theorem 
\ref{thm:3.5.2kkss}. We get
\begin{equation}
	\begin{split}
		&\mathrm{Tr_{s}}[(N^{\Lambda^{\bullet}(T^{\ast}Z)}-\frac{m}{2})\exp(-t\mathbf{D}^{Z,F,2}/2)]\\
		&=\sum_{[\gamma]\in[\Gamma]}\frac{\mathrm{Vol}(\Gamma\cap 
		Z(\gamma)\backslash X(\gamma))}{|S(\gamma)|} \mathrm{Tr}^{[\gamma]}[(N^{\Lambda^{\bullet}(T^{*}X)}-\frac{m}{2})\exp(-t\mathbf{D}^{X,F,2}/2)].
	\end{split}
	\label{eq:4.0.24kkss}
\end{equation}

As in \cite[Remark 
8.7]{BMZ2015toeplitz}, by \cite[Theorems 5.4 \& 5.5, Remark 5.6]{Ma2017bourbaki}, we have the following vanishing theorem on 
$\mathcal{T}(Z,F)$.
\begin{theorem}\label{thm:4.1.4paris20}
	If $m$ is even, or if $m$ is odd and $\delta(G)\geq 3$, then
	\begin{equation}
		\mathcal{T}(Z,F)=0.
		\label{eq:3.6.14s}
	\end{equation}
\end{theorem}
\begin{proof}
	By \cite[Theorem 7.9.1]{bismut2011hypoelliptic}, \cite[Theorem 5.4]{Ma2017bourbaki}, and use instead 
	\eqref{eq:4.0.24kkss}, we get that under the assumptions in this theorem, for $t>0$,
	\begin{equation}
		\mathrm{Tr_{s}}[(N^{\Lambda^{\bullet}(T^{*}Z)}-\frac{m}{2})\exp(-t\mathbf{D}^{Z,F,2}))]=0.
		\label{eq:3.6.19oct19}
	\end{equation}
	Then \eqref{eq:3.6.14s} follows from the definition of $\mathcal{T}(Z,F)$.
\end{proof}

Therefore, the only nontrivial case is that $\delta(G)=1$, so that 
$m$ is odd. If $\gamma\in G$ is of the form \eqref{eq:s4formst}.
Let $\kt(\gamma)\subset \kk(\gamma)$ be a Cartan subalgebra. 
Put
\begin{equation}
	\kb(\gamma)=\{v\in \pp(k)\,:\, [v,\kt(\gamma)]=0\},\; 
	\kh(\gamma)_{\pp}=\kb(\gamma)\cap \pp(\gamma).
\end{equation}
In particular, $a\in\kb(\gamma)$.
Then $\kh(\gamma)=\kh(\gamma)_{\pp}\oplus \kt(\gamma)$ is a Cartan subalgebra of $\z(\gamma)$.

Recall that $H$ is a maximally compact Cartan subgroup of $G$. The 
following result is just an analogue of \cite[Theorem 
4.12]{Shen_2016} and \cite[Theorem 7.9.1]{bismut2011hypoelliptic}. 
\begin{proposition}\label{prop:6.1.1ss}
	If $\delta(G)=1$, if $\gamma$ is semisimple and can not be 
	conjugated into $H$ by an element in $G$, then 
	\begin{equation}
		\mathrm{Tr_s}^{[\gamma]}\big[(N^{\Lambda^{\bullet}(T^*X)}-\frac{m}{2})\exp(-t\mathbf{D}^{X,F,2}/2)\big]=0.
		\label{eq:6.2.7pl}
	\end{equation}
\end{proposition}
\begin{proof}	
	Let $\kt$ be a Cartan subalgebra of $\kk$ containing $\kt(\gamma)$. 
	Then $\kb\subset \kb(\gamma)$.  If $a\notin \kb$, then $\dim 
	\kb(\gamma)\geq 2$. Therefore, by \cite[Eq. (4-44)]{Shen_2016}, for $Y^{\kk}_{0}\in\kk(\gamma)$, we have
	\begin{equation}
		\mathrm{Tr_s}^{\Lambda^{\bullet}(\pp^*)}\big[(N^{\Lambda^{\bullet}(\pp^*)}-\frac{m}{2})\mathrm{Ad}(k)\exp(-i\mathrm{ad}(Y^\kk_0))\big]=0.
		\label{eq:6.2.7cpl}
	\end{equation}
	This implies \eqref{eq:6.2.7pl}. The proof is completed.
\end{proof}

Set
\begin{equation}
	\g'=\z_{\kk}\oplus \g_{\mathrm{ss}}.
	\label{eq:4.1.26std}
\end{equation}
Then $\g'$ is an ideal of $\g$. Let $G'$ be the analytic subgroup of 
$G$ associated with $\g'$, which is closed and has a compact center 
(see \cite[Proposition 7.27]{knapp2002liegroupe}). The group $K$ 
is still a maximal subgroup of $G'$. Let $U'\subset 
U$ be the compact form of $G'$ with Lie algebra $\ku'$, then
\begin{equation}
	\ku=\sqrt{-1}\z_{\pp}\oplus\ku'.
	\label{eq:4.1.27bio}
\end{equation}

Now we assume that $\delta(G)=1$ and that $G$ has 
noncompact center, so that $\kb=\z_{\pp}$ has dimension $1$. 
Then $\delta(G')=0$. Under the hypothesis that $U$ is compact, then 
up to a finite cover, we may write 
\begin{equation}
	U\simeq \mathbb{S}^{1}\times U'.
	\label{eq:4.1.27std}
\end{equation}

We take $a_{1}\in\kb$ with $|a_{1}|=1$. If $(E,\rho^{E})$ is an 
irreducible unitary representation of $U$, then $\rho^{E}(a_{1})$ 
acts on $E$ by a real scalar operator. Let $\alpha_{E}\in\R$ be such 
that
\begin{equation}
	\rho^{E}(a_{1})=\alpha_{E}\mathrm{Id}_{E}.
	\label{eq:4.1.28std}
\end{equation}

Put $X'=G'/K$. Then $X'$ is an even-dimensional symmetric space (of 
noncompact type). We identify $\z_{\pp}$ with a 
real line $\R$, then
\begin{equation}
	G=\R\times G',\; X=\R\times X'.
	\label{eq:3.6.18oct19}
\end{equation}

In this case, the evaluation for analytic torsions can be made more 
explicit. If $\gamma\in G'$, let $X'(\gamma)$ denote the minimizing 
set of $d_{\gamma}(\cdot)$ in $X'$, so that
\begin{equation}
	X(\gamma)=\R\times X'(\gamma).
	\label{eq:ssok}
\end{equation}

Let $[\cdot]^{\mathrm{max}}$ denote the coefficient of a differential 
form (valued in $o(TX')$) on $X'$ of the corresponding Riemannian volume form. Similarly, for $k\in 
T$, let $[\cdot]^{\mathrm{max}(k)}$ denote the analog object on $X'(k)$.
The following results are the analogue of \cite[Proposition 
4.14]{Shen_2016}.  
\begin{proposition}\label{prop:3.6.5nov19}
	Assume that $G$ has noncompact center with $\delta(G)=1$ and that 
	$(E,\rho^{E})$ is irreducible. Then 
	\begin{equation}
		\mathrm{Tr_{s}}^{[1]}[(N^{\Lambda^{\bullet}(T^{*}X)}-\frac{m}{2})\exp(-t\mathbf{D}^{X,F,2}/2)]=-\frac{e^{-t\alpha^{2}_{E}/2}}{\sqrt{2\pi t}}[e(TX',\nabla^{TX'})]^{\mathrm{max}}\dim E.
		\label{eq:3.6.21novend}
	\end{equation}
	
	If $\gamma=e^{a}k$ is such that $a\in\kb$, 
	$k\in T$, then 
	\begin{equation}
		\begin{split}
			&\mathrm{Tr_{s}}^{[\gamma]}[(N^{\Lambda^{\bullet}(T^{*}X)}-\frac{m}{2})\exp(-t\mathbf{D}^{X,F,2}/2)]\\
			&=-\frac{1}{\sqrt{2\pi 
			t}}e^{-\frac{|a|^{2}}{2t}-t\alpha^{2}_{E}/2}[e(TX^{\prime}(k),\nabla^{TX^{\prime}(k)})]^{\mathrm{max}(k)}\mathrm{Tr}^{E}[\rho^{E}(k)].
		\end{split}
		\label{eq:3.6.22novend}
	\end{equation}
\end{proposition}
\begin{proof}
	Let $C^{\ku'}$ denote the Casimir operator of $\ku'$ associated 
	with $B|_{\ku'}$. Then we have
	\begin{equation}
		C^{\ku} = -a_{1}^{2} + C^{\ku'}.
		\label{eq:4.1.34std}
	\end{equation}
	Since $(E,\rho^{E})$ is an irreducible representation, by 
	\eqref{eq:4.1.28std} and \eqref{eq:4.1.34std}, we get
	\begin{equation}
		C^{\ku,E} = -\alpha_{E}^{2} + C^{\ku',E}.
		\label{eq:4.1.35std}
	\end{equation}
	
	Then by \eqref{eq:4.1.35std} and \cite[Theorem 
	8.5]{BMZ2015toeplitz}, a modification of the 
	proof to 
	\cite[Proposition 4.14]{Shen_2016} proves the identities in our 
	proposition.
\end{proof}

If we assembly the results in Proposition \ref{prop:3.6.5nov19}, it 
is enough to study the corresponding analytic torsions. We will get 
back to this point in Corollary \ref{cor:7.3.6kkss} for asymptotic 
analytic torsions.

%%%%%%%%%%%%%%%%%%%%%%%%%%%%%%%%%%%%%%%%%%%%%%%%%%%%%%%%%%%%%%%%%
\subsection{Symmetric spaces of noncompact type with fundamental rank 
$1$}\label{subs:4.2std}

In this subsection, we focus on the case where $\delta(G)=1$ 
and $G$ has compact center (i.e. $\z_{\pp}=0$), so that $X$ is a 
symmetric space of noncompact type \cite[Proposition 
6.18]{Shen_2016}. For simplicity, let us also assume that $G$ is 
linear semisimple in this subsection.

Note that the rank $\delta(X)$ of 
$X$ (see \cite[Section 2.7]{eberlein1996geometry}) is the 
same as $\delta(G)$, then $\delta(X)=1$. By the de Rham decomposition, we can write 
\begin{equation}
	X=X_{1}\times X_{2},
	\label{eq:4.0.33kk}
\end{equation}
where $X_{1}$ is an irreducible symmetric space 
of noncompact type with $\delta(X_{1})=1$, and $X_{2}$ is a 
symmetric space of noncompact type with $\delta(X_{2})=0$.

As in \cite[Remark 
7.9.2]{bismut2011hypoelliptic}, among the noncompact simple connected real linear groups such that 
$m$ is odd and $\dim \kb=1$, there are only $\mathrm{SL}_{3}(\R)$, 
$\mathrm{SL}_{4}(\R)$,  
$\mathrm{SL}_{2}(\mathbb{H})$, and $\mathrm{SO}^{0}(p,q)$ with $pq$ 
odd $>1$. Also, we have $\mathfrak{sl}_{4}(\R)=\mathfrak{so}(3,3)$ 
and $\mathfrak{sl}_{2}(\mathbb{H})=\mathfrak{so}(5,1)$. Therefore, $X_{1}$ is 
one of the following cases (see \cite[Proposition 6.19]{Shen_2016})
\begin{equation}
	X_{1}=\mathrm{SL}_{3}(\R)/\mathrm{SO}(3)\;\mathrm{or}\; 
	\mathrm{SO}^{0}(p,q)/\mathrm{SO}(p+q) \; \mathrm{with}\; 
	pq>1\;\mathrm{odd}\;.
	\label{eq:4.2.5std20}
\end{equation}

Since $\delta(G)=1$, we have the following decomposition of Lie 
algebras,
\begin{equation}
	\g=\g_{1}\oplus\g_{2},
	\label{eq:7.4.10kk20}
\end{equation}
where
\begin{equation}
	\g_{1}=\mathfrak{sl}_{3}(\R) \;\mathrm{or}\; \mathfrak{so}(p,q)
\end{equation}
with $pq>1$ odd, and $\g_{2}$ is semisimple with 
$\delta(\g_{2})=0$. 
The Cartan involution $\theta$ preserves the 
splitting \eqref{eq:7.4.10kk20} (see \cite[VII.6, p. 471]{knapp2002liegroupe}).

Let $G_{1}$ be the identity component of $Z_{G}(\g_{2})$, then 
$G_{1}$ is a connected linear semisimple closed subgroup of $G$ with Lie 
algebra of $\g_{1}$. Similarly, we can find a connected linear 
semisimple closed subgroup $G_{2}$ of $G$ with Lie 
algebra of $\g_{2}$ such that we have canonically $G_{1}\times 
G_{2}\rightarrow G$ a finite central extension. Let $\theta_{j}$ be 
the induced Cartan involution on $G_{j}$($j=1,2$) from $\theta$, set 
$K_{j}=G_{j}\cap K$, then
\begin{equation}
	X_{j}=G_{j}/K_{j}, j=1,2.
\end{equation}
Note that in general, $G_{1}$ is a just a finite central extension of 
$\mathrm{SL}_{3}(\R)$ or $\mathrm{SO}^{0}(p,q)$ ($pq>1$ odd). The 
invariant bilinear form $B$ also splits as $B_{1}\oplus B_{2}$ with 
respect to the splitting \eqref{eq:7.4.10kk20}.

\begin{remark}
	Let $G_{*}$, $G_{1,*}$
	, $G_{2,*}$ denote the identity components of the isometry groups of 
	$X$, $X_{1}$, $X_{2}$ respectively. Then we have
	\begin{equation}
		G_{*}=G_{1,*}\times G_{2,*}.
	\end{equation}
	By \cite[Proposition 6.19]{Shen_2016}, $G_{1,*}= \mathrm{SL}_{3}(\R)$ or 
	$\mathrm{SO}^{0}(p,q)$ with $pq>1$ odd, and $G_{2,*}$ is a semisimple Lie group 
	with Lie algebra $\g_{2}$ and trivial center. Also 
	$\delta(G_{2,*})=0$. If we consider $G_{*}$ instead of 
	$G$, then the factor $G_{1}$ is exactly $\mathrm{SL}_{3}(\R)$ or 
	$\mathrm{SO}^{0}(p,q)$ with $pq>1$ odd.
\end{remark}

Let $U_{1}$, $U_{2}$ be (connected linear) compact forms of $G_{1}$, 
$G_{2}$. Then $U_{1}\times U_{2}$ is a finite central extension of 
the compact form $U$ of $G$. Let $(E,\rho^{E})$ be an 
irreducible unitary representation of $U$, hence of $U_{1}\times U_{2}$, then 
\begin{equation}
	(E,\rho^{E})=(E_{1},\rho^{E_{1}})\otimes (E_{2},\rho^{E_{2}}),
\end{equation}
where $(E_{j},\rho^{E_{j}})$ is an irreducible unitary representation 
of $U_{j}$, $j=1,2$. Let $F$, $F_{1}$, $F_{2}$ be the homogeneous flat 
vector bundles on $X$, $X_{1}$, $X_{2}$ associated with these 
representations. Then we have
\begin{equation}
	F=F_{1}\boxtimes F_{2}:=\pi_{1}^{\ast}(F_{1})\otimes 
	\pi_{2}^{\ast}(F_{2}),
	\label{eq:4.0.39}
\end{equation}
where $\pi_{i}$ denote the projection $X\rightarrow X_{i}$, $i=1,2$.

Take $\gamma\in G$, let $(\gamma_{1},\gamma_{2})\in G_{1}\times 
G_{2}$ be one of its lifts. Then 
$\gamma$ is semisimple (resp. elliptic) if and only if both $\gamma_{1}$, $\gamma_{2}$ are 
semisimple (resp. elliptic). Set $m_{i}=\dim X_{i}$, then $m_{1}$ is 
odd, and $m_{2}$ is even.

\begin{proposition}\label{prop:4.2.2bio}
	If $\gamma\in G$ is semisimple, for 
	$t>0$, we have
	\begin{equation}
		\begin{split}
			&\mathrm{Tr_s}^{[\gamma]}\big[(N^{\Lambda^{\bullet}(T^*X)}-\frac{m}{2})\exp(-t\mathbf{D}^{X,F,2}/2)\big]\\
			&=	
			\mathrm{Tr_s}^{[\gamma_{1}]}\big[(N^{\Lambda^{\bullet}(T^*X_{1})}-\frac{m_{1}}{2})\exp(-t\mathbf{D}^{X_{1},F_{1},2}/2)\big]\cdot\mathrm{Tr_s}^{[\gamma_{2}]}\big[\exp(-t\mathbf{D}^{X_{2},F_{2},2}/2)\big]
		\end{split}
		\label{eq:4.2.8std}
	\end{equation}
	Then if $\gamma_{2}$ is nonelliptic, 
	\begin{equation}
		\mathrm{Tr_s}^{[\gamma]}\big[(N^{\Lambda^{\bullet}(T^*X)}-\frac{m}{2})\exp(-t\mathbf{D}^{X,F,2}/2)\big]=0.
		\label{eq:4.2.11bio}
	\end{equation}
	If $\gamma_{2}$ is elliptic, then
	\begin{equation}
		\begin{split}
			&\mathrm{Tr_s}^{[\gamma]}\big[(N^{\Lambda^{\bullet}(T^*X)}-\frac{m}{2})\exp(-t\mathbf{D}^{X,F,2}/2)\big]\\
			&=	
			[e(TX_{2}(\gamma_{2}),\nabla^{TX_{2}(\gamma_{2})})]^{\mathrm{max}_{2}(\gamma_{2})} \mathrm{Tr}^{E_{2}}[\rho^{E_{2}}(\gamma_{2})]\\
			&\qquad\cdot
			\mathrm{Tr_s}^{[\gamma_{1}]}\big[(N^{\Lambda^{\bullet}(T^*X_{1})}-\frac{m_{1}}{2})\exp(-t\mathbf{D}^{X_{1},F_{1},2}/2)\big],
		\end{split}
		\label{eq:4.2.12bio}
	\end{equation}
	where $[\cdot]^{\mathrm{max}_{2}(\gamma_{2})}$ is taking the 
	coefficient of the Riemannian volume element on $X_{2}(\gamma_{2})$.
\end{proposition}
\begin{proof}
	We write
	\begin{equation}
		N^{\Lambda^{\bullet}(T^*X)}-\frac{m}{2}=(N^{\Lambda^{\bullet}(T^*X_{1})}-\frac{m_{1}}{2})+ (N^{\Lambda^{\bullet}(T^*X_{2})}-\frac{m_{2}}{2}).
		\label{eq:4.2.11std}
	\end{equation}
	Note that, since $\delta(G_{1})=1$, then by \cite[Theorem 
	7.8.2]{bismut2011hypoelliptic}, we always have
	\begin{equation}
		\mathrm{Tr_s}^{[\gamma_{1}]}\big[\exp(-t\mathbf{D}^{X_{1},F_{1},2}/2)\big]=0.
		\label{eq:4.2.12std}
	\end{equation}
	Combining the definition of 
	orbital integrals \eqref{eq:3.3.8ffe} together with \eqref{eq:4.2.11std} and \eqref{eq:4.2.12std}, we 
	get \eqref{eq:4.2.8std}.
	
	The identities \eqref{eq:4.2.11bio}, \eqref{eq:4.2.12bio} follow from 
	applying the results in \cite[Theorem 
	7.8.2]{bismut2011hypoelliptic} to 
	$\mathrm{Tr_s}^{[\gamma_{2}]}\big[\exp(-t\mathbf{D}^{X_{2},F_{2},2}/2)\big]$. This completes the proof of our proposition.
\end{proof}

For studying $\mathcal{T}(Z,F)$, Proposition \ref{prop:4.2.2bio} helps us to reduce the 
computations on 
$\mathrm{Tr_{s}}[(N^{\Lambda^{\bullet}(T^{\ast}Z)}-\frac{m}{2})\exp(-t\mathbf{D}^{Z,F,2}/2)]$ to the model cases listed in \eqref{eq:4.2.5std20}. But it is far from enough to get an explicit evaluation. In Sections \ref{section4} \& \ref{section6paris}, we will carry out more tools, which allows us work out a proof to Theorem \ref{thm:maintheorem}.

%%%%%%%%%%%%%%%%%%%%%%%%%%%%%%%%%%%%%%%%%%%%%%%%%%%%%%%%%%%%%%%%%%%%
\section{Cartan subalgebra and root system of $G$ when $\delta(G)=1$}\label{section4}
We use the same notation as in Section \ref{section2} \& Subsection 
\ref{subsection4.1paris}. 
In Subsections \ref{section5.1ss20} - \ref{section5.3pl}, we always 
assume that $G$ is connected linear real reductive Lie group with 
compact center and with $\delta(G)=1$. But, as we 
will see in Remark \ref{rm:important}, the constructions and 
results in these subsections are still true (most of them are 
trivial) if $U$ is compact and if $G$ has noncompact center with 
$\delta(G)=1$.

Subsection \ref{section5.4} is independent from other 
subsections, where we introduce a generalized Kirillov formula for 
compact Lie groups.

Recall that $T$ is a 
maximal torus of $K$ with Lie algebra $\kt\subset \kk$, and that $\kb\subset 
\pp$ is defined in \eqref{eq:defofb}. Since $\delta(G)=1$, then $\kb$ 
is $1$-dimensional. We now fix a vector $a_{1}\in \kb$, $|a_{1}|=1$. Recall 
that $\kh=\kb\oplus\kt$ is a Cartan subalgebra of $\g$. Let 
$h^{\g_{\C}}$ be the Hermitian product  on $\g_\C$ induced by
the scalar product $-B(\cdot,\theta\cdot)$ on $\g$.
%%%%%%%%%%%%%%%%%%%%%%%%%%%%%%%%%%%%%%%%%%%%%%%%%%%%%%%%%%%%%%%%%%%%
%%%%%%%%%%%%%%%%%%%%%%%%%%%%%%%%%%%%%%%%%%%%%%%%%%%%%%%%%%%%%%%%%%%%%
\subsection{Reductive Lie algebra with fundamental rank $1$}\label{section5.1ss20}
Since $G$ has compact center, then $\kb\not\subset \z_{\g}$. Let $Z(\kb)$ be the centralizer of $\kb$ in $G$, and let $Z(\kb)^{0}$ 
be its identity component with Lie algebra 
$\z(\kb)=\pp(\kb)\oplus\kk(\kb)\subset\g$. Let $\km$ be the 
orthogonal subspace of $\kb$ in $\z(\kb)$ (with respect to $B$) such that
\begin{equation}
	\z(\kb)=\kb\oplus\km.
\end{equation}
Then $\km$ is a Lie subalgebra of $\z(\kb)$, which is invariant by 
$\theta$.

Put
\begin{equation}
	\pp_{\km}=\km\cap \pp,\; \kk_{\km}=\km\cap\kk.
\end{equation}
Then 
\begin{equation}
	\km=\pp_{\km}\oplus\kk_{\km}, \pp(\kb)=\kb\oplus\pp_{\km},\; 
	\kk(\kb)=\kk_{\km}.
	\label{eq:5.3.3paris19}
\end{equation}

Let $\z^{\perp}(\kb)$, $\pp^{\perp}(\kb)$, $\kk^{\perp}(\kb)$ be the 
orthogonal subspaces of $\z(\kb)$, $\pp(\kb)$, $\kk(\kb)$ in $\g$, 
$\pp$, $\kk$ respectively with respect to $B$. Then
\begin{equation}
	\z^{\perp}(\kb)=\pp^{\perp}(\kb)\oplus \kk^{\perp}(\kb).
	\label{eq:5.3.4paris19}
\end{equation}
Moreover,
\begin{equation}
	\pp=\kb\oplus\pp_{\km}\oplus\pp^{\perp}(\kb),\; 
	\kk=\kk(\kb)\oplus\kk^{\perp}(\kb).
	\label{eq:5.3.5paris19}
\end{equation}

Let $M\subset Z(\kb)^{0}$ be the analytic subgroup associated with 
$\km$. If we identify $\kb$ with $\R$, then
\begin{equation}
	Z(\kb)^{0}=\R\times M.
\end{equation}
Then $M$ is a Lie subgroup of $Z(\kb)^{0}$, i.e., it is closed in 
$Z(\kb)^{0}$. Let $K_{M}$ be the analytic subgroup of $M$ associated 
with the Lie subalgebra $\kk_{\km}$. Since $M$ is reductive, $K_{M}$ 
is a maximal compact subgroup of $M$. Then the splittings in 
\eqref{eq:5.3.3paris19}, \eqref{eq:5.3.4paris19}, 
\eqref{eq:5.3.5paris19} are invariant by the adjoint action of 
$K_{M}$.

Then $\kt$ is Cartan subalgebra of $\kk$, of $\kk_{\km}$, and of 
$\km$. Recall that $\kh=\kb\oplus\kt$ is a Cartan subalgebra of $\g$.
We fix $a_{1}\in\kb$ such that $B(a_{1},a_{1})=1$. The choice of 
$a_{1}$ fixes an orientation of $\kb$. Let $\kn\subset 
\z^{\perp}(\kb)$ be the direct sum of the eigenspaces of 
$\mathrm{ad}(a_{1})$ with the positive eigenvalues. Set 
$\bar{\kn}=\theta \kn$. Then
\begin{equation}
	\z^{\perp}(\kb)=\kn\oplus\bar{\kn}.
\end{equation}
By \cite[Subsection 6A]{Shen_2016}, $\dim \kn=\dim\pp 
-\dim\pp_{\km}-1$. Then $\dim \kn$ is even under our assumption 
$\delta(G)=1$. Put
\begin{equation}
	l=\frac{1}{2}\dim\kn.
\end{equation}

By \cite[Proposition 6.2]{Shen_2016}, there exists $\beta\in\kb^{*}$ 
such that if $a\in\kb$, $f\in\kn$, then
\begin{equation}
	[a,f]=\beta(a)f,\; [a,\theta(f)]=-\beta(a)\theta(f).
	\label{eq:5.1.9dec19}
\end{equation}

The map $f\in\kn\mapsto f-\theta(f)\in\pp^{\perp}(\kb)$ is an 
isomorphism of $K_{M}$-modules. Similarly, $f\in\kn\mapsto 
f+\theta(f)\in\kk^{\perp}(\kb)$ is also an isomorphism of 
$K_{M}$-modules. Since $\theta$ fixes $K_{M}$, $\kn\simeq \bar{\kn}$ 
as $K_{M}$-modules via $\theta$.

By \cite[Proposition 6.3]{Shen_2016}, we have
\begin{equation}
	[\kn, \bar{\kn}]\subset \z(\kb),\; 
	[\kn,\kn]=[\bar{\kn},\bar{\kn}]=0.
	\label{eq:5.1.10paris19}
\end{equation}
Also 
\begin{equation}
	B|_{\kn\times\kn}=0,\; B|_{\bar{\kn}\times\bar{\kn}}=0.
	\label{eq:5.1.11paris19}
\end{equation}
Then the bilinear form $B$ induces an isomorphism of $\kn^{*}$ and 
$\bar{\kn}$ as $K_{M}$-modules. Therefore, as $K_{M}$-modules, $\kn$ 
is isomorphic to $\kn^{*}$.

As a consequence of \eqref{eq:5.1.10paris19}, we get
\begin{equation}
	[\z(\kb),\z(\kb)]\;,[\z^{\perp}(\kb),\z^{\perp}(\kb)]\subset 
	\z(\kb),\; [\z(\kb),\z^{\perp}(\kb)]\subset\z^{\perp}(\kb).
	\label{eq:5.1.12paris19}
\end{equation}
Then $(\g,\z(\kb))$ is a symmetric pair.

If $k\in K_{M}$, let $M(k)$ be the centralizer of $k$ in $M$, and let 
$\km(k)$ be its Lie algebra. Let $M(k)^{0}$ be the identity component 
of $M(k)$. The Cartan involution $\theta$ acts on $M(k)$. The 
associated Cartan decomposition is
\begin{equation}
	\km(k)=\pp_{\km}(k)\oplus\kk_{\km}(k),
	\label{eq:5.1.13paris19}
\end{equation}
where $\pp_{\km}(k)=\pp_{\km}\cap \km(k)$, 
$\kk_{\km}(k)=\kk_{\km}\cap \km(k)$.

Recall that $Z(k)$ is the centralizer 
of $k$ in $G$ and that $Z(k)^{0}$ is the identity component of $Z(k)$ with Lie algebra $\z(k)\subset \g$.  Then
\begin{equation}
	M(k)=M\cap Z(k),\; \km(k)=\km \cap \z(k).
	\label{eq:5.1.14ss20k}
\end{equation}

Note that $Z(k)^{0}$ is still a reductive Lie group equipped with the 
Cartan involution induced by the action of $\theta$. By the 
assumption that $\delta(G)=1$, we have
\begin{equation}
	\delta(Z(k)^{0})=1.
\end{equation}
In particular,
\begin{equation}
	\kb\subset \pp(k).
	\label{eq:ss201301}
\end{equation}

Set
\begin{equation}
	\z_{\kb}(k)=\z(\kb)\cap \z(k),\; \pp_{\kb}(k)=\pp(\kb)\cap 
	\pp(k),\; \kk_{\kb}(k)=\kk(\kb)\cap \kk(k).
	\label{eq:5.1.17ss20m}
\end{equation}
Then
\begin{equation}
	\z_{\kb}(k)=\kb\oplus \km(k)=\pp_{\kb}(k)\oplus \kk_{\kb}(k).
	\label{eq:5.1.18ss20m}
\end{equation}
We also have the following identities,
\begin{equation}
	\pp_{\kb}(k)=\kb\oplus\pp_{\km}(k),\; \kk_{\kb}(k)=\kk_{\km}(k).
	\label{eq:5.1.19ss20m}
\end{equation}

Let $\pp_{\kb}^{\perp}(k)$, $\kk_{\kb}^{\perp}(k)$, 
$\z^{\perp}_{\kb}(k)$ be the orthogonal 
spaces of $\pp_{\kb}(k)$, $\kk_{\kb}(k)$, $\z_{\kb}(k)$ in $\pp(k)$, 
$\kk(k)$, $\z(k)$ with 
respect to $B$, so that
\begin{equation}
	\pp(k)=\pp_{\kb}(k)\oplus \pp_{\kb}^{\perp}(k),\; 
	\kk(k)=\kk_{\kb}(k)\oplus \kk_{\kb}^{\perp}(k),\; 
	\z(k)=\z_{\kb}(k)\oplus \z_{\kb}^{\perp}(k).
\end{equation}
Then
\begin{equation}
	\z_{\kb}^{\perp}(k)=\pp_{\kb}^{\perp}(k) \oplus 
	\kk_{\kb}^{\perp}(k)= \z^{\perp}(\kb)\cap\z(k).
\end{equation}

Put
\begin{equation}
	\kn(k)=\z(k)\cap \kn, \;\bar{\kn}(k)=\z(k)\cap\bar{\kn}.
	\label{eq:5.1.14ssbis}
\end{equation}
Then
\begin{equation}
	\z_{\kb}^{\perp}(k)=\kn(k)\oplus \bar{\kn}(k).
	\label{eq:5.1.21ss20m}
\end{equation}

By \eqref{eq:5.1.17ss20m}, \eqref{eq:5.1.21ss20m}, we get
\begin{equation}
	\z(k)=\pp_{\kb}(k)\oplus\kk_{\kb}(k)\oplus 
	\kn(k)\oplus\bar{\kn}(k).
\end{equation}
Since $\delta(\km(k))=0$, $\dim \kn(k)$ is even. We set
\begin{equation}
	l(k)=\frac{1}{2}\dim \kn(k).
	\label{eq:5.1.24ss20m}
\end{equation}

Let $K_{M}(k)$ denote the centralizer of $k$ in $K_{M}$.
The map $f\in\kn(k)\mapsto f-\theta(f)\in\pp_{\kb}^{\perp}(k)$ is an 
isomorphism of $K_{M}(k)$-modules, similar for 
$\kk_{\kb}^{\perp}(k)$. Since $\theta$ fixes 
$K_{M}(k)$, 
$\kn(k)\simeq \bar{\kn}(k)$ 
as $K_{M}(k)$-modules via $\theta$.

%%%%%%%%%%%%%%%%%%%%%%%%%%%%%%%%%%%%%%%%%%%%%%%%%%%%%%%%%%%%%%%%%%%%%%%
\subsection{A compact Hermitian symmetric space 
$Y_{\kb}$}\label{subsection5.2ss}
Recall that $\ku=\sqrt{-1}\pp\oplus\kk$ is the compact form of $\g$. 

Let $\ku(\kb)\subset \ku$, $\ku_{\km}\subset \ku$ be the compact 
forms of $\z(\kb)$, $\km$. Then
\begin{equation}
	\ku(\kb)=\sqrt{-1}\kb\oplus \ku_{\km},\; 
	\ku_{\km}=\sqrt{-1}\pp_{\km}\oplus\kk_{\km}.
	\label{eq:5.1.17nov19ss}
\end{equation}

Since $M$ has compact center, let $U_{M}$ be the analytic subgroup of 
$U$ associated with $\ku_{\km}$. Then $U_{M}$ is the compact form of 
$M$. Let $U(\kb)\subset U$, $A_{0}\subset U$ be the connected 
subgroups of $U$ associated with Lie algebras $\ku(\kb)$, 
$\sqrt{-1}\kb$. Then $A_{0}$ is in the center of $U(\kb)$. By 
\cite[Proposition 6.6]{Shen_2016}, $A_{0}$ is closed in $U$ and is 
diffeomorphic to a circle $\mathbb{S}^{1}$. Moreover, we have
\begin{equation}
	U(\kb)=A_{0}U_{M}.
	\label{eq:5.1.18dec19}
\end{equation}

The bilinear form $-B$ induces an $\mathrm{Ad}(U)$-invariant metric 
on $\ku$. Let $\ku^{\perp}(\kb)\subset\ku$ be the orthogonal subspace 
of $\ku(\kb)$. Then
\begin{equation}
	\ku^{\perp}(\kb)=\sqrt{-1}\pp^{\perp}(\kb)\oplus\kk^{\perp}(\kb).
	\label{eq:5.2.3dec19}
\end{equation}
By \eqref{eq:5.1.12paris19}, we get
\begin{equation}
	[\ku(\kb),\ku(\kb)]\;,[\ku^{\perp}(\kb),\ku^{\perp}(\kb)]\subset 
	\ku(\kb),\; [\ku(\kb),\ku^{\perp}(\kb)]\subset\ku^{\perp}(\kb).
	\label{eq:5.2.4dec19}
\end{equation}
Then $(\ku,\ku(\kb))$ is a symmetric pair.

Put $a_{0}=a_{1}/\beta(a_{1})\in\kb$. Set
\begin{equation}
	J=\sqrt{-1}\mathrm{ad}(a_{0})|_{\ku^{\perp}(\kb)}\in\mathrm{End}(\ku^{\perp}(\kb)).
	\label{eq:5.2.5dec19}
\end{equation}
By \eqref{eq:5.1.9dec19}, $J$ is an $U(\kb)$-invariant complex 
structure on $\ku^{\perp}(\kb)$ which preserves 
$B|_{\ku^{\perp}(\kb)}$. The spaces $\kn_{\C}=\kn\otimes_{\R}\C$, 
$\bar{\kn}_{\C}=\bar{\kn}\otimes_{\R}\C$ are exactly the eigenspaces 
of $J$ associated with eigenvalues $\sqrt{-1}$, $-\sqrt{-1}$.

The following proposition is just the summary of the results in 
\cite[Section 6B]{Shen_2016}.
\begin{proposition}\label{prop:5.2.1dec19}
	Set
	\begin{equation}
		Y_{\kb}=U/U(\kb).
		\label{eq:5.2.6dec19}
	\end{equation}
	Then $Y_{\kb}$ is a compact symmetric space, and $J$ induces an integrable complex structure on $Y_{\kb}$ 
	such that
	\begin{equation}
		T^{(1,0)}Y_{\kb}=U\times_{U(\kb)}\kn_{\C},\; 
		T^{(0,1)}Y_{\kb}=U\times_{U(\kb)}\bar{\kn}_{\C}.
		\label{eq:5.2.7dec19}
	\end{equation}
	The form $-B(\cdot, J\cdot)$ induces a K\"{a}hler form 
	$\omega^{Y_{\kb}}$ on $Y_{\kb}$. 
\end{proposition}

Let $\omega^{\ku}$ be the canonical left-invariant $1$-form on $U$ 
with values in $\ku$. Let $\omega^{\ku(\kb)}$ and
$\omega^{\ku^{\perp}(\kb)}$ be the $\ku(\kb)$ and $\ku^{\perp}(\kb)$ 
components of $\omega^{\ku}$, so that
\begin{equation}
	\omega^{\ku}=\omega^{\ku(\kb)}+\omega^{\ku^{\perp}(\kb)}.
	\label{eq:5.2.8dec19}
\end{equation}
Moreover, $\omega^{\ku(\kb)}$ defines a connection form on the principal 
$U(\kb)$-bundle $U\rightarrow Y_{\kb}$. Let $\Omega^{\ku(\kb)}$ be 
the curvature form, then
\begin{equation}
	\Omega^{\ku(\kb)}=-\frac{1}{2}[\omega^{\ku^{\perp}(\kb)},\omega^{\ku^{\perp}(\kb)}].
	\label{eq:5.2.9dec19}
\end{equation}

Note that the real tangent bundle of $Y_{\kb}$ is given
\begin{equation}
	TY_{\kb}=U\times_{U(\kb)}\ku^{\perp}(\kb).
	\label{eq:5.2.7decbis}
\end{equation}
Then $-B|_{\ku^{\perp}(\kb)}$ induces a Riemannian metric 
$g^{TY_{\kb}}$ on $Y_{\kb}$.  The 
corresponding Levi-Civita 
connection is induced by $\omega^{\ku(\kb)}$. 

Recall that the first splitting in \eqref{eq:5.1.17nov19ss} is 
orthogonal with respect to $-B$. Let $\Omega^{\ku_{\km}}$ be the 
$\ku_{\km}$-component of $\Omega^{\ku(\kb)}$. 
Since the K\"{a}hler form $\omega^{Y_{\kb}}$ is invariant  under 
the left action of $U$ on $Y_{\kb}$, we also can view 
$\omega^{Y_{\kb}}$ as an element in 
$\Lambda^{2}(\ku_{\kb}^{\perp})^{*})$.
By 
\cite[Eq.(6-48)]{Shen_2016},
\begin{equation}
	\Omega^{\ku(\kb)}=\beta(a_{1})\omega^{Y_{\kb}}\otimes\sqrt{-1}a_{1}+\Omega^{\ku_{\km}}.
	\label{eq:5.2.11dec19}
\end{equation}
Moreover, by \cite[Proposition 6.9]{Shen_2016}, we have
\begin{equation}
	B(\Omega^{\ku(\kb)}, \Omega^{\ku(\kb)})=0,\; 
	B(\Omega^{\ku_{\km}},\Omega^{\ku_{\km}})=\beta(a_{1})^{2}\omega^{Y_{\kb},2}.
	\label{eq:5.2.12dec19}
\end{equation}

\begin{remark}
	By \cite[Proposition 6.20]{Shen_2016}, if $G$ has compact center, 
	then as symmetric spaces, the K\"{a}hler manifold 
	$Y_{\kb}$ is isomorphic either to $\mathrm{SU}(3)/\mathrm{U}(2)$ or to 
	$\mathrm{SO}(p+q)/\mathrm{SO}(p+q-2)\times \mathrm{SO}(2)$ with $pq>1$ odd. This way, the 
	computations on $Y_{\kb}$ can be made more explicit.	
\end{remark}

Now we fix $k\in K_{M}$. Let $U(k)$ be the centralizer of $k$ in $U$, 
and let $U(k)^{0}$ be its identity component. Let $\ku(k)$ be the Lie 
algebra of $U(k)^{0}$. Then $\ku(k)$ is the compact form of $\z(k)$, and 
$U(k)^{0}$ is the compact form of $Z(k)^{0}$.

We will use the same notation as in Subsection 
\ref{section5.1ss20}. Then the compact form of $\km(k)$ is given by
\begin{equation}
	\ku_{\km}(k)=\sqrt{-1}\pp_{\km}(k)\oplus \kk_{\km}(k).
	\label{eq:5.1.13ssss20s}
\end{equation}
Let $\ku_{\kb}(k)$ be the compact form of $\z_{\kb}(k)$. Then
\begin{equation}
	\ku_{\kb}(k)=\sqrt{-1}\kb\oplus\ku_{\km}(k).
	\label{eq:5.2.14ss20s}
\end{equation}
Let $U_\kb(k)$ be the analytic subgroup associated with 
$\ku_{\kb}(k)$. Then
\begin{equation}
	U_{\kb}(k)=U(\kb)\cap U(k)^{0}.
\end{equation}

Set
\begin{equation}
	Y_{\kb}(k)=U(k)^{0}/ U_{\kb}(k).
\end{equation}
As in Proposition \ref{prop:5.2.1dec19}, $Y_{\kb}(k)$ is a connected 
complex manifold equipped with a K\"{a}hler form 
$\omega^{Y_{\kb}(k)}$. 

Let $\ku^{\perp}_{\kb}(k)$ be the orthogonal space of $\ku_{\kb}(k)$ 
in $\ku(k)$ with respect to $B$. Then
\begin{equation}
	\ku_{\kb}^{\perp}(k)=\sqrt{-1}\pp_{\kb}^{\perp}(k)\oplus 
	\kk_{\kb}^{\perp}(k).
\end{equation}
Then the real tangent bundle of 
$Y_{\kb}(k)$ is given by
\begin{equation}
	TY_{\kb}(k)=U(k)^{0}\times_{U_{\kb}(k)}\ku_{\kb}^{\perp}(k).
	\label{eq:5.2.10ss20}
\end{equation}
Moreover,
\begin{equation}
	T^{(1,0)}Y_{\kb}(k)=U(k)^{0}\times_{U_\kb(k)}\kn(k)_{\C},\; 
	T^{(0,1)}Y_{\kb}(k)=U(k)^{0}\times_{U_\kb(k)}\bar{\kn}(k)_{\C}.
	\label{eq:5.2.11ss20}
\end{equation}

Let $\Omega^{\ku_{\kb}(k)}$ be the curvature form as in 
\eqref{eq:5.2.9dec19} for the pair $(U(k)^{0}, U_{\kb}(k))$, which can be 
viewed as an element in $\Lambda^{2}(\ku_{\kb}^{\perp}(k)^{*})\otimes 
\ku_{\kb}(k)$. Using the splitting \eqref{eq:5.2.14ss20s}, let 
$\Omega^{\ku_{\km}(k)}$ be the $\ku_{\km}(k)$-component of 
$\Omega^{\ku_{\kb}(k)}$. Then as in \eqref{eq:5.2.11dec19} and 
\eqref{eq:5.2.12dec19}, we have
\begin{equation}
	\Omega^{\ku_{\kb}(k)}=\beta(a_{1})\omega^{Y_{\kb}(k)}\otimes\sqrt{-1}a_{1}+\Omega^{\ku_{\km}(k)},
	\label{eq:5.2.20dec19}
\end{equation}
and
\begin{equation}
	B(\Omega^{\ku_{\kb}(k)}, \Omega^{\ku_{\kb}(k)})=0,\; 
	B(\Omega^{\ku_{\km}(k)},\Omega^{\ku_{\km}(k)})=\beta(a_{1})^{2}\omega^{Y_{\kb}(b),2}.
	\label{eq:5.2.21dec19}
\end{equation}

%%%%%%%%%%%%%%%%%%%%%%%%%%%%%%%%%%%%%%%%%%%%%%%%%%%%%%%%%%%%%%%%%%%%%
\subsection{Positive root system and character formula}\label{section5.3pl}
Recall that $\kt$ is Cartan subalgebra of $\kk$, of $\kk_{\km}$, and of 
$\km$. Recall that $\kh=\kb\oplus\kt$ is a Cartan subalgebra of $\g$, 
and $H$ is the associated maximally compact Cartan subgroup of $G$.

Put 
\begin{equation}
	\kt_{U}=\sqrt{-1}\kb\oplus\kt\subset \ku.
	\label{eq:5.3.1th19}
\end{equation}
Then $\kt_{U}$ is a 
Cartan subalgebra of $\ku$. Let $T_{U}\subset U$ be the corresponding 
maximal torus. Then $A_{0}$ is a circle in $T_{U}$. Then $\kt$ is a 
Cartan subalgebra of $\ku_{\km}$, and the corresponding maximal torus 
is $T$.

Let $R(\ku,\kt_{U})$ be the real root system for the pair 
$(U,T_{U})$ \cite[Chapter V]{TTrepresentation1985}. The root system 
for the complexified pair $(\ku_{\C},\kt_{U,\C})=(\g_{\C},\kh_{\C})$ 
is given by $2\pi i R(\ku,\kt_{U})$. Similarly, let $R(\ku(\kb),\kt_{U})$, $R(\ku_{\km},\kt)$ denote the 
real root systems for the pairs $(\ku(\kb),\kt_{U})$, 
$(\ku_{\km},\kt)$. When we embed $\kt^*$ into $\kt_{U}^{*}$ by the 
splitting in \eqref{eq:5.3.1th19}, then 
\begin{equation}
	R(\ku(\kb),\kt_{U})=R(\ku_{\km},\kt).
	\label{eq:5.3.7dec19s}
\end{equation}

For a root $\alpha\in R(\ku,\kt_{U})$, if $\alpha(\sqrt{-1}a_{1})=0$, then 
$\alpha\in R(\ku_{\km},\kt)$. 
Fix a positive root system $R^{+}(\ku_{\km},\kt)$, we get a positive 
root system $R^{+}(\ku,\kt_{U})$ consisting of element $\alpha$ such 
that $\alpha(\sqrt{-1}a_{1})>0$ and the elements in $R^{+}(\ku_{\km},\kt)$.

Let $W(\ku,\kt_{U})$ denote the algebraic Weyl group associated with 
$R(\ku,\kt_{U})$. If $\omega\in W(\ku,\kt_{U})$, let $l(\omega)$ denote the length of $\omega$ with 
respect to $R^{+}(\ku,\kt_{U})$. Set
\begin{equation}
	\varepsilon(\omega)=(-1)^{l(\omega)}.
\end{equation}
Let $W(U,T_{U})$ be the analytic Weyl group, then $W(\ku,\kt_{U})=W(U,T_{U})$.

Put
\begin{equation}
	W_{u}=\{\omega\in W(U,T_{U})\;:\; \omega^{-1}\cdot \alpha >0, 
	\mathrm{for\;all\;} \alpha\in R^{+}(\ku_{\km}, \kt)\}.
	\label{eq:5.3.6ss20}
\end{equation}

Put
\begin{equation}
	\rho_{\ku}=\frac{1}{2}\sum_{\alpha^{0}\in 
	R^{+}(\ku,\kt_{U})}\alpha^{0}\in\kt_{U}^{*},\; 
	\rho_{\ku_{\km}}=\frac{1}{2}\sum_{\alpha^{0}\in 
	R^{+}(\ku_{\km},\kt)}\alpha^{0}\in\kt^{*}.
\end{equation}
Then $\rho_{\ku}|_{\kt}=\rho_{\ku_\km}$.

Let $P_{++}(U)\subset \kt_{U}^{*}$ be the set of dominant weights of 
$(U,T_{U})$ with respect to $R^{+}(\ku,\kt_{U})$. If $\lambda\in 
P_{++}(U)$, let $(E_{\lambda},\rho^{E_{\lambda}})$ be the irreducible 
unitary representation of $U$ with the highest weight 
$\lambda$, which by the unitary trick extends to an irreducible 
representation of $G$.

By 
\cite[Lemmas 1.1.2.15 \& 2.4.2.1]{MR0498999}, if $\omega\in W_{u}$, 
then $\omega(\lambda+\rho_{\ku})-\rho_{\ku}$ is a dominant weight for 
$R^{+}(\ku(\kb),\kt_{U})$. Let $V_{\lambda,\omega}$ denote the 
representation of $U(\kb)$ with the highest weight 
$\omega(\lambda+\rho_{\ku})-\rho_{\ku}$.

Recall that $U(\kb)$ acts on $\kn_{\C}$.
Let $H^{\cdot}(\kn_{\C},E_{\lambda})$ be the Lie algebra cohomology of 
$\kn_{\C}$ with coefficients in $E_{\lambda}$ (see 
\cite{Kostant1961coh}).
By \cite[Theorem 2.5.1.3]{MR0498999}, for $i=0,\cdots, 2l$, we have 
the identification of $U(\kb)$-modules,
\begin{equation}
	H^{i}(\kn_{\C},E_{\lambda})\simeq \oplus_{\substack{\omega\in W_{u}\\ l(\omega)=i} } 
	V_{\lambda,\omega}. 
	\label{eq:5.3.11pap}
\end{equation}

By \eqref{eq:5.3.11pap} and the Poincar\'{e} duality, we get the 
following identifications as $U(\kb)$-modules,
\begin{equation}
	\oplus_{i=0}^{2l}(-1)^{i}\Lambda^{i}\kn_{\C}^{*}\otimes 
	E_{\lambda}=\oplus_{\omega\in 
	W_{u}}\varepsilon(\omega)V_{\lambda,\omega}.
	\label{eq:5.3.10dec19ss}
\end{equation}
Note that if we apply the unitary trick, the above identification 
also holds as $Z(\kb)^{0}$-modules.

\begin{definition}\label{def:5.3.1kk}
	Let $P_{0}:\kt^{*}_{U}\rightarrow \kt^{*}$ 
	denote the orthogonal projection with respect to 
	$B^{*}|_{\kt_{U}^{*}}$. For $\omega\in W_{u}$, $\lambda\in 
	P_{++}(U)$, put
	\begin{equation}
		\eta_{\omega}(\lambda)=P_{0}\big(\omega(\lambda+\rho_{\ku})-\rho_{\ku}\big)\in \kt^{*}.
		\label{eq:5.3.11final}
	\end{equation}
\end{definition}

Note that 
\begin{equation}
	P_{0}\rho_{\ku}=\rho_{\ku_{\km}}.
	\label{eq:5.3.12decfinal19}
\end{equation}
Then
\begin{equation}
	\eta_{\omega}(\lambda)=P_{0}\big(\omega(\lambda+\rho_{\ku})\big)-\rho_{\ku_\km}.
	\label{eq:5.3.13finalkk}
\end{equation}

\begin{proposition}\label{prop:5.3.2dec19}
	If $\lambda\in 
	P_{++}(U)$, for $\omega\in W_{u}$, then $\eta_{\omega}(\lambda)$ is a dominant weight of $(U_{M}, 
	T)$ with respect to $R^{+}(\ku_{\km},\kt)$. Moreover, the restriction of  the $U(\kb)$-representation $V_{\lambda,\omega}$ to the 
	subgroup $U_{M}$ is irreducible, which has the highest weight $\eta_{\omega}(\lambda)$.
\end{proposition}
\begin{proof}
	Since $\omega(\lambda+\rho_{\ku})-\rho_{\ku}$ is analytically 
	integrable, then $\eta_{\omega}(\lambda)$ is also analytically 
	integrable as a weight associated with $(U_{M},T)$. By  
	\eqref{eq:5.3.7dec19s}  and the corresponding
	identification of positive root systems, we know that 
	$\eta_{\omega}(\lambda)$ is dominant with respect to $R^{+}(\ku_{\km},\kt)$.
	
	Recall that $A_{0}\simeq \mathbb{S}^{1}$ is defined in Subsection 
	\ref{subsection5.2ss}. By \eqref{eq:5.1.18dec19}, we get that $A_{0}$ acts on $V_{\lambda,\omega}$ 
	as scalars given by its character, and then $U_{M}$ act irreducibly on 
	$V_{\lambda,\omega}$, which clearly has the highest weight 
	$\eta_{\omega}(\lambda)$. This completes the proof of our 
	proposition.
\end{proof}

\begin{remark}\label{rm:important}
	In general, $U$ is just the analytic subgroup of $G_{\C}$ with 
	Lie algebra $\ku$. If $U$ is compact but $G$ has noncompact 
	center, i.e., 
	$\z_{\pp}=\kb$, then $\kn=\bar{\kn}=0$, so that $l=0$. Recall 
	that in this case, $G'$, $U'$ are defined in Subsection 
	\ref{subsection4.1paris}. Then
	\begin{equation}
		M=G', U_{M}=U'.
		\label{eq:5.3.14conf}
	\end{equation}
	The compact symmetric space $Y_{\kb}$ now reduces to one 
	point. 
	
	Moreover, in \eqref{eq:5.3.6ss20}, the set $W_{u}=\{1\}$, so 
	that  $V_{\lambda,\omega}$ becomes just $E_{\lambda}$ itself. The 
	identities \eqref{eq:5.3.11pap}, \eqref{eq:5.3.10dec19ss} are 
	trivially true, so is Proposition \ref{prop:5.3.2dec19}.
\end{remark}

%%%%%%%%%%%%%%%%%%%%%%%%%%%%%%%%%%%%%%%%%%%%%%%%%%%%%%%%%%%%%%%%%%%%%%
\subsection{Kirillov character formula for compact Lie 
groups}\label{section5.4}

In this subsection, we recall the Kirillov character formula for 
compact Lie groups. We only use the group $U_{M}$ as an explanatory example. 
We fix the maximal torus $T$ and the positive (real) root system 
$R^{+}(\ku_{\km},\kt)$. 

Let $\lambda\in\kt^{*}$ be a dominant (analytically integrable) 
weight of $U_{M}$ with respect to the above root system. Let 
$(V_{\lambda},\rho^{V_\lambda})$ be the irreducible unitary 
representation of $U_{M}$ with the highest weight $\lambda$.

Put
\begin{equation}
	\mathcal{O}=\mathrm{Ad}^{*}(U_{M})(\lambda+\rho_{\ku_{\km}})\subset \ku_{\km}^{*}.
	\label{eq:5.4.1dec19s}
\end{equation}
Then $\mathcal{O}$ is an even-dimensional 
closed manifold.

Since $\lambda+\rho_{\ku_{\km}}$ is regular, then we have the 
following identifications of $U_{M}$-manifolds,
\begin{equation}
	\mathcal{O}\simeq U_{M}/T.
	\label{eq:5.4.2dec19s}
\end{equation}

For $u\in\ku_{\km}$, an associated vector field $\widetilde{u}$ on 
$\mathcal{O}$ is defined as follows, if 
$f\in\mathcal{O}$, then
\begin{equation}
	\widetilde{u}_{f}=-\mathrm{ad}^{*}(u)f\in T_{f} \mathcal{O}.
	\label{eq:5.4.3dec19s}
\end{equation}
Such vector fields span the 
whole tangent space at each point. Let $\omega_{L}$ denote the real 
$2$-form on $\mathcal{O}$ such that if $u,v\in\ku_{\km}$, $f\in \mathcal{O}$,
\begin{equation}
	\omega_{L}(\widetilde{u},\widetilde{v})_{f}=-\langle f, [u,v]\rangle.
	\label{eq:5.4.4dec19s}
\end{equation}

Then $\omega_{L}$ is a 
$U_{M}$-invariant symplectic form on 
$\mathcal{O}$. Put $r^{+}=\frac{1}{2}\dim 
\ku_{\km}/\kt$. In fact, if we can define an almost complex structure 
on $T\mathcal{O}$ such that the 
holomorphic tangent bundle is given by the positive root system 
$R^{+}(\ku_{\km},\kt)$. Then 
$(\mathcal{O},\omega_{L})$ become a closed 
K\"{a}hler manifold, and $r^{+}$ is its complex dimension.

The Liouville measure on 
$\mathcal{O}$ is defined as follows,
\begin{equation}
	d\mu_{L}=\frac{(\omega_{L})^{r^{+}}}{(r^+)!}.
	\label{eq:5.4.5dec19s}
\end{equation}
It is invariant by the left action of $U_{M}$. Let $\mathrm{Vol_L}(\mathcal{O})$ denote the symplectic volume of $\mathcal{O}$ with respect to 
the Liouville measure. Then we have (see \cite[Proposition 7.26]{berline2003heat})
\begin{equation}
	\mathrm{Vol_{L}}(\mathcal{O})=\Pi_{\alpha^{0}\in 
	R^{+}(\ku_{\km},\kt)}\frac{\langle \alpha^{0}, \lambda+
	\rho_{\ku_{\km}}\rangle}{\langle\alpha^{0},\rho_{\ku_{\km}}\rangle}=\dim V_{\lambda}.
	\label{eq:5.4.40qq20}
\end{equation}
The second identity is the Weyl dimension formula (see 
\cite[Theorem 4.48]{knapp1986representation}).

By the Kirillov formula (see \cite[Theorem 8.4]{berline2003heat}), if $y\in\ku_{\km}$, we have
\begin{equation}
	\widehat{A}^{-1}(\mathrm{ad}(y)|_{\ku_{\km}})\mathrm{Tr}^{V_{\lambda}}[\rho^{V_{\lambda}}(e^{y})]=\int_{f\in\mathcal{O}} e^{2\pi i\langle f, y\rangle} d\mu_{L}.
	\label{eq:5.4.6dec19s}
\end{equation}

To shorten the notation here, if $k\in T$, put $Y=U_{M}(k)^{0}$ with Lie algebra $\y=\ku_{\km}(k)$. 
Then $T\subset Y$, and it also a maximal torus of $Y$. 

In the sequel, we will give a generalized version of 
\eqref{eq:5.4.6dec19s} for describing the function 
$\mathrm{Tr}^{V_{\lambda}}[\rho^{V_{\lambda}}(ke^{y})]$ with $y\in\y$.

Let $\kq$ be the orthogonal space of $\y$ in 
$\ku_{\km}$ with respect to $B$, so that
\begin{equation}
	\ku_{\km}=\y\oplus\kq.
	\label{eq:5.4.7ss20kk}
\end{equation}
Since the adjoint action of $T$ preserves the splitting 
in \eqref{eq:5.4.7ss20kk}. Then $R(\ku_{\km},\kt)$ splits into two 
disjoint parts
\begin{equation}
	R(\ku_{\km},\kt)=R(\y,\kt)\cup R(\kq,\kt),
\end{equation}
where $R(\kq,\kt)$ is just the set of real roots for the adjoint 
action of $\kt$ on $\kq_{\C}$.

The positive root system $R^{+}(\ku_{\km},\kt)$ induces a positive 
root system $R^{+}(\y,\kt)$. Set
\begin{equation}
	R^{+}(\kq,\kt)=R^{+}(\ku_{\km},\kt)\cap R(\kq,\kt).
	\label{eq:5.4.9bath}
\end{equation}
Then we have the disjoint union as follows,
\begin{equation}
	R^{+}(\ku_{\km},\kt)=R^{+}(\y,\kt)\cup R^{+}(\kq,\kt).
\end{equation}

Put
\begin{equation}
	\rho_{\y}=\frac{1}{2}\sum_{\alpha^{0}\in R^{+}(\y,\kt)} 
	\alpha^{0},\; \rho_{\kq}=\frac{1}{2}\sum_{\alpha^{0}\in R^{+}(\kq,\kt)} 
	\alpha^{0}.
\end{equation}
Then
\begin{equation}
	\rho_{\ku_{\km}}=\rho_{\y}+\rho_{\kq}\in\kt^{*}.
\end{equation}

Let $\mathcal{C}\subset \kt^{*}$ denote the Weyl chamber 
corresponding to $R^{+}(\ku_{\km},\kt)$, and let $\mathcal{C}_{0}\subset \kt^{*}$ denote the Weyl chamber 
corresponding to $R^{+}(\y,\kt)$. Then $\mathcal{C}\subset 
\mathcal{C}_{0}$.

Let $W(U_{M},T)$, $W(Y,T)$ be the Weyl groups associated 
with the pairs $(U_{M},T)$, $(Y,T)$ respectively. Then 
$W(Y,T)$ is canonically a subgroup of $W(U_{M},T)$. Put
\begin{equation}
	W^{1}(k)=\{\omega\in W(U_{M},T)\;|\; 
	\omega(\mathcal{C})\subset \mathcal{C}_{0}\}.
	\label{eq:5.4.13ss20}
\end{equation}
Note that the set $W^{1}(k)$ is similar to the set $W_{u}$ defined in 
\eqref{eq:5.3.6ss20}.

\begin{lemma}\label{lm:5.4.1vogel}
	The inclusion $W^{1}(k)\hookrightarrow W(U_{M},T)$ induces 
	a bijection between $W^{1}(k)$ and the quotient 
	$W(Y,T)\backslash W(U_M,T)$.
\end{lemma}
\begin{proof}
	This lemma follows from that $W(Y, T)$ acts simply 
	transitively on the Weyl chambers associated with $(\y,\kt)$.
\end{proof}

Let $\mathcal{O}^{k}$ denote the fixed point set of the holomorphic 
action of $k$ on $\mathcal{O}$. We embeds $\y^{*}$ in $\ku_{\km}^{*}$ 
by the splitting \eqref{eq:5.4.7ss20kk}. Then
\begin{equation}
	\mathcal{O}^{k}=\mathcal{O}\cap \y^{*}.
\end{equation}

\begin{lemma}[{see \cite[I.2, Lemma (7)]{DufloHeckmanVergne1984}, \cite[Lemmas 6.1.1, 7.2.2]{BOUAZIZ19871}}]
	As subsets of $\y^{*}$, we have the following identification,
	\begin{equation}
		\mathcal{O}^{k}=\cup_{\sigma\in 
		W^{1}(k)}\mathrm{Ad}^{*}(Y)(\sigma(\lambda+\rho_{\ku_{\km}}))\subset \y^{*},
	\end{equation}
	where the union is disjoint.
\end{lemma}

For each $\sigma\in W^{1}(k)$, put
\begin{equation}
	\mathcal{O}^{k}_{\sigma(\lambda+\rho_{\ku_{\km}})}=\mathrm{Ad}^{*}(Y)(\sigma(\lambda+\rho_{\ku_{\km}}))\subset \y^{*}.
	\label{eq:5.4.17qq20}
\end{equation}
Let $d\mu^{k}_{\sigma}$ denote the 
Liouville measure on 
$\mathcal{O}^{k}_{\sigma(\lambda+\rho_{\ku_{\km}})}$ as defined in 
\eqref{eq:5.4.5dec19s}.

If $\delta\in \kt^{*}$ is (real) analytically integrable, let 
$\xi_{\delta}$ denote the character of $T$ with differential $2\pi 
i\delta$. Note that for $\sigma\in W^{1}(k)$, 
$\sigma\rho_{\ku_{\km}}+\rho_{\ku_{\km}}$ is analytically integrable 
even $\rho_{\ku_{\km}}$ may not be analytically integrable. 

\begin{definition}\label{def:5.4.3ss20}
	For $\sigma\in W^{1}(k)$, set
	\begin{equation}
		\varphi_{k}(\sigma,\lambda)=\varepsilon(\sigma)\frac{\xi_{\sigma(\lambda+\rho_{\ku_{\km}})+\rho_{\ku_{\km}}}(k)}{\Pi_{\alpha^{0}\in R^{+}(\kq,\kt)}(\xi_{\alpha^{0}}(k)-1)}.
		\label{eq:5.4.18kk20}
	\end{equation}
\end{definition}

Note that if $y\in\y$, the following analytic function 
\begin{equation}
	\frac{\det (1-e^{\mathrm{ad}(y)}\mathrm{Ad}(k))|_{\kq}}{\det (1-\mathrm{Ad}(k))|_{\kq}}
\end{equation}
has a square root which is analytic in $y\in\y$ and equals to $1$ at 
$y=0$. We denote this square 
root by
\begin{equation}
	\big[\frac{\det (1-e^{\mathrm{ad}(y)}\mathrm{Ad}(k))|_{\kq}}{\det 
	(1-\mathrm{Ad}(k))|_{\kq}}\big]^{\frac{1}{2}}.
\end{equation}

The following theorem is a special case of a generalized Kirillov formula obtained by Duflo, 
Heckman and Vergne \cite[II. 3, Theorem (7)]{DufloHeckmanVergne1984}. We 
will also include a simpler proof for the sake of completeness.
\begin{theorem}[Generalized Kirillov formula]\label{thm:5.4.4ss20}
	For $y\in\y$, we have the following identity of analytic 
	functions,
	\begin{equation}
		\begin{split}
			&\widehat{A}^{-1}(\mathrm{ad}(y)|_{\y})\big[\frac{\det (1-
			e^{\mathrm{ad}(y)}\mathrm{Ad}(k))|_{\kq}}{\det(1-\mathrm{Ad}(k))
			_{\kq}}\big]^{\frac{1}{2}} 
			\mathrm{Tr}^{V_{\lambda}}[\rho^{V_{\lambda}}(ke^{y})]\\
			&=\sum_{\sigma\in W^{1}(k)} 
			\varphi_{k}(\sigma,\lambda)\int_{f\in\mathcal{O}^{k}
			_{\sigma(\lambda+\rho_{\ku_{\km}})}}e^{2\pi i \langle f,y\rangle}
			d\mu^{k}_{\sigma}.
		\end{split}
		\label{eq:haha20imp}
	\end{equation}
	If $k=1$, \eqref{eq:haha20imp} is reduced to 
	\eqref{eq:5.4.6dec19s}.
\end{theorem}
\begin{proof}
	Let $\kt'$ denote the set of regular element in $\kt$ associated 
	with the root $R(\ku_{\km},\kt)$, which is an open dense subset 
	of $\kt$. Since both sides of \eqref{eq:haha20imp} are analytic 
	and invariant 
	by the adjoint action of $Y$, then we only need to prove 
	\eqref{eq:haha20imp} for $y\in\kt'$.

	We firstly compute the left-hand side of \eqref{eq:haha20imp}. 	
	
	For $y\in\kt'$, we have
	\begin{equation}
		\widehat{A}^{-1}(\mathrm{ad}(y)|_{\y})=\Pi_{\alpha^{0}\in 
		R^{+}(\y,\kt)}\frac{e^{\pi i \langle \alpha^{0},y\rangle}-e^{-\pi i 
		\langle \alpha^{0},y\rangle}}{\langle 2\pi 
		i\alpha^{0},y\rangle}.
		\label{eq:5.4.23qq20}
	\end{equation}
	
	Let $y_{0}\in\kt$ be such that $k=\exp(y_{0})$. Then
	\begin{equation}
		\begin{split}
			&\big[\frac{\det (1-e^{\mathrm{ad}(y)}\mathrm{Ad}(k))|_{\kq}}{\det 
			(1-\mathrm{Ad}(k))|_{\kq}}\big]^{\frac{1}{2}}\\
			&=\Pi_{\alpha^{0}\in 
			R^{+}(\kq,\kt)}\frac{e^{\pi i \langle 
			\alpha^{0},y+y_{0}\rangle}-e^{-\pi i 
			\langle \alpha^{0},y+y_{0}\rangle}}{e^{\pi i \langle
			\alpha^{0},y_{0}\rangle}-e^{-\pi i 
			\langle \alpha^{0},y_{0}\rangle}}.
		\end{split}
		\label{eq:5.4.24qq20}
	\end{equation}
	
	By the Weyl character formula for $(U_{M},T)$, we get
	\begin{equation}
		\begin{split}
			&\mathrm{Tr}^{V_{\lambda}}[\rho^{V_{\lambda}}(ke^{y})]\\
			&=\mathrm{Tr}^{V_{\lambda}}[\rho^{V_{\lambda}}(e^{y+y_{0}})]\\
			&=\frac{\sum_{\omega\in 
			W(\ku_{\km,\C},\kt_{\C})}\varepsilon(\omega)e^{2\pi i \langle 
			\omega(\lambda+\rho_{\ku_{\km}}),y+y_{0}\rangle}}{\Pi_{\alpha^{0}\in 
			R^{+}(\ku_{\km},\kt)} \big(e^{\pi i \langle 
			\alpha^{0},y+y_{0}\rangle}-e^{-\pi i 
			\langle \alpha^{0},y+y_{0}\rangle}\big)}.
		\end{split}
		\label{eq:5.4.25qq20}
	\end{equation}

	Note that we have 
	$\xi_{\alpha_{0}}(k)=1$ for $\alpha_{0}\in R^{+}(\y,\kt)$, then
	\begin{equation}
		\begin{split}
			\Pi_{\alpha^{0}\in 
			R^{+}(\y,\kt)}\frac{e^{\pi i \langle 
			\alpha^{0},y+y_{0}\rangle}-e^{-\pi i 
			\langle \alpha^{0},y+y_{0}\rangle}}{e^{\pi i \langle
			\alpha^{0},y\rangle}-e^{-\pi i 
			\langle \alpha^{0},y\rangle}}=e^{-2\pi i\langle 
			\rho_{\y},y_{0}\rangle}.
		\end{split}
		\label{eq:5.4.26qq20}
	\end{equation}

	Combining \eqref{eq:5.4.23qq20} - \eqref{eq:5.4.26qq20}, we get 
	the left-hand side of \eqref{eq:haha20imp} is equal the following 
	function,
	\begin{equation}
		\frac{e^{2\pi i\langle \rho_{\y},y_{0}\rangle}}{ \Pi_{\alpha^{0}\in 
		R^{+}(\y,\kt)} \langle 2\pi 
		i\alpha^{0},y\rangle}\frac{\sum_{\omega\in 
		W(\ku_{\km,\C},\kt_{\C})}\varepsilon(\omega)e^{2\pi i \langle 
		\omega(\lambda+\rho_{\ku_{\km}}),y+y_{0}\rangle}}{\Pi_{\alpha^{0}\in 
		R^{+}(\kq,\kt)} \big(e^{\pi i \langle \alpha^{0},y\rangle}-e^{-\pi i 
		\langle \alpha^{0},y\rangle}\big)}.
		\label{eq:5.4.27qq20}
	\end{equation}
	
	Now we show that the right-hand side of \eqref{eq:haha20imp} is 
	also equal to \eqref{eq:5.4.27qq20}.

	Note that for $\omega\in W(Y, T)$, 
	$\omega\rho_{\ku_{\km}}-\rho_{\ku_{\km}}$ is analytically 
	integrable. We claim that if $\omega\in W(Y, T)$, then
	\begin{equation}
		\xi_{\omega\rho_{\ku_{\km}}-\rho_{\ku_{\km}}}(k)=e^{2\pi i\langle 
		\omega\rho_{\ku_{\km}}-\rho_{\ku_{\km}},y_{0}\rangle }=1.
		\label{eq:5.4.30ss20k}
	\end{equation}
	Actually, we have 
	$\xi_{2\rho_{\ku_{\km}}}(k)=\xi_{2\omega\rho_{\ku_{\km}}}(k)=1$.
	Then, after taking the square roots, we get that 
	$\xi_{\omega\rho_{\ku_{\km}}-\rho_{\ku_{\km}}}(k)=\xi_{\omega\rho_{\ku_{\km}}-\rho_{\ku_{\km}}}(e^{y_{0}})=\pm 1$. The continuity of the character implies exactly \eqref{eq:5.4.30ss20k}.

	As a consequence of \eqref{eq:5.4.30ss20k}, we get that for 
	$\sigma\in W^{1}(k)$, if $\omega\in W(Y, T)$, then
	\begin{equation}
		e^{2\pi i \langle 
		\omega\sigma(\lambda+\rho_{\ku_{\km}}),y_{0}\rangle}=e^{2\pi i \langle 
		\sigma(\lambda+\rho_{\ku_{\km}}),y_{0}\rangle}.
		\label{eq:5.4.37qq20}
	\end{equation}
	
	For $\sigma\in W^{1}(k)$, since 
	$\sigma(\lambda+\rho_{\ku_{\km}})\in \mathcal{C}_{0}$ and $y$ is 
	regular, by \cite[Corollary 7.25]{berline2003heat}, we have
	\begin{equation}
		\begin{split}
			&\int_{f\in\mathcal{O}^{k}_{\sigma(\lambda+\rho_{\ku_{\km}})}}e^{2\pi i \langle f,y\rangle}d\mu^{k}_{\sigma} \\
			&=\frac{1}{\Pi_{\alpha^{0}\in R^{+}(\y,\kt)} \langle 2\pi
			i\alpha^{0},y\rangle}\sum_{\omega\in W(Y, T)}
			\varepsilon(\omega)e^{2\pi i \langle \omega\sigma(\lambda+
			\rho_{\ku_{\km}}),y\rangle}.
		\end{split}
		\label{eq:5.4.38qq20}
	\end{equation}
	
	We rewrite $\varphi_{k}(\sigma,\lambda)$ as follows,
	\begin{equation}
		\varepsilon(\sigma)\frac{e^{2\pi i\langle \rho_{\y},y_{0}\rangle}}{\Pi_{\alpha^{0}\in 
		R^{+}(\kq,\kt)} \big(e^{\pi i \langle \alpha^{0},y\rangle}-e^{-\pi i 
		\langle \alpha^{0},y\rangle}\big)}e^{2\pi i \langle 
		\sigma(\lambda+\rho_{\ku_{\km}}),y_{0}\rangle}.
		\label{eq:5.4.39qq20}
	\end{equation}
	
	Combining together Lemma \ref{lm:5.4.1vogel} and \eqref{eq:5.4.37qq20} - 
	\eqref{eq:5.4.39qq20}, a direct computation shows that the 
	right-hand side of  \eqref{eq:haha20imp} is given exactly by 
	\eqref{eq:5.4.27qq20}. This completes the proof of our theorem. 
\end{proof}

\begin{remark}\label{rm:5.4.5s}
	Let $C^{0}$ denote the identity component of the center of 
	$Y$, and let $Y_{\mathrm{ss}}$ be the closed analytic subgroup of 
	$Y$ associated with $\y_{\mathrm{ss}}=[\y,\y]$. By Weyl's theorem \cite[Theorem 
	4.26]{knapp1986representation}, the universal covering group of 
	$Y_{\mathrm{ss}}$ is compact, which we denote by 
	$\widetilde{Y}_{\mathrm{ss}}$. Put
	\begin{equation}
		\widetilde{Y}=C^{0}\times \widetilde{Y}_{\mathrm{ss}}.
	\end{equation}
	Then $\widetilde{Y}$ is clearly a finite central extension of 
	$Y$. Let $\widetilde{T}$ be the maximal torus of 
	$\widetilde{Y}$ associated with the Cartan subalgebra $\kt$, 
	which is also a finite extension of $T$. By \cite[Corollary 
	4.25]{knapp1986representation}, the weight $\rho_{\ku_{\km}}, 
	\rho_{\y}$ are analytically integrable with respect to 
	$\widetilde{T}$, since they are algebraically integrable 
	(\cite[Propositions 4.15 \& 4.33]{knapp1986representation}).

	Note that for $\sigma\in W^{1}(k)$, 
	$\sigma(\lambda+\rho_{\ku_{\km}})$ is regular and positive with 
	respect to $R^{+}(\y,\kt)$, thus 
	$\sigma(\lambda+\rho_{\ku_{\km}})-\rho_{\y}$ is nonnegative with 
	respect to $R^{+}(\y,\kt)$ by the property of $\rho_{\y}$ 
	(\cite[Proposition 4.33]{knapp1986representation}). Since now
	$\sigma(\lambda+\rho_{\ku_{\km}})-\rho_{\y}$ is also analytically integrable 
	with respect to $\widetilde{T}$,
	then it is a dominant weight 
	for $(\widetilde{Y},\widetilde{T})$ with respect to 
	$R^{+}(\y,\kt)$. In this case, let 
	$V^{k}_{\lambda,\sigma}$ be the 
	irreducible unitary representation of $\widetilde{Y}$ with 
	highest weight $\sigma(\lambda+\rho_{\ku_{\km}})-\rho_{\y}$. Then 
	by \eqref{eq:5.4.6dec19s}, \eqref{eq:haha20imp}, we get that for 
	$y\in\y$,
	\begin{equation}
		\begin{split}
			&\big[\frac{\det (1-
			e^{\mathrm{ad}(y)}\mathrm{Ad}(k))|_{\kq}}{\det(1-\mathrm{Ad}(k))
			_{\kq}}\big]^{\frac{1}{2}} 
			\mathrm{Tr}^{V_{\lambda}}[\rho^{V_{\lambda}}(ke^{y})]\\
			&=\sum_{\sigma\in W^{1}(k)} 
			\varphi_{k}(\sigma,\lambda)\mathrm{Tr}^{V^{k}_{\lambda,\sigma}}[\rho^{V^{k}_{\lambda,\sigma}}(e^{y})].
		\end{split}
		\label{eq:haha20imp12}
	\end{equation}
\end{remark}

%%%%%%%%%%%%%%%%%%%%%%%%%%%%%%%%%%%%%%%%%%%%%%%%%%%%%%%%%%%%%%%%%%
\section{A geometric localization formula for orbital 
integrals}\label{section6paris}

Recall that $G_\C$ is the complexification of $G$ with Lie algebra 
$\g_\C$, 
and that  $G$, $U$ are the analytic subgroups of $G_\C$ with Lie 
algebra $\g$, $\ku$ respectively. In this section, we always assume 
that $U$ is compact, we do not require that $G$ has compact center. 
We need not to assume $\delta(G)=1$ either.

Under the settings in Subsection \ref{subsection4.1paris}, for 
$t>0$ and semisimple $\gamma\in G$, we set
\begin{equation}
	\mathcal{E}_{X,\gamma}(F,t)=\mathrm{Tr_s}^{[\gamma]}\big[(N^{\Lambda^{\bullet}(T^*X)}-\frac{m}{2})\exp(-t\mathbf{D}^{X,F,2}/2)\big].
	\label{eq:6.0.1kk}
\end{equation}
The indice $X$, $F$ in this notation indicate precisely  the symmetric space 
and the flat vector bundle which are concerned for defining the orbital 
integrals.  

If $\gamma\in G$ is semisimple, then there exists a unique elliptic 
element $\gamma_{e}$ and a unique hyperbolic element $\gamma_{h}$ in 
$G$, such that $\gamma=\gamma_{e}\gamma_{h}=\gamma_{h}\gamma_{e}$. 
Here, we will show that $\mathcal{E}_{X,\gamma}(F,t)$ becomes a sum of the orbital integrals associated with 
$\gamma_{h}$, but defined for the centralizer of $\gamma_{e}$ instead 
of $G$. This suggests that the elliptic part of $\gamma$ should lead to 
a localization for the geometric orbital integrals.

We still fix a maximal torus $T$ of $K$ with Lie algebra $\kt$. For simplicity, if 
$\gamma\in G$ is semisimple, we may and we will assume that
\begin{equation}
	\gamma=e^ak, k\in T, a\in \pp, \mathrm{Ad}(k^{-1})a=a.
	\label{eq:6.2.1ss20}
\end{equation}
In this case,
\begin{equation}
	\gamma_{e}=k\in T,\;\gamma_{h}=e^{a}.
\end{equation}

Recall that $Z(\gamma_{e})^{0}$ is the identity component of the 
centralizer of $\gamma_{e}$ in $G$. Then
\begin{equation}
	\gamma_{h}\in Z(\gamma_{e})^{0}.
\end{equation}
The Cartan involution $\theta$ preserves $Z(\gamma_{e})^{0}$ such that 
$Z(\gamma_{e})^{0}$ is a connected linear reductive Lie group. Then 
we have the following diffeomorphism
\begin{equation}
	Z(\gamma_{e})^{0}=K(\gamma_{e})^{0}\exp(\pp(\gamma_{e})).
\end{equation}
It is clear that  $\delta(Z(\gamma_{e})^{0})=\delta(G)$.

Recall that $T_{U}$ is a maximal torus of $U$ with Lie algebra 
$\kt_{U}=\ii\kb\oplus \kt\subset\ku$. Let $R^{+}(\ku,\kt_{U})$ be a 
positive root system for $R(\ku,\kt_{U})$, which is not necessarily 
the same as in Subsection \ref{section5.3pl} when $\delta(G)=1$.

Since $U$ is the compact form of $G$, then $U(\gamma_{e})^{0}$ is the compact 
form for $Z(\gamma_{e})^{0}$. Moreover, $T_{U}$ is also a maximal torus of 
$U(\gamma_{e})^{0}$. Let $R(\ku(\gamma_{e}),\kt_{U})$ be the corresponding 
real root system with the positive root system 
$R^{+}(\ku(\gamma_{e}),\kt_{U})=R(\ku(\gamma_{e}),\kt_{U})\cap 
R^{+}(\ku,\kt_{U})$. Let $\rho_{\ku}$, $\rho_{\ku(\gamma_{e})}$ be the 
corresponding half sums of positive roots.

Let $\widetilde{U}(\gamma_{e})$ be a connected finite covering group of 
$U(\gamma_{e})^{0}$ such that $\rho_{\ku}$, $\rho_{\ku(\gamma_{e})}$ are 
analytically integrable with respect to the maximal torus $\widetilde{T}_{U}$ of $\widetilde{U}(\gamma_{e})$ associated with 
$\kt_{U}$. It always exists by a similar construction as in 
Remark \ref{rm:5.4.5s}.

Let $\widetilde{K}(\gamma_{e})$ be the analytic subgroup of $\widetilde{U}(\gamma_{e})$ associated with 
Lie algebra $\kk(\gamma_{e})$. By \cite[Proposition 
7.12]{knapp2002liegroupe}, $\widetilde{U}(\gamma_{e})$ has a unique 
complexification $\widetilde{U}(\gamma_{e})_{\C}$ which is a connected 
linear reductive Lie group. Let $\widetilde{Z}(\gamma_{e})$ be the 
analytic subgroup of $\widetilde{U}(\gamma_{e})_{\C}$ associated with 
$\z(\gamma_{e})\subset \ku(\gamma_{e})_{\C}=\z(\gamma_{e})_{\C}$. Then we have 
the following Cartan decomposition
\begin{equation}
	\widetilde{Z}(\gamma_{e})= 
	\widetilde{K}(\gamma_{e})\exp(\pp(\gamma_{e})).
\end{equation}
We still denote by $\theta$ the corresponding Cartan involution on 
$\widetilde{Z}(\gamma_{e})$.

The Lie group $\widetilde{Z}(\gamma_{e})$ is a finite covering group of 
$Z(\gamma_{e})^{0}$. Moreover, we have the identification of symmetric spaces
\begin{equation}
	X(\gamma_{e})\simeq \widetilde{Z}(\gamma_{e})/\widetilde{K}(\gamma_{e}).
\end{equation}
Note that even under an additional assumption that $G$ has compact 
center, $\widetilde{Z}(\gamma_{e})$ may still have noncompact center.

Let $\lambda$ be a dominant weight for $(U, T_{U})$ with respect to 
$R^{+}(\ku,\kt_{U})$. Let 
$(E_{\lambda},\rho^{E_{\lambda}})$ be the associated irreducible unitary representation 
of $U$. As before, let $(F_{\lambda}, h^{F_{\lambda}})$ be the 
corresponding homogeneous vector 
bundle on $X$ with the $G$-invariant flat connection $\nabla^{F_{\lambda},f}$. Let 
$\mathbf{D}^{X,F_{\lambda},2}$ denote the 
associated de Rham-Hodge Laplacian.

Let $W^{1}_{U}(\gamma_{e})\subset W(U,T_{U})$ be the set defined as 
in \eqref{eq:5.4.13ss20} but with respect to the group $U$ and to 
$\gamma_{e}=k\in T\subset T_{U}$. As in Definition \ref{def:5.4.3ss20}, for $\sigma\in W^{1}_{U}(\gamma_{e})$, set
\begin{equation}
	\varphi^{U}_{\gamma_{e}}(\sigma,\lambda)=\varepsilon(\sigma)\frac{\xi_{\sigma(\lambda+\rho_{\ku})+\rho_{\ku}}(\gamma_{e})}{\Pi_{\alpha^{0}\in R^{+}(\ku^{\perp}(\gamma_{e}),\kt_{U})}(\xi_{\alpha^{0}}(\gamma_{e})-1)}.
	\label{eq:6.2.7kk20}
\end{equation}

As explained in Remark \ref{rm:5.4.5s}, if $\sigma\in 
W^{1}_{U}(\gamma_{e})$, then 
$\sigma(\lambda+\rho_{\ku})-\rho_{\ku(\gamma_{e})}$ is a dominant 
weight of $\widetilde{U}(\gamma_{e})$ with respect to 
$R^{+}(\ku(\gamma_{e}),\kt_{U})$. Let $E_{\sigma,\lambda}$ be the 
irreducible unitary representation of $\widetilde{U}(\gamma_{e})$ with 
highest weight $\sigma(\lambda+\rho_{\ku})-\rho_{\ku(\gamma_{e})}$. 

We extend $E_{\sigma,\lambda}$ to an irreducible representation 
of $\widetilde{Z}(\gamma_{e})$ by the unitary trick. Then 
$F_{\sigma,\lambda}=\widetilde{Z}(\gamma_{e})\times_{\widetilde{K}(\gamma_{e})} E_{\sigma,\lambda}$ is a homogeneous vector bundle on $X(\gamma_{e})$ with an invariant flat connection $\nabla^{F_{\sigma,\lambda},f}$ as explained in Subsection \ref{section3.6}. Let 
$\mathbf{D}^{X(\gamma_{e}),F_{\sigma,\lambda},2}$ denote the 
associated de Rham-Hodge Laplacian acting on 
$\Omega^{\cdot}(X(\gamma_{e}), F_{\sigma,\lambda})$. 	

We also view $\gamma_{h}=e^{a}$ 
as a hyperbolic element in $\widetilde{Z}(\gamma_{e})$. For $\sigma\in W^{1}_{U}(\gamma_{e})$, as in \eqref{eq:6.0.1kk}, we 
set
\begin{equation}
	\mathcal{E}_{X(\gamma_{e}),\gamma_{h}}(F_{\sigma,\lambda},t)=\mathrm{Tr_s}^{[\gamma_{h}]}\big[(N^{\Lambda^{\bullet}(T^*X(\gamma_{e}))}-\frac{p'}{2})\exp(-t\mathbf{D}^{X(\gamma_{e}),F_{\sigma,\lambda},2}/2)\big].
\end{equation}
Note that we use $B|_{\z(\gamma_{e})}$ on $\z(\gamma_{e})$ to define this orbital integral for $\widetilde{Z}(\gamma_{e})$.

Set
\begin{equation}
	c(\gamma)=\bigg|\frac{\det 
	(1-\mathrm{Ad}(\gamma_{e}))|_{\z^{\perp}(\gamma_{e})}}{\det 
	(1-\mathrm{Ad}(\gamma))|_{\z^{\perp}(\gamma_{e})}} \bigg|^{1/2}>0.
	\label{eq:6.2.9kk20}
\end{equation}
In particular, $c(\gamma_{e})=1$.

The following theorem is essentially a consequence of the generalized 
Kirillov formula in Theorem \ref{thm:5.4.4ss20}. 
\begin{theorem}\label{thm:6.2.1ss}
	Let $\gamma\in G$ be given as in \eqref{eq:6.2.1ss20}. For $t>0$, we 
	have the following identity,
	\begin{equation}
		\mathcal{E}_{X,\gamma}(F_{\lambda},t)=c(\gamma)\sum_{\sigma\in 
		W^{1}_{U}(\gamma_{e})}\varphi^{U}_{\gamma_{e}}(\sigma,\lambda) \mathcal{E}_{X(\gamma_{e}),\gamma_{h}}(F_{\sigma,\lambda},t).
		\label{eq:aoyess20}
	\end{equation}
	We call \eqref{eq:aoyess20} a localization formula for the
	geometric orbital integral.
\end{theorem}
\begin{proof}
	Set $p'=\dim \pp(\gamma_{e})=\dim X(\gamma_{e})$. At first, if $m$ is even, then $p'$ is even. Then the both sides of 
	\eqref{eq:aoyess20} are $0$ by \cite[Theorem 
	7.9.1]{bismut2011hypoelliptic}.
	
	If $m$ is odd, then $p'$ is odd, and 
	$\delta(G)=\delta(Z(\gamma_{e})^{0})$ is odd. If $\delta(G)\geq 3$, 
	then the both sides of \eqref{eq:aoyess20} are $0$ by \cite[Theorem 
	7.9.1]{bismut2011hypoelliptic}.
	
	Now we consider the case where 
	$\delta(G)=\delta(Z(\gamma_{e})^{0})=1$. If $\gamma$ can not be 
	conjugated into $H$ by an element in $G$, then $\gamma_{h}$ can 
	not be conjugated into $H$ by an element in $Z(\gamma_{e})^{0}$. Then 
	the both
	sides of \eqref{eq:aoyess20} are $0$ by Proposition 
	\ref{prop:6.1.1ss}.
	
	Now we assume that $\delta(G)=1$ and $a\in\kb$. Note that 
	$\z(\gamma)$ is the centralizer of $\gamma_{h}$ in 
	$\z(\gamma_{e})$. We will prove \eqref{eq:aoyess20} using \eqref{eq:6.2.8pl} 
	
	For $y\in \kk(\gamma)$, let $J^{\sim}_{\gamma_{h}}(y)$ be the function 
	defined in \eqref{def:3.3.1ss20} for $\gamma_{h}=e^{a}\in 
	\widetilde{Z}(\gamma_{e})$, 
	\begin{equation}
		\begin{split}
			J^{\sim}_{\gamma_{h}}(y)=\frac{1}{|\det 
			(1-\mathrm{Ad}(\gamma_{h}))|_{\z^{\perp}_0 \cap \z(\gamma_{e})}|^{1/2}} 
			\frac{\widehat{A}(i\mathrm{ad}(y)|_{\pp(\gamma)})}{\widehat{A}(i\mathrm{ad}(y)|_{\kk(\gamma)})}.
		\end{split}
		\label{eq:6.2.11jan20}
	\end{equation}

	The Casimir operator $C^{\ku(\gamma_{e}),E_{\sigma,\lambda}}$ acts on $E_{\sigma,\lambda}$ by the scalar 
	given 
	\begin{equation}
		-4\pi^{2} 
		(|\lambda+\rho_{\ku}|^{2}-|\rho_{\ku(\gamma_{e})}|^{2}).
		\label{eq:6.2.14ss20}
	\end{equation}
	Similar to \eqref{eq:4.1.13vogel}, set
	\begin{equation}
		\beta_{\z(\gamma_{e})}=\frac{1}{16}\mathrm{Tr}^{\pp(\gamma_{e})}[C^{\kk(\gamma_{e}),\pp(\gamma_{e})}]+\frac{1}{48}\mathrm{Tr}^{\kk(\gamma_{e})}[C^{\kk(\gamma_{e}),\kk(\gamma_{e})}].
	\end{equation}
	Then by \cite[Propositions 2.6.1 \& 7.5.1]{bismut2011hypoelliptic},
	\begin{equation}
		2\pi^{2}|\rho_{\ku(\gamma_{e})}|^{2}=-\beta_{\z(\gamma_{e})}.
		\label{eq:6.2.15s20}
	\end{equation}
	
	By \eqref{eq:6.2.8pl}, \eqref{eq:6.2.14ss20}, 
	\eqref{eq:6.2.15s20}, for $\sigma\in W^{1}_{U}(\gamma_{e})$, we get
	\begin{equation}
		\begin{split}
			&\mathcal{E}_{X(\gamma_{e}),\gamma_{h}}(F_{\sigma,\lambda},t)=\frac{e^{-\frac{|a|^2}{2t}}}{(2\pi t)^{p/2}} 
			\exp\big(-2\pi^{2} t|\lambda+\rho_{\ku}|^{2}\big)\\
			&\qquad\qquad\cdot \int_{\kk(\gamma)} 
			J^{\sim}_{\gamma_{h}}(y)\mathrm{Tr_s}^{\Lambda^{\bullet}(\pp(\gamma_{e})^*)}\big[(N^{\Lambda^{\bullet}(\pp(\gamma_{e})^*)}-\frac{p'}{2})\exp(-i\mathrm{ad}(y))\big]\\
			&\qquad\qquad\qquad
			\cdot\mathrm{Tr}^{E_{\sigma,\lambda}}[\exp(-i\rho^{E_{\sigma,\lambda}}(y))] e^{-|y|^2/2t}\frac{dy}{(2\pi t)^{q/2}}.
		\end{split}
		\label{eq:6.2.16ss20}
	\end{equation}
	
	Note that $\dim \pp^{\perp}(\gamma_{e})$ is even. We claim that if $y\in\kk(\gamma)$, then
	\begin{equation}
		\begin{split}
			&\mathrm{Tr_s}^{\Lambda^{\bullet}(\pp^*)}\big[(N^{\Lambda^{\bullet}(\pp^*)}-\frac{m}{2})\exp(-i\mathrm{ad}(y))\mathrm{Ad}(k^{-1})\big]\\
			&=\mathrm{Tr_s}^{\Lambda^{\bullet}(\pp(\gamma_{e})^*)}\big[(N^{\Lambda^{\bullet}(\pp(\gamma_{e})^*)}-\frac{p'}{2})e^{-i\mathrm{ad}(y)}\big]\det(1-e^{-i\mathrm{ad}(y)}\mathrm{Ad}(k^{-1}))|_{\pp^{\perp}(\gamma_{e})}.
		\end{split}
		\label{eq:6.2.18ss20}
	\end{equation}
	Indeed, we can verify \eqref{eq:6.2.18ss20} for $y\in\kt$. Since both sides of \eqref{eq:6.2.18ss20} are invariant by adjoint 
	action of $K(\gamma_{e})^{0}$, then \eqref{eq:6.2.18ss20} holds in full 
	generality.
	
	Also $K(\gamma)^{0}$ preserves the splitting
	\begin{equation}
		\pp^{\perp}(\gamma_{e})=\pp^{\perp}_{0}(\gamma)\oplus 
		(\pp^{\perp}(\gamma_{e})\cap \pp^{\perp}_{0}).
		\label{eq:6.2.19kk20}
	\end{equation}
	The action $\mathrm{ad}(a)$ gives an isomorphism between 
	$\pp^{\perp}(\gamma_{e})\cap \pp^{\perp}_{0}$ and 
	$\kk^{\perp}(\gamma_{e})\cap \kk^{\perp}_{0}$ as $K(\gamma)$-vector 
	spaces. 
	
	Note that
	\begin{equation}
		\z^{\perp}(\gamma_{e})\cap \z^{\perp}_{0}=(\pp^{\perp}(\gamma_{e})\cap 
		\pp^{\perp}_{0})\oplus (\kk^{\perp}(\gamma_{e})\cap \kk^{\perp}_{0}).
		\label{eq:6.2.20kk20}
	\end{equation}
	Then 
	\begin{equation}
		\begin{split}		
			&\det(1-e^{-i\mathrm{ad}(y)}\mathrm{Ad}(\gamma_{e}))|_{\pp^{\perp}(\gamma_{e})}\\
			&=\det(1-e^{-i\mathrm{ad}(y)}\mathrm{Ad}(\gamma_{e}))|_{\pp_{0}^{\perp}(\gamma_{e})}[\det(1-e^{-i\mathrm{ad}(y)}\mathrm{Ad}(\gamma_{e}))|_{\z^{\perp}(\gamma_{e})\cap \z^{\perp}_{0}}]^{\frac{1}{2}}.
		\end{split}
		\label{eq:6.2.21kk20}
	\end{equation}
	Here the square root is taken to be positive at $y=0$.
	
	By Definition \ref{def:3.3.1ss20} and \eqref{eq:6.2.11jan20}, for 
	$y\in\kk(\gamma)$,
	\begin{equation}
		\begin{split}
			J_{\gamma}(y)=&J^{\sim}_{\gamma_{h}}(y)\frac{1}{|\det 
			(1-\mathrm{Ad}(\gamma))|_{\z^{\perp}_0\cap 
			\z^{\perp}(\gamma_{e})}|^{1/2}}\\
			&\cdot\bigg[ \frac{1}{\det 
			(1-\mathrm{Ad}(\gamma_{e}))|_{\z^{\perp}_{0}(\gamma)}} \frac{\det 
			(1-\exp(-i\mathrm{ad}(y))\mathrm{Ad}(\gamma_{e}))|_{\kk^{\perp}_{0}(\gamma)}}{\det (1-\exp(-i\mathrm{ad}(y))\mathrm{Ad}(\gamma_{e}))|_{\pp^{\perp}_{0}(\gamma)}}     \bigg]^{1/2}.
		\end{split}
		\label{eq:6.2.22kk20}
	\end{equation}
	
	Combining \eqref{eq:6.2.18ss20}, \eqref{eq:6.2.21kk20} and 
	\eqref{eq:6.2.22kk20}, we get
	\begin{equation}
		\begin{split}
			&J_{\gamma}(y)\mathrm{Tr_s}^{\Lambda^{\bullet}(\pp^*)}\big[(N^{\Lambda^{\bullet}(\pp^*)}-\frac{m}{2})\exp(-i\mathrm{ad}(y))\mathrm{Ad}(\gamma_{e})\big]\\
			=&c(\gamma)J^{\sim}_{\gamma_{h}}(y)\mathrm{Tr_s}^{\Lambda^{\bullet}(\pp(\gamma_{e})^*)}\big[(N^{\Lambda^{\bullet}(\pp(\gamma_{e})^*)}-\frac{p'}{2})e^{-i\mathrm{ad}(y)}\big]\\
			&\cdot\bigg[\frac{\det 
			(1-\exp(-i\mathrm{ad}(y))\mathrm{Ad}(\gamma_{e}))|_{\z^{\perp}(\gamma_{e})}}{\det 
			(1-\mathrm{Ad}(\gamma_{e}))|_{\z^{\perp}(\gamma_{e})}}     \bigg]^{1/2}.
		\end{split}
		\label{eq:6.2.23kk20}
	\end{equation}
	
	Note that for $y\in\kk(\gamma)$, 
	\begin{equation}
		\begin{split}
			&\bigg[\frac{\det(1-\exp(-i\mathrm{ad}(y))\mathrm{Ad}(\gamma_{e}))|_{\z^{\perp}(\gamma_{e})}}{\det(1-\mathrm{Ad}(\gamma_{e}))|_{\z^{\perp}(\gamma_{e})}}     \bigg]^{1/2}\\
			&=\bigg[\frac{\det(1-\exp(-i\mathrm{ad}(y))\mathrm{Ad}(\gamma_{e}))|_{\ku^{\perp}(\gamma_{e})}}{\det(1-\mathrm{Ad}(\gamma_{e}))|_{\ku^{\perp}(\gamma_{e})}}     \bigg]^{1/2}.
		\end{split}
		\label{eq:6.2.24kk20}
	\end{equation}

	By \eqref{eq:6.2.8pl}, \eqref{eq:6.2.14ss20}, \eqref{eq:6.2.15s20}, 
	\eqref{eq:6.2.23kk20} and \eqref{eq:6.2.24kk20},
	we get
	\begin{equation}
		\begin{split}
			&\mathcal{E}_{X,\lambda}(F_{\lambda},t)=c(\gamma)\frac{e^{-\frac{|a|^2}{2t}}}{(2\pi t)^{p/2}} 
			\exp\big(-2\pi^{2} t|\lambda+\rho_{\ku}|^{2}\big)\\
			&\cdot \int_{\kk(\gamma)} 
			J^{\sim}_{\gamma_{h}}(y)\mathrm{Tr_s}^{\Lambda^{\bullet}(\pp(\gamma_{e})^*)}\big[(N^{\Lambda^{\bullet}(\pp(\gamma_{e})^*)}-\frac{p'}{2})e^{-i\mathrm{ad}(y)}\big]\\
			&\bigg[\frac{\det(1-e^{-i\mathrm{ad}(y)}\mathrm{Ad}(\gamma_{e}))|_{\ku^{\perp}(\gamma_{e})}}{\det(1-\mathrm{Ad}(\gamma_{e}))|_{\ku^{\perp}(\gamma_{e})}}     \bigg]^{1/2}\mathrm{Tr}^{E_{\lambda}}[\rho^{E_{\lambda}}(\gamma_{e})e^{-i\rho^{E_{\lambda}}(y)}] e^{-|y|^2/2t}\frac{dy}{(2\pi t)^{q/2}}.
		\end{split}
		\label{eq:6.2.25ss20}
	\end{equation}
	Then \eqref{eq:aoyess20} follows from the \eqref{eq:haha20imp12}, 
	\eqref{eq:6.2.16ss20} and \eqref{eq:6.2.25ss20}. This completes the 
	proof of our theorem.
\end{proof}
\begin{remark}
	A similar consideration can be made for 
	$\mathrm{Tr_s}^{[\gamma]}\big[\exp(-t\mathbf{D}^{X,F_{\lambda},2})\big]$, where \eqref{eq:aoyess20} will become an analogue of the index theorem for orbifolds as in \eqref{eq:1.2.12alpha}. The related computation can be found in \cite[Subsection 10.4]{bismut2019geometric}.
\end{remark}

\section{Full asymptotics of elliptic orbital 
integrals}\label{section7paris}
In this section, we always assume that $\delta(G)=1$ and that $U$ is 
compact. We also use the notation and settings as in Subsections 
\ref{section5.1ss20}, \ref{subsection5.2ss} and \ref{section5.3pl}.

In this section, given a irreducible unitary representation $E$ of 
$U$ with certain nondegenerate highest weight $\Lambda$, and for elliptic 
$\gamma$, we will compute explicitly 
$\mathcal{E}_{X,\gamma}(F=G\times_{K} E,t)$ and its Mellin transform in terms 
of the root systems. Note that when $\gamma=1$, $\mathcal{E}_{X,\gamma}(F_{d},t)$ is 
already computed by Bergeron-Venkatesh \cite{BV2013torsion} and by M\"{u}ller-Pfaff \cite{MR3128980} using the 
Plancherel formula for identity orbital integral. We here give a 
different approach via Bismut's formula as in 
\eqref{eq:6.2.8pl}.

Then in Subsection \ref{section4.3}, we apply these results to a sequence of flat 
vector bundles $\{F_{d}\}_{d\in\bN}$ on $X$ defined by a sequence of 
nondegenerate dominant weights $\Lambda=d\lambda+\lambda_{0}$. This 
way, we show that the Mellin transforms of the elliptic orbital 
integrals are exponential polynomials in $d$.

%%%%%%%%%%%%%%%%%%%%%%%%%%%%%%%%%%%%%%%%%%%%%%%%%%%%%%%%%%%%%%%%%%%%%%%%%%%%%%%%
\subsection{Estimates of elliptic orbital integrals for small time $t$}
Recall that $T$ is a maximal torus of $K$, $T_{U}$ is a maximal 
torus of $U$, and that $W(U,T_{U})$ denote the (analytic) Weyl group 
of $(U,T_{U})$. The positive root system $R^{+}(\ku,\kt_{U})$ is given in Subsection 
\ref{section5.3pl}. Recall that $P_{++}(U)$ is the set of dominant 
weights of $(U,T_U)$ with respect to $R^{+}(\ku,\kt_{U})$.

Let $(E,\rho^{E})$ be the irreducible 
unitary representation of $U$ associated with the highest weight 
$\Lambda\in P_{++}(U)$. We will prove our main result of this 
subsection and next subsection for this $(E,\rho^{E})$.

Our homogeneous flat vector bundle concerned here is given by 
$F=G\times_{K} E$. Let $\mathbf{D}^{X,F,2}$ denote the associated 
de Rham-Hodge Laplacian.

For $t>0$, if $\gamma\in G$ is semisimple, as in 
\eqref{eq:6.0.1kk}, set
\begin{equation}
	\cE_{X,\gamma}(F,t)=\mathrm{Tr_s}^{[\gamma]}\big[(N^{\Lambda^{\bullet}(T^*X)}-\frac{m}{2})\exp(-\frac{t}{2}\mathbf{D}^{X,F,2})\big].
	\label{eq:6.3.1pl}
\end{equation}
It is clear that $\cE_{X,\gamma}(F_{d},t)$ only depends on the conjugacy 
class $[\gamma]$ in $G$. 
If $\gamma=1$, we also write
\begin{equation}
	\cI_{X}(F,t)=\cE_{X,1}(F,t).
	\label{eq:6.3.2pl}
\end{equation}
In the sequel, we only consider the case of elliptic $\gamma$.

By \eqref{eq:6.2.8pl}, \eqref{eq:6.2.14ss20}, \eqref{eq:6.2.15s20}, 
if $\gamma=k\in K$, we have
\begin{equation}
	\begin{split}
		\cE_{X,\gamma}(F,t)&=\frac{1}{(2\pi t)^{p/2}} \exp\big(-2\pi^{2} 
		t|\Lambda+\rho_{\ku}|^{2}\big)\\
		&\int_{\kk(\gamma)} J_\gamma(Y^\kk_0)\mathrm{Tr_s}^{\Lambda^{\bullet}(\pp^*)}\big[(N^{\Lambda^{\bullet}(\pp^*)}-\frac{m}{2})\mathrm{Ad}(k)\exp(-i\mathrm{ad}(Y^\kk_0))\big]\\
		& 
		\qquad\cdot\mathrm{Tr}^{E}[\rho^{E}(k)\exp(-i\rho^{E}(Y^\kk_0))] e^{-|Y^\kk_0|^2/2t}\frac{dY^\kk_0}{(2\pi t)^{q/2}}.
	\end{split}
	\label{eq:7.1.5pl20}
\end{equation}

By \eqref{Jfunction}, we have the formula for 
$J_{\gamma}(Y^{\kk}_{0})$, $Y_{0}^{\kk}\in\kk(\gamma)$,
\begin{equation}\label{Jfunctionell}
	\begin{split}
		&J_{\gamma}(Y_0^{\kk})= \frac{\widehat{A}(i\mathrm{ad}(Y_0^\kk)|_{\pp(\gamma)})}{\widehat{A}(i\mathrm{ad}(Y_0^\kk)|_{\kk(\gamma)})}\\
		&\qquad\cdot\bigg[ \frac{1}{\det (1-\mathrm{Ad}(k))|_{\z^{\perp}(\gamma)}} \frac{\det (1-\exp(-i\mathrm{ad}(Y_0^\kk))\mathrm{Ad}(k))|_{\kk^{\perp}(\gamma)}}{\det (1-\exp(-i\mathrm{ad}(Y_0^\kk))\mathrm{Ad}(k))|_{\pp^{\perp}(\gamma)}}     \bigg]^{1/2}.
	\end{split}
\end{equation}

\begin{proposition}\label{prop:6.3.1est}
	For an elliptic element $\gamma\in G$, there exists a 
	constant $C_{\gamma}>0$ (depending on $\Lambda$) 
	such that for $t\in\, ]0,1]$
	\begin{equation}
		\begin{split}
			&|\sqrt{t}\cE_{X,\gamma}(F,t)|\leq C_{\gamma},\\
			&|(1+2t\frac{\partial}{\partial t})\cE_{X,\gamma}(F,t)|\leq 
			C_{\gamma}\sqrt{t}.
		\end{split}
		\label{eq:6.3.3pl}
	\end{equation}
	
	As $t\rightarrow 0$, $\cE_{X,\gamma}(E,t)$ has the asymptotic expansion in the form of 
	\begin{equation}
		\frac{1}{\sqrt{t}}\sum_{j=0}^{+\infty} a^\gamma_j t^j,
		\label{eq:6.3.4pl}
	\end{equation}
	with $a^\gamma_j\in\mathbb{C}$.
\end{proposition}
\begin{proof}
	If $\gamma$ is elliptic, up to a conjugation, we assume that 
	$\gamma=k\in T$. Thus the subgroup $H$ defined in Subsection 
	\ref{subsection4.1paris} is also a Cartan subgroup of $Z(\gamma)^{0}$, 
	then $\kb(\gamma)=\kb$. Let $\kb^{\perp}(\gamma)$ be the 
	orthogonal complementary space of $\kb(\gamma)$ in $\pp(\gamma)$, whose 
	dimension is $p-1$. Note that the similar estimates have been 
	proved in \cite[Theorem 
	4.4.1]{LIU2021109117}, here we only sketch a proof to \eqref{eq:6.3.3pl}. 
	
	By \eqref{eq:7.1.5pl20}, we have
	\begin{equation}
		\begin{split}
			\cE_{X,\gamma}(F,t)&=\frac{1}{(2\pi t)^{p/2}} \exp\big(-2\pi^{2} 
			t|\Lambda+\rho_{\ku}|^{2}\big)\\
			&\cdot\int_{\kk(k)} 
			J_k(\sqrt{t}Y^\kk_0)\mathrm{Tr_s}^{\Lambda^{\bullet}(\pp^*)}\big[(N^{\Lambda^{\bullet}(\pp^*)}-\frac{m}{2})\mathrm{Ad}(k)\exp(-i\mathrm{ad}(\sqrt{t}Y^\kk_0))\big]\\
			& 
			\qquad\quad\mathrm{Tr}^{E}[\rho^{E}(k)\exp(-i\rho^{E}(\sqrt{t}Y^\kk_0))] e^{-|Y^\kk_0|^2/2}\frac{dY^\kk_0}{(2\pi )^{q/2}},
		\end{split}
		\label{eq:7.2.5ss20jan}
	\end{equation}
	where the integral is rescaled by $\sqrt{t}$.
	
	In this proof, we denote by 
	$C$ or $c$ a positive constant independent of the variables $t$ and 
	$Y^{\kk}_{0}$. We use the symbol $\mathcal{O}_{\mathrm{ind}}$ to 
	denote the big-O convention which does not depend on $t$ and 
	$Y^{\kk}_{0}$.
	
	The same computations as in \cite[Eqs. (4.4.8) - 
	(4.4.10)]{LIU2021109117} shows that for $Y^{\kk}_{0}\in 
	\kt$,
	\begin{equation}
		\begin{split}
			&J_{k}(\sqrt{t}Y^\kk_0)=\frac{1}{\det(1-\mathrm{Ad}(k))|_{\pp^\perp(k)}}+\mathcal{O}_{\mathrm{ind}}(\sqrt{t}|Y^{\kk}_{0}|e^{C \sqrt{t}|Y^{\kk}_{0}|}),\\
			&\frac{1}{t^{(p-1)/2}}\mathrm{Tr_s}^{\Lambda^{\bullet}(\pp^*)}[\big(N^{\Lambda^{\bullet}(\pp^*)}-\frac{m}{2}\big)\rho^{\Lambda^{\bullet}(\pp^*)}(k)\exp(-i\rho^{\Lambda^{\bullet}(\pp^*)}(\sqrt{t}Y^\kk_0))]\\
			&=-\det(i\mathrm{ad}(Y^\kk_0))|_{\kb^\perp(k)}\det(1-\mathrm{Ad}(k))|_{\pp^\perp(k)}+\mathcal{O}_{\mathrm{ind}}(\sqrt{t}|Y^{\kk}_{0}|e^{C \sqrt{t}|Y^{\kk}_{0}|}).
		\end{split}
		\label{eq:8.4.19ugcs}
	\end{equation}
	Using the adjoint invariance, the further estimates on the above 
	quantities by a function in $|Y^{\kk}_{0}|$ hold for all $Y^{\kk}_{0}\in 
	\kk(k)$.
	
	It is clear that 
	\begin{equation}
		|\mathrm{Tr}^{E}[\rho^{E}(k)\exp(-i\rho^{E}(\sqrt{t}Y^\kk_0))]|\leq C\exp(C\sqrt{t}|Y^{\kk}_{0}|).
		\label{eq:6.3.7pl}
	\end{equation}
	
	Combining \eqref{eq:8.4.19ugcs} and \eqref{eq:6.3.7pl}, we see that 
	there exists a number $N\in\mathbb{N}$ big enough, if $t\in ]0,1]$
	\begin{equation}
		|\sqrt{t}\cE_{X,\gamma}(F,t)|\leq C'_{k}\int_{\kk(k)} 
		(1+|Y^{\kk}_{0}|)^{N}\exp(C|Y^{\kk}_{0}|-|Y^\kk_0|^2/2)dY^\kk_0.
		\label{eq:6.3.8pl}
	\end{equation}
	
	The second estimate in \eqref{eq:6.3.3pl} can be proved using the 
	same arguments as in \cite[Eqs. (4.4.24) - 
	(4.4.29)]{LIU2021109117}.
	
	The asymptotic expansion in \eqref{eq:6.3.4pl} is just a 
	consequence of \eqref{eq:6.3.3pl} and \eqref{eq:7.2.5ss20jan}. This completes the proof of our proposition.
\end{proof}

%%%%%%%%%%%%%%%%%%%%%%%%%%%%%%%%%%%%%%%%%%%%%%%%%%%%%%%%%%%%%%%%%%
\subsection{Elliptic orbital integrals for Hodge 
Laplacians}\label{sec7.2bath}
In this subsection, we explain how to use Bismut's formula 
\eqref{eq:6.2.8pl} to compute explicitly the expansion of 
$\cE_{X,\gamma}(F,t)$ in $t>0$ when
$\gamma\in G$ is elliptic. Then we study the corresponding Mellin 
transform. After 
conjugation, we may and we will 
assume that $\gamma=k\in T$. Then $T$ is also a maximal torus 
for $K(\gamma)^0$, and $\kb(\gamma)=\kb$.

Recall that $\omega^{Y_{\kb}(\gamma)}$, $\Omega^{\ku_{\kb}(\gamma)}$, 
$\Omega^{\ku_{\km}(\gamma)}$ are defined in Subsection 
\ref{subsection5.2ss}. Note that $\dim \ku^{\perp}_\kb(\gamma)=4l(\gamma)$. If 
$\nu\in\Lambda^{\bullet}(\ku^{\perp}_\kb(\gamma)^{*})$, let 
$[\nu]^{\mathrm{max}(\gamma)}\in\R$ be such that 
\begin{equation}
	\nu - 
	[\nu]^{\mathrm{max}(\gamma)}\frac{\omega^{Y_{\kb}(\gamma),2l(\gamma)}}{(2l(\gamma))!}
	\label{eq:7.5.19sskk20}
\end{equation}
is of degree strictly smaller than $4l(\gamma)$.

Recall that $-B(\cdot,\theta\cdot)$ is an Euclidean product on $\g$. 
Let $\kn^{\perp}(\gamma)$, $\bar{\kn}^{\perp}(\gamma)$ be the 
orthogonal spaces of $\kn(\gamma)$, $\bar{\kn}(\gamma)$ in $\kn$, 
$\bar{\kn}$ respectively. As $T$-modules, 
$\kn^{\perp}(\gamma)\simeq\bar{\kn}^{\perp}(\gamma)$.

Since $\kt\subset \kk(\gamma)\subset \kk$, then $R(\kk(\gamma),\kt)$ 
is a sub-root system of $R(\kk,\kt)$. Let $R^{+}(\kk(\gamma),\kt)$ be the 
positive root system for $(\kk(\gamma),\kt)$ induced by 
$R^{+}(\kk,\kt)$. We use the notation in Subsections \ref{section5.1ss20}, 
\ref{subsection5.2ss}. Then $\kt$ is a Cartan subalgebra for
$\kk_{\km}(\gamma)$, $\ku_{\km}(\gamma)$, $\km(\gamma)$. Let $R(\kk_{\km}(\gamma),\kt)$, $R(\ku_{\km}(\gamma),\kt)$ be the 
corresponding root systems. 

Similar to \eqref{eq:5.4.9bath}, we have the following disjoint union
\begin{equation}
	R(\ku_{\km}(\gamma),\kt)=R(\sqrt{-1}\pp_{\km}(\gamma),\kt)\cup 
	R(\kk_{\km}(\gamma),\kt).
\end{equation}
Since $R(\ku_{\km}(\gamma),\kt)\subset R(\ku_{\km},\kt)$, then by 
intersecting with $R^{+}(\ku_{\km},\kt)$, we get a positive root 
system $R^{+}(\ku_{\km}(\gamma),\kt)$. Moreover,
\begin{equation}
	R^{+}(\ku_{\km}(\gamma),\kt)=R^{+}(\sqrt{-1}\pp_{\km}(\gamma),\kt)\cup 
	R^{+}(\kk_{\km}(\gamma),\kt).
\end{equation}

Let $\mathrm{Vol}(K/T)$, $\mathrm{Vol}(U_{M}/T)$ be the Riemannian 
volumes 
of $K/T$, $U_{M}/T$ with respect to the restriction of $-B$ to $\kk$, 
$\ku_{\km}$ respectively. We have explicit formulae for them in terms 
of the roots, for example,
\begin{equation}
	\mathrm{Vol}(U_{M},T)=\Pi_{\alpha^{0}\in R^{+}(\ku_{\km},\kt)}\frac{1}{2\pi\langle 
	\alpha^{0},\rho_{\ku_{\km}}\rangle}.
\end{equation}

For $\gamma=k\in T$, set
\begin{equation}
	c_{G}(\gamma)=\frac{(-1)^{\frac{p-1}{2} 
	+1}\mathrm{Vol}(K(\gamma)^{0}/T)|W(U_{M}(\gamma)^{0}, T)|}{\mathrm{Vol}(U_{M}(\gamma)^{0}/T)|W(K(\gamma)^{0},T)|}\frac{1}{\det(1-\mathrm{Ad}(\gamma))|_{\kn^{\perp}(\gamma)}}.
	\label{eq:7.5.5kk20s}
\end{equation}
If $\gamma=1$, we denote
\begin{equation}
	c_{G}=c_{G}(1)=\frac{(-1)^{\frac{m-1}{2}+1}{\mathrm{Vol}}(K/T)|W(U_{M},T)|}{\mathrm{Vol}(U_{M}/T)|W(K,T)|}.
	\label{eq:eq:6.4.47ss20}
\end{equation}

We will use the same notation as in Subsections \ref{section5.3pl} \& \ref{section5.4}. In 
particular, $W_{u}$ is defined by \eqref{eq:5.3.6ss20} as a subset of 
$W(U,T_{U})$, and $W^{1}(\gamma)$ is defined by \eqref{eq:5.4.13ss20} as a 
subset of $W(U_{M},T)$. As explain in Remark 
\ref{rm:5.4.5s}, for $\omega\in W_{u}$, $\sigma\in W^{1}(\gamma)$, 
let $E^{\gamma}_{\omega,\sigma}$ denote the irreducible unitary 
representation of $Y=U_{M}(\gamma)^{0}$ or its finite central extension 
with highest weight 
$\sigma(\eta_{\omega}(\Lambda)+\rho_{\ku_{\km}})-\rho_{\y}$.

\begin{definition}
	For $j=0,1,\cdots, l(\gamma)$, $\omega\in 
	W_{u}$, $\sigma\in W^{1}(\gamma)$, set
	\begin{equation}
		\begin{split}
			Q^{\gamma}_{j,\omega,\sigma}(\Lambda)=&\frac{(-1)^{j}\beta(a_{1})^{2j}}{j!(2l(\gamma)-2j)! (8\pi^{2})^{j}}\dim E^{\gamma}_{\omega,\sigma} \\
			&\cdot\big[\omega^{Y_{\kb}(\gamma),2j}\langle 
			\omega(\Lambda+\rho_{\ku}), 
			\Omega^{\ku_{\km}(\gamma)}\rangle^{2l(\gamma)-2j}\big]^{\mathrm{max}(\gamma)}.
		\end{split}
		\label{eq:6.4.46ss20}
	\end{equation}
	
	In particular, if $l(\gamma)\geq 1$, we have
	\begin{equation}
		\begin{split}
			&Q^{\gamma}_{0,\omega,\sigma}(\Lambda)=\frac{1}{(2l)! }\dim 
			E^{\gamma}_{\omega,\sigma}\big[\langle 
			\omega(\Lambda+\rho_{\ku}), 
			\Omega^{\ku_{\km}(\gamma)}\rangle^{2l(\gamma)}\big]^{\mathrm{max}(\gamma)},\\
			&Q^{\gamma}_{l(\gamma),\omega,\sigma}(\Lambda)=\frac{(-1)^{l(\gamma)}\beta(a_{1})^{2l(\gamma)}(2l(\gamma)-1)!!}{(4\pi^{2})^{l(\gamma)}}\dim E^{\gamma}_{\omega,\sigma}.
		\end{split}
		\label{eq:6.4.46ss20copiesbath}
	\end{equation}
\end{definition}

Recall that $a_{1}\in \kb$ is such that $B(a_{1},a_{1})=1$. For 
$\omega\in W_{u}$, set
\begin{equation}
	b_{\Lambda,\omega}=\langle 
	\omega\cdot(\Lambda+\rho_{\ku}), \ii 
	a_{1}\rangle\in\R.
	\label{eq:7.2.9vogel}
\end{equation}
Then we have
\begin{equation}
	|\eta_{\omega}(\Lambda)+\rho_{\ku_{\km}}|^{2}-|\Lambda+\rho_{\ku}|^{2}=-b_{\Lambda,\omega}^{2}.
	\label{eq:6.4.8ssss20}
\end{equation}

Note that $\varphi_{\gamma}(\sigma,\eta_{\omega}(\Lambda))$ is 
defined in Definition \ref{def:5.4.3ss20}.
\begin{theorem}\label{thm:7.3.1bath}
	For $t>0$, we have the following identity
	\begin{equation}
		\cE_{X,\gamma}(F,t)=\frac{c_{G}(\gamma)}{\sqrt{2\pi 
		t}}\sum_{j=0}^{l(\gamma)} 
		t^{-j}\sum_{\genfrac{}{}{0pt}{2}{\omega\in W_{u}}{\sigma \in 
		W^{1}(\gamma)}}\varepsilon(\omega) 
		\varphi_{\gamma}(\sigma,\eta_{\omega}(\Lambda))e^{-2\pi^{2}tb_{\Lambda,\omega}^{2}}Q^{\gamma}_{j,\omega,\sigma}(\Lambda).
		\label{eq:6.4.48ss20bath}
	\end{equation}
\end{theorem}

\begin{remark}
	The formula \eqref{eq:6.4.48ss20bath} is compatible with the estimate 
	\eqref{eq:6.3.3pl}. For example, we take $\gamma=1$, then 
	$W^{1}(\gamma)$ reduce to $\{1\}$, the representation 
	$E^{\gamma}_{\omega,\sigma}$ is just $V_{\Lambda,\omega}$ 
	introduced in \eqref{eq:5.3.11pap}, and $l(\gamma)=l$, 
	$\varphi_{\gamma}(\sigma,\eta_{\omega}(\Lambda))=1$. Then we take the asymptotic 
	expansion of the right-hand side of \eqref{eq:6.4.48ss20} as 
	$t\rightarrow 0$, the coefficient of $t^{-l-1/2}$ is given by
	\begin{equation}
		\frac{c_{G}}{\sqrt{2\pi }}\sum_{\omega\in 
		W_{u}} \varepsilon(\omega) Q^{\gamma=1}_{l,\omega,1}(\Lambda).
		\label{eq:7.3.7ss20}
	\end{equation}
	By \eqref{eq:5.3.10dec19ss}, if $l\geq 1$, we get
	\begin{equation}
		\sum_{\omega\in W_{u}}\varepsilon(\omega)\dim 
		V_{\Lambda,\omega}=\mathrm{Tr_{s}}^{\Lambda^{\bullet}(\kn^{*}_{\C})}[1]\dim E=0.
		\label{eq:7.3.8ss20k}
	\end{equation}
	Then by \eqref{eq:6.4.46ss20copiesbath} and \eqref{eq:7.3.8ss20k}, 
	the quantity in \eqref{eq:7.3.7ss20} is $0$ (provided $l\geq 1$).
\end{remark}

Before proving Theorem \ref{thm:7.3.1bath}, we need some preparation work. 

\begin{definition}\label{def:6.4.1ss20}
	For $y\in\kt$, put
	\begin{equation}
		\begin{split}
			&\pi_{\ku_{\km}(\gamma)/\kt}(y)=\prod_{\alpha^{0}\in 
			R^{+}(\ku_{\km}(\gamma),\kt)}\langle 2\pi\ii \alpha^{0}, y\rangle.\\
			&\pi_{\ii\pp_{\km}(\gamma)/\kt}(y)=\prod_{\alpha^{0}\in 
			R^{+}(\ii\pp_{\km}(\gamma),\kt)}\langle 2\pi\ii \alpha^{0}, y\rangle.\\
			&\pi_{\kk_{\km}(\gamma)/\kt}(y)=\prod_{\alpha^{0}\in 
			R^{+}(\kk_{\km}(\gamma),\kt)}\langle 2\pi\ii \alpha^{0}, y\rangle.\\
		\end{split}
		\label{eq:6.4.11final}
	\end{equation}

	For $y\in\kt$, put
	\begin{equation}
		\begin{split}
			&\sigma_{\ku_{\km}(\gamma)/\kt}(y)=\prod_{\alpha^{0}\in 
			R^{+}(\ku_{\km}(\gamma),\kt)}\big(\exp(\langle \pi\ii \alpha^{0}, 
			y\rangle)-\exp(-\langle \pi\ii \alpha^{0}, y\rangle)\big).\\
			&\sigma_{\ii\pp_{\km}(\gamma)/\kt}(y)=\prod_{\alpha^{0}\in 
			R^{+}(\ii\pp_{\km}(\gamma),\kt)}\big(\exp(\langle \pi\ii \alpha^{0}, 
			y\rangle)-\exp(-\langle \pi\ii \alpha^{0}, y\rangle)\big).\\
			&\sigma_{\kk_{\km}(\gamma)/\kt}(y)=\prod_{\alpha^{0}\in 
			R^{+}(\kk_{\km}(\gamma)t,\kt)}\big(\exp(\langle \pi\ii \alpha^{0}, 
			y\rangle)-\exp(-\langle \pi\ii \alpha^{0}, y\rangle)\big).\\
		\end{split}
		\label{eq:6.4.12final}
	\end{equation}
	
	We can always extend analytically the above functions to 
	$y\in\kt_{\C}$. If $\gamma=1$, the above functions become 
	$\pi_{\ku_{\km}/\kt}(y)$, 
	$\pi_{\sqrt{-1}\pp_{\km}/\kt}(y)$, $\pi_{\kk_{\km}/\kt}(y)$, $\sigma_{\ku_{\km}/\kt}(y)$, $\sigma_{\ii\pp_{\km}/\kt}(y)$, $\sigma_{\kk_{\km}/\kt}(y)$. 
	
	If the adjoint action of $T$ preserves certain orthogonal 
	splittings of $\ku_{\km}$, $\ku_{\km}(\gamma)$, etc, so that we 
	have the corresponding splitting of the root systems, then we can 
	also define the associated $\pi$-function or $\sigma$-function as 
	above. 
\end{definition}

It is clear that if $y\in\kt_{\C}$, 
\begin{equation}
	\begin{split}
		&\pi_{\ku_{\km}(\gamma)/\kt}(y)=\pi_{\ii\pp_{\km}(\gamma)/\kt}(y)\pi_{\kk_{\km}(\gamma)/\kt}(y), \\ 
		&\sigma_{\ku_{\km}(\gamma)/\kt}(y)=\sigma_{\ii\pp_{\km}(\gamma)/\kt}(y)\sigma_{\kk_{\km}/\kt}(y)
	\end{split}
	\label{eq:6.4.13final}
\end{equation}

Set
\begin{equation}
	\begin{split}		
		&\kk^{\prime}_{\km}(\gamma)=\kk^{\perp}(\gamma)\cap\kk_{\km},\; 
		\pp^{\prime}_{\km}(\gamma)=\pp^{\perp}(\gamma)\cap\pp_{\km};\\
		&\kk^{\prime\prime}_{\km}(\gamma)=\kk^{\perp}(\gamma)\cap\kk^{\perp}(\kb),\; 
		\pp^{\prime\prime}_{\km}(\gamma)=\pp^{\perp}(\gamma)\cap\pp^{\perp}(\kb).
	\end{split}
\end{equation}
Let $\km^{\perp}(\gamma)$ be the orthogonal space of $\km(\gamma)$ in 
$\km$ with respect to $B$. Then
\begin{equation}
	\km^{\perp}(\gamma)=\pp^{\prime}_{\km}(\gamma)\oplus 
	\kk^{\prime}_{\km}(\gamma).
\end{equation}
We also have
\begin{equation}
	\kk_{\km}=\kk_{\km}(\gamma)\oplus \kk^{\prime}_{\km}(\gamma),\; 
	\pp_{\km}=\pp_{\km}(\gamma)\oplus \pp^{\prime}_{\km}(\gamma), 
	\label{eq:6.6.11ss20}
\end{equation}
and
\begin{equation}
	\kk^{\perp}(\gamma)=\kk^{\prime}_{\km}(\gamma)\oplus 
	\kk^{\prime\prime}_{\km}(\gamma),\; 
	\pp^{\perp}(\gamma)=\pp^{\prime}_{\km}(\gamma)\oplus 
	\pp^{\prime\prime}_{\km}(\gamma). 
	\label{eq:6.6.12ss20}
\end{equation}

Set
\begin{equation}
	\ku^{\perp}_{\km}(\gamma)=\ii\pp'_{\km}(\gamma)\oplus \kk'_{\km}(\gamma).
\end{equation}
Then it is the orthogonal space of $\ku_{\km}(\gamma)$ in $\ku_{\km}$ 
with respect to $B$.

\begin{lemma}\label{lm:7.5.7skk}
	The following spaces are isomorphic to each other as modules of 
	$T$ by the adjoint actions,
	\begin{equation}
		\kn^{\perp}(\gamma)\simeq\bar{\kn}^{\perp}(\gamma)\simeq 
		\kk^{\prime\prime}_{\km}(\gamma) \simeq 
		\pp^{\prime\prime}_{\km}(\gamma).
	\end{equation}
\end{lemma}
\begin{proof}
	Note that
	\begin{equation}
		\dim \kn=\dim\kk -\dim \kk_{\km},\; \dim \kn(\gamma)=\dim 
		\kk(\gamma)-\dim \kk_{\km}(\gamma).
	\end{equation}
	Together with the splittings \eqref{eq:6.6.11ss20}, 
	\eqref{eq:6.6.12ss20}, we get
	\begin{equation}
		\dim \kk^{\prime\prime}_{\km}(\gamma)=\dim \kn^{\perp}(\gamma).
	\end{equation}
	Similarly, $\dim \pp^{\prime\prime}_{\km}(\gamma)=\dim 
	\kn^{\perp}(\gamma)$.

	If $f\in\kn^{\perp}(\gamma)$, then $f+\theta(f)\in\kk$, we can verify 
	directly that $f+\theta(f)\in \kk^{\prime\prime}_{\km}(\gamma)$. Then 
	the map $f\in\kn^{\perp}(\gamma)\mapsto f+\theta(f)\in 
	\kk^{\prime\prime}_{\km}(\gamma)$ defines an isomorphisms of 
	$T$-modules. Similar for $\kn^{\perp}(\gamma)\simeq 
	\pp^{\prime\prime}_{\km}(\gamma)$.
\end{proof}

Since $\gamma=k\in T$, let $y_0\in \kt$ be such that $\exp(y_0)=\gamma$. 
Note that $y_{0}$ is not unique.
\begin{lemma}\label{lm:5.4.1}
	If $y\in \kt$ is regular with respect to $R(\kk_{\km}(\gamma),\kt)$, then we have
	\begin{equation}
		\begin{split}		
			&J_{\gamma}(y)\mathrm{Tr_{s}}^{\Lambda^{\bullet}(\pp^{*})}[(N^{\Lambda^{\bullet}(\pp^{*})}-\frac{m}{2})\mathrm{Ad}(k)\exp(-\mathrm{i}\mathrm{ad}(y))]\\
			=&\frac{(-1)^{\dim\pp_{\km}(\gamma)/2 +1}}{\det(1-\mathrm{Ad}(k))|_{\kn^{\perp}(\gamma)}}\mathrm{Tr_{s}}^{\Lambda^{\bullet}(\kn^{*}_{\C})}[e^{-i\mathrm{ad}(y)}\mathrm{Ad}(k)]\\
			&\cdot\frac{\pi_{\ii\pp_{\km}(\gamma)/\kt}(iy)}{\pi_{\kk_{\km}(\gamma)/\kt}(iy)}\frac{\sigma_{\ku_{\km}(\gamma)/\kt}(iy)\sigma_{\ku^{\perp}_{\km}(\gamma)/\kt}(-iy+y_{0})}{\sigma_{\ku^{\perp}_{\km}(\gamma)/\kt}(y_{0})}.
		\end{split}
		\label{eq:5.4.5bonn}
	\end{equation}
\end{lemma}
\begin{proof}
	Using \eqref{eq:5.4.24qq20}, \eqref{eq:6.6.12ss20} and Lemma 
	\ref{lm:7.5.7skk}, we get that for $y\in\kt$, 
	\begin{equation}
		\begin{split}
			&\bigg[\frac{1}{\det 
			(1-\mathrm{Ad}(k))|_{\z^\perp(\gamma)}}\frac{\det(1-e^{-\mathrm{i}\mathrm{ad}(y)}\mathrm{Ad}(k))|_{\kk^\perp(\gamma)}}{\det(1-e^{-\mathrm{i}\mathrm{ad}(y)}\mathrm{Ad}(k))|_{\pp^\perp(\gamma)}}\bigg]^{1/2}\\
			&=\frac{(-1)^{\frac{\dim 
			\pp^{\prime}_{\km}(\gamma)}{2}}}{\det(1-\mathrm{Ad}(k))|_{\kn^{\perp}(\gamma)}}\frac{1}{\sigma_{\ku^{\perp}_{\km}(\gamma)/\kt}(y_{0})}\frac{\sigma_{\kk^{\prime}_{\km}(\gamma)/\kt}(-iy+y_{0})}{\sigma_{\sqrt{-1}\pp^{\prime}_{\km}(\gamma)/\kt}(-iy+y_{0})}.
		\end{split}
		\label{eq:7.5.20kks}
	\end{equation}
	
	Recall that in Subsection \ref{section5.1ss20}, as $K_{M}$-modules, we have the following isomorphism
	\begin{equation}\label{eq:7.3.22vogel}
		\pp\simeq \kb\oplus \pp_{\km}\oplus \kn.
	\end{equation}
	Note that 
	\begin{equation}
		\mathrm{Ad}(k)=e^{\mathrm{ad}(y_{0})}.
	\end{equation}
	If $y\in\kt$, when acting on $\pp$, we have
	\begin{equation}
		\mathrm{Ad}(k)\exp(-i\mathrm{ad}(y))=\exp(\mathrm{ad}(-iy+y_{0})).
	\end{equation}
	
	Note that $\dim \kb=1$. Then for $y\in \kt$, we get
	\begin{equation}\label{eq:6.4.22sasy}
		\begin{split}
			&\mathrm{Tr_s}^{\Lambda^{\bullet}(\pp^*)}\big[(N^{\Lambda^{\bullet}(\pp^*)}-\frac{m}{2})\mathrm{Ad}(k)\exp(-\mathrm{i}\mathrm{ad}(y))\big]\\
			&=-\mathrm{Tr_{s}}^{\Lambda^{\bullet}(
			\pp_{\km}^{*})}[\mathrm{Ad}(k)e^{-i\mathrm{ad}(y)}]\mathrm{Tr_{s}}^{\Lambda^{\bullet}(\kn^{*}_{\C})}[\mathrm{Ad}(k)e^{-i\mathrm{ad}(y)}]\\
			&=-\det(1- 
			\mathrm{Ad}(k^{-1})e^{i\mathrm{ad}(y)})|_{\pp_{\km}}\mathrm{Tr_{s}}^{\Lambda^{\bullet}(\kn^{*}_{\C})}[\mathrm{Ad}(k)e^{-i\mathrm{ad}(y)}],
		\end{split}
	\end{equation}
	where we have the identity
	\begin{equation}
		\det(1- 
		\mathrm{Ad}(k^{-1})e^{i\mathrm{ad}(y)})|_{\pp_{\km}}=(-1)^{(\dim\pp_{\km})/2}\sigma_{\ii\pp'_{\km}(\gamma)/\kt}(-iy+y_{0})^2\sigma_{\ii\pp_{\km}(\gamma)/\kt}(iy)^2.
		\label{eq:6.4.23dec19s}
	\end{equation}
	
	Note that analogue to \eqref{eq:7.3.22vogel}, we have 
	$\pp(\gamma)\simeq \kb\oplus\pp_{\km}(\gamma)\oplus\kn(\gamma)$, using 
	\cite[Eq.(7.5.24)]{bismut2011hypoelliptic}, if $y\in\kt$, we have
	\begin{equation}
		\begin{split}
			&\widehat{A}(i\mathrm{ad}(y)|_{\ii\pp(\gamma)})=\frac{\pi_{\ii\pp_{\km}(\gamma)/\kt}(iy)}{\sigma_{\ii\pp_{\km}(\gamma)/\kt}(iy)}\widehat{A}(i\mathrm{ad}(y)|_{\kn(\gamma)}),\\
			&\widehat{A}(i\mathrm{ad}(y)|_{\kk(\gamma)})=\frac{\pi_{\kk(\gamma)/\kt}(iy)}{\sigma_{\kk(\gamma)/\kt}(iy)}=\frac{\pi_{\kk_{\km}(\gamma)/\kt}(iy)}{\sigma_{\kk_{\km}(\gamma)/\kt}(iy)}\widehat{A}(i\mathrm{ad}(y)|_{\kn(\gamma)}).
		\end{split}
		\label{eq:6.4.24kk19dss}
	\end{equation}
	
	Combining \eqref{Jfunctionell}, \eqref{eq:7.5.20kks} and 
	\eqref{eq:6.4.22sasy} - \eqref{eq:6.4.24kk19dss}, we get \eqref{eq:5.4.5bonn}. 
\end{proof}

Now we prove Theorem \ref{thm:7.3.1bath}.
\begin{proof}[Proof to Theorem \ref{thm:7.3.1bath}]
	Put
	\begin{equation}
		\begin{split}
			F_{\gamma}(\Lambda,t)=\frac{1}{(2\pi t)^{p/2}}\int_{\kk(\gamma)} 
			J_\gamma(Y^\kk_0)\mathrm{Tr_s}^{\Lambda^{\bullet}(\pp^*)}\big[(N^{\Lambda^{\bullet}(\pp^*)}-\frac{m}{2})\mathrm{Ad}(k)e^{-i\mathrm{ad}(Y^\kk_0)}\big]&\\
			\cdot\mathrm{Tr}^{E}[\rho^{E}(k)e^{-i\rho^{E}(Y^\kk_0)}] e^{-|Y^\kk_0|^2/2t}\frac{dY^\kk_0}{(2\pi t)^{q/2}}.&
		\end{split}
		\label{eq:5.4.6bonn}
	\end{equation}
	By \eqref{eq:7.1.5pl20}, we have
	\begin{equation}
		\cE_{X,\gamma}(F,t)=\exp\big(-2\pi^{2} t 
		|\Lambda+\rho_{\ku}|^{2}\big)F_{\gamma}(\Lambda,t).
		\label{eq:7.5.22decjan19}
	\end{equation}

	Recall that $r=p+q=\dim_{\R}\z(\gamma)$. By Weyl integration formula, then
	\begin{equation}\label{eq:7.3.29vogel}
		\begin{split}
			F_{\gamma}(\Lambda,t)=&\frac{\mathrm{Vol}(K(\gamma)^{0}/T)}{(2\pi 
			t)^{r/2} |W(K(\gamma)^{0},T)|}\\
			&\cdot\int_\kt |\pi_{\kk(\gamma)/\kt}(y)|^2 
			J_\gamma(y)\mathrm{Tr_s}^{\Lambda^{\bullet}(\pp^*)}\big[(N^{\Lambda^{\bullet}(\pp^*)}-\frac{m}{2})\mathrm{Ad}(k)e^{-i\mathrm{ad}(y)})\big]\\
			&\qquad\qquad\mathrm{Tr}^{E}[\rho^{E}(k)\exp(-i\rho^{E}(y))] e^{-|y|^2/2t}dy.
		\end{split}
	\end{equation}
	
	Recall that $l(\gamma)=\frac{1}{2}\dim \kn(\gamma)$. We can verify directly that if $y\in\kt$, 
	\begin{equation}
		\pi_{\kk(\gamma)/\kt}(iy)^{2}=(-1)^{l(\gamma)}\pi_{\kk_{\km}(\gamma)/\kt}(iy)^{2}\det(i\mathrm{ad}(y))|_{\kn(\gamma)_{\C}}.
		\label{eq:7.5.25ss20ks}
	\end{equation}

	Moreover, if $y\in\kt$ is such that $\pi_{\ku_{\km}(\gamma)/\kt}(y)\neq 0$,
	\begin{equation}\label{eq:4.3.6asy}
		\frac{|\pi_{\kk(\gamma)/\kt}(y)|^2}{|\pi_{\ku_{\km}(\gamma)/\kt}(y)|^2}=\frac{\pi_{\kk(\gamma)/\kt}(iy)^2}{\pi_{\ku_{\km}(\gamma)/\kt}(iy)^2}.
	\end{equation}
	Then by Lemma \ref{lm:5.4.1} and \eqref{eq:7.5.5kk20s}, 
	\eqref{eq:6.4.24kk19dss}, \eqref{eq:7.5.25ss20ks}, 
	we get
	\begin{equation}
		\begin{split}
			&\frac{|\pi_{\kk(\gamma)/\kt}(y)|^2}{|\pi_{\ku_{\km}(\gamma)/\kt}(y)|^2}J_{\gamma}(y)\mathrm{Tr_{s}}^{\Lambda^{\bullet}(\pp^{*})}[(N^{\Lambda^{\bullet}(\pp^{*})}-\frac{m}{2})\mathrm{Ad}(k)\exp(-\mathrm{i}\mathrm{ad}(y))]\\
			&=\frac{(-1)^{l(\gamma)+\dim\pp_{\km}(\gamma)/2 +1}}{\det(1-\mathrm{Ad}(k))|_{\kn^{\perp}(\gamma)}}\mathrm{Tr_{s}}^{\Lambda^{\bullet}(\kn^{*}_{\C})}[e^{-i\mathrm{ad}(y)}\mathrm{Ad}(k)]\\
			&\cdot \det(i\mathrm{ad}(y))|_{\kn(\gamma)_{\C}}\widehat{A}^{-1}(i\mathrm{ad}(y)|_{\ku_{\km}(\gamma)}) \big[\frac{\det (1-
			e^{-i\mathrm{ad}(y)}\mathrm{Ad}(k))|_{\ku^{\perp}_{\km}(\gamma)}}{\det(1-\mathrm{Ad}(k))
			_{\ku^{\perp}_{\km}(\gamma)}}\big]^{\frac{1}{2}}.
		\end{split}
	\end{equation}
	
	Note that we have the even number
	\begin{equation}
		p-1=\dim\pp_{\km}(\gamma)+2l(\gamma).
		\label{eq:7.3.33bath}
	\end{equation}	
	Now we can rewrite \eqref{eq:7.3.29vogel} as follows,
	\begin{equation}\label{eq:7.3.34bath}
		\begin{split}
			F_{\gamma}(\Lambda,t)&=\frac{(-1)^{\frac{p-1}{2}+1}\mathrm{Vol}(K(\gamma)^{0}/T)}{(2\pi 
			t)^{r/2} |W(K(\gamma)^{0},T)|}\frac{1}{\det(1-\mathrm{Ad}(k))|_{\kn^{\perp}(\gamma)}}\\
			&\cdot\int_\kt |\pi_{\ku_{\km}(\gamma)/\kt}(y)|^2 
			\det(i\mathrm{ad}(y))|_{\kn(\gamma)_{\C}}\cdot\widehat{A}^{-1}(i\mathrm{ad}(y)|_{\ku_{\km}(\gamma)})\\
			&\qquad\cdot\big[\frac{\det (1-
			e^{-i\mathrm{ad}(y)}\mathrm{Ad}(k))|_{\ku^{\perp}_{\km}(\gamma)}}{\det(1-\mathrm{Ad}(k))
			_{\ku^{\perp}_{\km}(\gamma)}}\big]^{\frac{1}{2}}\\
			&\qquad\cdot\mathrm{Tr_{s}}^{\Lambda^{\bullet}(\kn^{*}_{\C})\otimes 
			E}[e^{-i\rho^{\Lambda^{\bullet}(\kn^{*}_{\C})\otimes 
			E}(y)}\rho^{\Lambda^{\bullet}(\kn^{*}_{\C})\otimes E}(k)] e^{-|y|^2/2t}dy.
		\end{split}
	\end{equation}
	
	Note that the function in $y\in\kt$
	\begin{equation}
		\begin{split}
			\det(i\mathrm{ad}(y))|_{\kn(\gamma)_{\C}}\cdot\widehat{A}^{-1}(i\mathrm{ad}(y)|_{\ku_{\km}(\gamma)}) \big[\frac{\det (1-
			e^{-i\mathrm{ad}(y)}\mathrm{Ad}(k))|_{\ku^{\perp}_{\km}(\gamma)}}{\det(1-\mathrm{Ad}(k))
			_{\ku^{\perp}_{\km}(\gamma)}}\big]^{\frac{1}{2}}&\\
			\cdot\mathrm{Tr_{s}}^{\Lambda^{\bullet}(\kn^{*}_{\C})\otimes 
			E}[e^{-i\rho^{\Lambda^{\bullet}(\kn^{*}_{\C})\otimes 
			E}(y)}\rho^{\Lambda^{\bullet}(\kn^{*}_{\C})\otimes E}(k)]&
		\end{split}
		\label{eq:6.4.29dec19s}
	\end{equation}
	can be extended directly to a $U_{M}(\gamma)^{0}$-invariant function in 
	$y\in\ku_{\km}(\gamma)$. Since $\kt$ is a Cartan subalgebra of 
	$\ku_{\km}(\gamma)$, 
	we can apply the Weyl integration formula for the pair 
	$(\ku_{\km}(\gamma),\kt)$, we get
	\begin{equation}
		\begin{split}
			F_{\gamma}(\Lambda,t)=&\frac{c_{G}(\gamma)}{(2\pi 
			t)^{r/2} 
			}\int_{y\in 
			\ku_{\km}(\gamma)}\det(i\mathrm{ad}(y))|_{\kn(\gamma)_{\C}}\cdot\widehat{A}^{-1}(i\mathrm{ad}(y)|_{\ku_{\km}(\gamma)})\\
			&\qquad\qquad\cdot\big[\frac{\det (1-
			e^{-i\mathrm{ad}(y)}\mathrm{Ad}(k))|_{\ku^{\perp}_{\km}(\gamma)}}{\det(1-\mathrm{Ad}(k))
			_{\ku^{\perp}_{\km}(\gamma)}}\big]^{\frac{1}{2}}\\
			&\qquad\qquad\cdot\mathrm{Tr_{s}}^{\Lambda^{\bullet}(\kn^{*}_{\C})\otimes 
			E}[e^{-i\rho^{\Lambda^{\bullet}(\kn^{*}_{\C})\otimes 
			E}(y)}\rho^{\Lambda^{\bullet}(\kn^{*}_{\C})\otimes E}(k)]e^{-|y|^2/2t}dy.
		\end{split}
		\label{eq:7.5.26ss20}
	\end{equation}
	The constant $c_{G}(\gamma)$ is defined by \eqref{eq:7.5.5kk20s}.

	Note that
	\begin{equation}
		r=\dim\ku_{\km}(\gamma)+4l(\gamma)+1.
		\label{eq:7.5.33decss19s}
	\end{equation}	
	If $y\in \ku_{\km}(\gamma)$, then
	\begin{equation}
		B(y,\frac{\Omega^{\ku_{\km}(\gamma)}}{2\pi})\in\Lambda^{2}(\ku^{\perp}_\kb(\gamma)^{*}).
		\label{eq:7.4.31dec19}
	\end{equation}
	If $y\in\ku_{\km}(\gamma)$, by \cite[Eq. (7-27)]{Shen_2016}, we have
	\begin{equation}
		\frac{\det(i\mathrm{ad}(y))|_{\kn(\gamma)_{\C}}}{(2\pi 
		t)^{2l(\gamma)}}=[\exp\big(\frac{1}{t}B(y,\frac{\Omega^{\ku_{\km}(\gamma)}}{2\pi})\big)]^{\mathrm{max}(\gamma)}.
		\label{eq:7.5.32decss19s}
	\end{equation}

	Combining \eqref{eq:7.5.26ss20} - \eqref{eq:7.5.32decss19s}, we get
	\begin{equation}
		\begin{split}
			&F_{\gamma}(\Lambda,t)=\frac{c_{G}(\gamma)}{\sqrt{2\pi t}}\bigg[\int_{y\in 
			\ku_{\km}(\gamma)}\widehat{A}^{-1}(i\mathrm{ad}(y)|_{\ku_{\km}(\gamma)})\big[\frac{\det (1-
			e^{-i\mathrm{ad}(y)}\mathrm{Ad}(k))|_{\ku^{\perp}_{\km}(\gamma)}}{\det(1-\mathrm{Ad}(k))
			_{\ku^{\perp}_{\km}(\gamma)}}\big]^{\frac{1}{2}}\\
			&\cdot\mathrm{Tr_{s}}^{\Lambda^{\bullet}(\kn^{*}_{\C})\otimes 
			E}[\rho^{\Lambda^{\bullet}(\kn^{*}_{\C})\otimes 
			E}(e^{-iy}k)]e^{\frac{1}{t}B(y,\frac{\Omega^{\ku_{\km}(\gamma)}}{2\pi})-|y|^2/2t}\frac{dy}{(2\pi t)^{\dim \ku_{\km}(\gamma)/2}}\bigg]^{\mathrm{max}(\gamma)}.
		\end{split}
		\label{eq:7.5.32kk20}
	\end{equation}

	By \eqref{eq:5.2.21dec19}  if $y\in\ku_{\km}(\gamma)$, then
	\begin{equation}
		B(y,\frac{\Omega^{\ku_{\km}(\gamma)}}{2\pi})-\frac{|y|^{2}}{2}=\frac{1}{2}B(y+\frac{\Omega^{\ku_{\km}(\gamma)}}{2\pi},y+\frac{\Omega^{\ku_{\km}(\gamma)}}{2\pi})-\frac{\beta(a_{1})^{2}}{8\pi^{2}}\omega^{Y_{\kb}(\gamma),2}.
		\label{eq:6.4.36dec19s}
	\end{equation}

	Let $\Delta^{\ku_{\km}(\gamma)}$ be the standard negative Laplace 
	operator on the Euclidean space $(\ku_{\km}(\gamma), 
	-B|_{\ku_{\km}(\gamma)})$. Then by considering the heat kernel of 
	$-\Delta^{\ku_{\km}(\gamma)}$, we can rewrite \eqref{eq:7.5.32kk20} as follows,
	\begin{equation}
		\begin{split}
			F_{\gamma}(\Lambda,t)=&\frac{c_{G}(\gamma)}{\sqrt{2\pi t}}\bigg[\exp(-\frac{\beta(a_{1})^{2} 
			\omega^{Y_{\kb}(\gamma),2}}{8\pi^{2}t})\\
			&\exp(\frac{t}{2}\Delta^{\ku_{\km}(\gamma)})\bigg\{\widehat{A}^{-1}(i\mathrm{ad}(y)|_{\ku_{\km}(\gamma)})\big[\frac{\det (1-
			e^{-i\mathrm{ad}(y)}\mathrm{Ad}(k))|_{\ku^{\perp}_{\km}(\gamma)}}{\det(1-\mathrm{Ad}(k))
			_{\ku^{\perp}_{\km}(\gamma)}}\big]^{\frac{1}{2}}\\
			&\mathrm{Tr_{s}}^{\Lambda^{\bullet}(\kn^{*}_{\C})\otimes 
			E}\left[\rho^{\Lambda^{\bullet}(\kn^{*}_{\C})\otimes 
			E}(e^{-iy}k)\right]\bigg\}|_{y=-\frac{\Omega^{\ku_{\km}(\gamma)}}{2\pi}}\bigg]^{\mathrm{max}(\gamma)}.
		\end{split}
		\label{eq:7.5.33kk20}
	\end{equation}
	
	Recall that $V_{\Lambda,\omega}$ is an irreducible unitary 
	representation of $U_{M}$ with highest weight 
	$\eta_{\omega}(\Lambda)$. By \eqref{eq:5.3.10dec19ss}, 
	for $y\in\ku_{\km}(\gamma)$, then
	\begin{equation}
		\mathrm{Tr_{s}}^{\Lambda^{\bullet}(\kn^{*}_{\C})\otimes 
		E}[\rho^{\Lambda^{\bullet}(\kn^{*}_{\C})\otimes 
		E}(e^{-iy}k)]=\sum_{\omega\in 
		W_{u}}\varepsilon(\omega)\mathrm{Tr}^{V_{\Lambda,\omega}}\big[\rho^{V_{\Lambda,\omega}}(e^{-iy}k)\big].
		\label{eq:6.4.39sssks}
	\end{equation}
	Then
	we apply the generalized Kirillov formula \eqref{eq:haha20imp} to 
	each term in the right-hand side of \eqref{eq:6.4.39sssks}, we 
	conclude that, for $\omega\in W_{u}$, the function in $y\in\ku_{\km}(\gamma)$
	\begin{equation}
		\widehat{A}^{-1}(i\mathrm{ad}(y)|_{\ku_{\km}(\gamma)})\big[\frac{\det (1-
		e^{-i\mathrm{ad}(y)}\mathrm{Ad}(k))|_{\ku^{\perp}_{\km}(\gamma)}}{\det(1-\mathrm{Ad}(k))
		_{\ku^{\perp}_{\km}(\gamma)}}\big]^{\frac{1}{2}}\mathrm{Tr}^{V_{\Lambda,\omega}}[\rho^{V_{\Lambda,\omega}}(e^{-iy}k)]
		\label{eq:7.5.36conf}
	\end{equation}
	is an eigenfunction of $\Delta^{\ku_{\km}(\gamma)}$ associated with the 
	eigenvalue 
	$4\pi^{2}|\eta_{\omega}(\Lambda)+\rho_{\ku_{\km}}|^{2}$.
	Then the heat operator 
	$\exp\big(\frac{t}{2}\Delta^{\ku_{\km}(\gamma)}\big)$ acts on the
	function \eqref{eq:7.5.36conf} as a scalar 
	$e^{2\pi^{2}t|\eta_{\omega}(\Lambda)+\rho_{\ku_{\km}}|^{2}}$. By \eqref{eq:5.3.11final}, \eqref{eq:5.3.12decfinal19}, for 
	$\omega\in W_{u}$, we get
	\begin{equation}
		\eta_{\omega}(\Lambda)+\rho_{\ku_{\km}}=P_{0}(\omega(\Lambda+\rho_{\ku})).
		\label{eq:5.4.7jan20}
	\end{equation}
	Combing the above computation with the term $e^{-2\pi^{2} t 
	|\Lambda+\rho_{\ku}|^{2}}$ in \eqref{eq:7.5.22decjan19}, by 
	\eqref{eq:6.4.8ssss20}, we get the factor $e^{-2\pi^{2}t 
	b^{2}_{\Lambda,\omega}}$ in \eqref{eq:6.4.48ss20bath}.
	
	Now we deal with the main part in \eqref{eq:7.5.33kk20} 
	after removing the heat operator 
	$\exp\big(\frac{t}{2}\Delta^{\ku_{\km}(\gamma)}\big)$. We will use the same notation as in Subsection \ref{section5.4}. The orbit 
	$\mathcal{O}^{\gamma}_{\sigma(\eta_{\omega}(\Lambda)+
	\rho_{\ku_{\km}})}$ is defined in \eqref{eq:5.4.17qq20} equipped 
	with a Liouville measure $d\mu^{\gamma}_{\sigma}$. We claim the following identity,
	\begin{equation}
		\begin{split}
			&\bigg[\exp(-\frac{\beta(a_{1})^{2} 
			\omega^{Y_{\kb}(\gamma),2}}{8\pi^{2}t})\bigg\{\widehat{A}^{-1}(i\mathrm{ad}(y)|_{\ku_{\km}(\gamma)})\big[\frac{\det (1-
			e^{-i\mathrm{ad}(y)}\mathrm{Ad}(k))|_{\ku^{\perp}_{\km}(\gamma)}}{\det(1-\mathrm{Ad}(k))
			_{\ku^{\perp}_{\km}(\gamma)}}\big]^{\frac{1}{2}}\\
			&\mathrm{Tr_{s}}^{V_{\Lambda,\omega}}[\rho^{V_{\Lambda,\omega}}(e^{-iy}k)]\bigg\}|_{y=-\frac{\Omega^{\ku_{\km}(\gamma)}}{2\pi}}\bigg]^{\mathrm{max}(\gamma)}\\
			=& \sum_{\sigma\in 
			W^{1}(\gamma)}\varphi_{\gamma}(\sigma,\eta_{\omega}(\Lambda))\cdot\dim E^{\gamma}_{\omega,\sigma}\\
			&\qquad\cdot\bigg[\exp\left(-\frac{\beta(a_{1})^{2} 
			\omega^{Y_{\kb}(\gamma),2}}{8\pi^{2}t}-\langle 
			\sigma(\eta_{\omega}(\Lambda)+
			\rho_{\ku_{\km}}), 
			\Omega^{\ku_{\km}(\gamma)}\rangle\right) 
			\bigg]^{\mathrm{max}(\gamma)}.
		\end{split}
		\label{eq:7.5.34ss20}
	\end{equation}
	
	Indeed, by \eqref{eq:haha20imp}, we have the following 
	identity as elements in 
	$\Lambda^{\bullet}(\ku^{\perp}_{\kb}(\gamma)^{*})$, 
	\begin{equation}
		\begin{split}
			&\bigg\{\widehat{A}^{-1}(i\mathrm{ad}(y)|_{\ku_{\km}(\gamma)})\big[\frac{\det (1-
			e^{-i\mathrm{ad}(y)}\mathrm{Ad}(k))|_{\ku^{\perp}_{\km}(\gamma)}}{\det(1-\mathrm{Ad}(k))
			_{\ku^{\perp}_{\km}(\gamma)}}\big]^{\frac{1}{2}}\\
			&\mathrm{Tr_{s}}^{V_{\Lambda,\omega}}[\rho^{V_{\Lambda,\omega}}(e^{-iy}k)]\bigg\}|_{y=-\frac{\Omega^{\ku_{\km}(\gamma)}}{2\pi}}\\
			&=\sum_{\sigma\in W^{1}(\gamma)} 
			\varphi_{\gamma}(\sigma,\eta_{\omega}(\Lambda))\int_{f\in\mathcal{O}^{\gamma}
			_{\sigma(\eta_{\omega}(\Lambda)+\rho_{\ku_{\km}})}}e^{- \langle f,\Omega^{\ku_{\km}(\gamma)}\rangle}
			d\mu^{\gamma}_{\sigma}.
		\end{split}
		\label{eq:7.3.48ss20}
	\end{equation}
	
	Recall that the curvature form $\Omega^{\ku_{\kb}(\gamma)}$ is invariant by 
	the action of $U_{M}(\gamma)^{0}$ on $Y_{\kb}(\gamma)$. Since $a_{1}$ and 
	$\omega^{Y_{\kb}(\gamma)}$ are invariant by $U_{M}(\gamma)^{0}$-action, so is 
	$\Omega^{\ku_{\km}(\gamma)}$. Therefore, for $f\in\ku_{\km}(\gamma)^{*}$, $u\in 
	U_{M}(\gamma)^{0}$, then
	\begin{equation}
		\begin{split}
			&\big[\exp(-\frac{\beta(a_{1})^{2} 
			\omega^{Y_{\kb}(\gamma),2}}{8\pi^{2}t})\exp(-\langle \mathrm{Ad}^{*}(u)f, 
			\Omega^{\ku_{\km}(\gamma)}\rangle)\big]^{\mathrm{max}(\gamma)}\\
			&=\det 
			\mathrm{Ad}(u)|_{\ku^{\perp}_{\kb}(\gamma)} \big[\exp(-\frac{\beta(a_{1})^{2} 
			\omega^{Y_{\kb}(\gamma),2}}{8\pi^{2}t})\exp(-\langle f, 
			\Omega^{\ku_{\km}(\gamma)}\rangle)\big]^{\mathrm{max}(\gamma)}.
		\end{split}
		\label{eq:7.3.49ss20}
	\end{equation}
	Since $U_{M}(\gamma)^{0}$ acts on $\ku^{\perp}_{\kb}(\gamma)$ isometrically with respect 
	to $-B|_{\ku^{\perp}_{\kb}(\gamma)}$, then 
	\begin{equation}
		\det \mathrm{Ad}(u)|_{\ku^{\perp}_{\kb}(\gamma)} =1.
		\label{eq:7.3.50ss20}
	\end{equation}
	Then \eqref{eq:7.5.34ss20} follows from \eqref{eq:5.4.40qq20} and \eqref{eq:7.3.48ss20} - \eqref{eq:7.3.50ss20}.

	The right-hand side of \eqref{eq:7.5.34ss20} is a 
	polynomial in $t^{-1}$. Recall that $\dim 
	\ku^{\perp}_{\kb}(\gamma)=4l(\gamma)$. Then for each $\sigma\in 
	W^{1}(\gamma)$, we can rewrite the term 
	$\big[\cdots\big]^{\mathrm{max}(\gamma)}$ in the right-hand side of \eqref{eq:7.5.34ss20} as follows,
	\begin{equation}
		\sum_{j=0}^{l(\gamma)}\frac{1}{t^{j}}\frac{(-1)^{j}\beta(a_{1})^{2j}}{j!(2l(\gamma)-2j)! (8\pi^{2})^{j}}\big[\omega^{Y_{\kb}(\gamma),2j}\langle 
		\omega(\Lambda+\rho_{\ku}), 
		\Omega^{\ku_{\km}(\gamma)}\rangle^{2l(\gamma)-2j}\big]^{\mathrm{max}(\gamma)}.
		\label{eq:6.4.45ss20}
	\end{equation}
	
	Finally, we put together \eqref{eq:6.4.46ss20}, \eqref{eq:7.5.22decjan19}, 
	\eqref{eq:7.5.33kk20}, \eqref{eq:6.4.39sssks}, 
	\eqref{eq:7.5.34ss20}, and \eqref{eq:6.4.45ss20}, we get 
	\eqref{eq:6.4.48ss20bath}. This completes the 
	proof of our theorem.
\end{proof}

The Mellin transform of 
$\cE_{X,\gamma}(F,t)$ (if applicable) is defined by the 
following formula as a function in $s\in\C$ with $\Re(s)\gg 0$,
\begin{equation}
	\mathcal{ME}_{X,\gamma}(F,s)=-\frac{1}{\Gamma(s)}\int_0^{+\infty} 
	\cE_{X,\gamma}(F,t)t^{s-1}dt.
	\label{eq:7.3.53conf}
\end{equation}
If $\mathcal{ME}_{X,\gamma}(F,s)$ admits a meromorphic extension on $\C$ 
which is holomorphic at $s=0$, we will set 
\begin{equation}
	\cP\cE_{X,\gamma}(F)=\frac{\partial}{\partial 
	s}|_{s=0}\mathcal{ME}_{X,\gamma}(F,s).
\end{equation}

\begin{theorem}\label{thm:7.3.6bath}
	Suppose that the dominant weight $\Lambda$ is such that for every 
	$\omega\in W_{u}$, $b_{\Lambda,\omega}\neq 0$. 
	Then for $s\in\C$ with $\Re(s) > l(\gamma)+1$ , 
	$\mathcal{ME}_{X,\gamma}(F,s)$ is 
	well-defined and holomorphic, which admits a meromorphic extension to 
	$s\in\C$.
	
	Moreover, we have the following identity,
	\begin{equation}
		\begin{split}
			\mathcal{ME}_{X,\gamma}(F,s)=&-\frac{c_{G}(\gamma)}{\sqrt{2\pi}}\sum_{j=0}^{l(\gamma)}\frac{\Gamma(s-j-\frac{1}{2})}{\Gamma(s)}\\
			&\qquad\qquad\cdot\big[\sum_{\genfrac{}{}{0pt}{2}{\omega\in W_{u}}{\sigma \in 
			W^{1}(\gamma)}}\varepsilon(\omega) 
			\varphi_{\gamma}(\sigma,\eta_{\omega}(\Lambda))Q^{\gamma}_{j,\omega,\sigma}(\Lambda)(2\pi^2b_{\Lambda,\omega}^2)^{j+\frac{1}{2}-s}\big].
		\end{split}
		\label{eq:7.3.53kk20bath}
	\end{equation}
	Then $\mathcal{ME}_{X,\gamma}(F,s)$ is holomorphic at 
	$s=0$. We have
	\begin{equation}
		\begin{split}
			\cP\cE_{X,\gamma}(F)=&-\frac{c_{G}(\gamma)}{\sqrt{2}}\sum_{j=0}^{l(\gamma)}\frac{(-4)^{j+1}(j+1)!}{(2j+2)!}\\
			&\qquad\qquad\cdot\big[\sum_{\genfrac{}{}{0pt}{2}{\omega\in W_{u}}{\sigma \in 
			W^{1}(\gamma)}}\varepsilon(\omega) 
			\varphi_{\gamma}(\sigma,\eta_{\omega}(\Lambda))Q^{\gamma}_{j,\omega,\sigma}(\Lambda)(2\pi^2b_{\Lambda,\omega}^2)^{j+\frac{1}{2}}\big].
		\end{split}
		\label{eq:7.3.49kk20bath}
	\end{equation}
\end{theorem}
\begin{proof}
	By Theorem \ref{thm:7.3.1bath}, the assumption on $\Lambda$
	implies that $\cE_{X,\gamma}(F,t)$ decays exponential as 
	$t\rightarrow +\infty$. By \eqref{eq:6.3.4pl} and 
	\eqref{eq:6.4.48ss20bath}, we get \eqref{eq:7.3.53kk20bath}. This proves the first part of this theorem.
	
	The equation \eqref{eq:7.3.49kk20bath} is a direct 
	consequence of \eqref{eq:7.3.53kk20bath} by taking its derivative at 
	$0$. This 
	completes the proof of our theorem.
\end{proof}

The formula in the right-hand side of \eqref{eq:7.3.49kk20bath} still looks 
complicated, we can rewrite it in a neat way as follows. Let's introduce the 
following functions.
\begin{definition}
	
	Let $a^{1}\in \kb^{*}$ be which takes value $-1$ at $a_{1}$.
	Note that $\gamma\in T$. For $\omega\in W_{u}$, $\sigma\in W^{1}(\gamma)$, if $\Lambda\in P_{++}(U)$, for $z\in\C$, set
	\begin{equation}
		\begin{split}
			P^{\gamma}_{\omega,\sigma,\Lambda}(z)
			=\dim E^{\gamma}_{\omega,\sigma}\cdot\bigg[\exp\big(\langle 
			\Omega^{\ku_{\kb}(\gamma)},\sigma(\eta_{\omega}(\Lambda)+\rho_{\ku_{\km}})+z 
			\sqrt{-1}a^{1}\rangle\big)\bigg]^{\mathrm{max}(\gamma)}.
		\end{split}
		\label{eq:7.5.35ss20s}
	\end{equation}
\end{definition}

Since $\theta$ fix $\Omega^{\ku_{\kb}(\gamma)}$, by the fact that $\det 
\theta|_{\ku^{\perp}_{\kb}(\gamma)}=1$, then $P^{\gamma}_{\omega,\sigma,\Lambda}(z)$ is an even 
polynomial in $z$. Moreover, by the dimension formula 
\eqref{eq:5.4.40qq20}, the coefficients of $z^{j}$, 
$j\in \mathbb{N}$ in 
$P^{\gamma}_{\omega,\sigma,\Lambda}(z)$ are polynomials in $\Lambda$. 
Such polynomials are related to the Plancherel measures in the
representation theory.

\begin{lemma}\label{lm:7.4.2kk}
	We have the following identity
	\begin{equation}
		\begin{split}
			&\sum_{j=0}^{l(\gamma)}\frac{(-4)^{j+1}(j+1)!}{\sqrt{2}(2j+2)!}Q^{\gamma}_{j,\omega,\sigma}(\Lambda)(2\pi^2(b_{\Lambda,\omega})^2)^{j+\frac{1}{2}}\\
			&=- 2\pi \int_{0}^{|b_{\Lambda,\omega}|} 
			P^{\gamma}_{\omega, \sigma,\Lambda}(t)dt.
		\end{split}
		\label{eq:7.4.2kk20s}
	\end{equation}
\end{lemma}
\begin{proof}
	We have
	\begin{equation}
		\begin{split}
			&\langle \eta_{\omega}(\Lambda)+\rho_{\ku_{\km}}+z 
			\sqrt{-1}a^{1},\Omega^{\ku_{\kb}(\gamma)}\rangle\\
			&=z\beta(a_{1})\omega^{Y_{\kb}(\gamma)}+\langle 
			\omega(\Lambda+\rho_{\ku}),\Omega^{\ku_{\km}(\gamma)}\rangle.
		\end{split}
		\label{eq:7.4.3kk20s}
	\end{equation}
	Since $P^{\gamma}_{\omega,\sigma,\Lambda}(z)$ is an even function in 
	$z$, then
	\begin{equation}
		\begin{split}
			P^{\gamma}_{\omega,\sigma,\Lambda}(z)=&\dim 
			E^{\gamma}_{\omega,\sigma} \cdot
			\frac{1}{(2l(\gamma))!}\bigg[\left(z\beta(a_{1})\omega^{Y_{\kb}(\gamma)}+\langle 
			\omega(\Lambda+\rho_{\ku}),\Omega^{\ku_{\km}(\gamma)}\rangle\right)^{2l(\gamma)}\bigg]^{\mathrm{max}(\gamma)}\\
			=&\dim 
			E^{\gamma}_{\omega,\sigma} \cdot
			\sum_{j=0}^{l(\gamma)}\frac{\beta(a_{1})^{2j} 
			z^{2j}}{(2l(\gamma)-2j)!(2j)!}\\
			&\qquad\qquad\qquad\qquad\cdot\bigg[\omega^{Y_{\kb}(\gamma),2j}\langle 
			\omega(\Lambda+\rho_{\ku}),\Omega^{\ku_{\km}(\gamma)}\rangle^{2l(\gamma)-2j}\bigg]^{\mathrm{max}(\gamma)}.
		\end{split}
		\label{eq:7.4.4kk20s}
	\end{equation}
	
	Note that for $j=0,1,\cdots, l(\gamma)$,
	\begin{equation}
		\int_{0}^{|b_{\Lambda,\omega}|} 
		t^{2j}dt=\frac{1}{2j+1} 
		|b_{\Lambda,\omega}|^{2j+1}.
		\label{eq:7.4.5kk20s}
	\end{equation}
	Then \eqref{eq:7.4.2kk20s} is a consequence of 
	\eqref{eq:6.4.46ss20}, \eqref{eq:7.4.4kk20s} and 
	\eqref{eq:7.4.5kk20s}.
\end{proof}

As a consequence, we get the following formula for 
$\cP\cE_{X,\gamma}(F)$.
\begin{theorem}\label{thm:7.3.4bath}
	Suppose that the dominant weight $\Lambda$ is such that for every 
	$\omega\in W_{u}$, $b_{\Lambda,\omega}\neq 0$. Then
	\begin{equation}
		\cP\cE_{X,\gamma}(F)=2\pi c_{G}(\gamma)\cdot
		\sum_{\genfrac{}{}{0pt}{2}{\omega\in W_{u}}{\sigma \in 
		W^{1}(\gamma)}}\varepsilon(\omega) 
		\varphi_{\gamma}(\sigma,\eta_{\omega}(\Lambda))\int_{0}^{|b_{\Lambda,\omega}|} P^{\gamma}_{\omega,\sigma,\Lambda}(t)dt.
		\label{eq:7.5.35ss20}
	\end{equation}
\end{theorem}

%%%%%%%%%%%%%%%%%%%%%%%%%%%%%%%%%%%%%%%%%%%%%%%%%%%%%%%%%%%%%
\subsection{A family of representations of $G$}\label{section4.3}

We recall a 
definition of nondegeneracy of $\lambda$ in \cite[Definition 1.13 \& Proposition 8.12]{BMZ2015toeplitz}.

\begin{definition}\label{def:nondegdef}
	A dominant weight $\Lambda\in P_{++}(U)$ is said to be 
	nondegenerate with respect to the Cartan involution $\theta$ if 
	\begin{equation}
		W(U,T_{U})\cdot\Lambda\cap \kt^*=\emptyset.
	\end{equation}
	It is equivalent to 
	\begin{equation}
		\mathrm{Ad}^{*}(U)\Lambda\cap \kk^*=\emptyset.
	\end{equation}
	Note that if such dominant weight exists, we must have $\delta(G)>0$.
\end{definition}

Let $(E,\rho^{E})$ be the irreducible unitary
representation of $U$ with highest weight $\Lambda\in P_{++}(U)$. By the unitary 
trick, it extends to an irreducible representation of $G$, which we 
still denote by $(E,\rho^{E})$. Then $\Lambda$ 
being 
nondegenerate is equivalent to say that 
$(E,\rho^{E})$ is not isomorphic to 
$(E,\rho^{E}\circ \theta)$ as $G$-representation 
(as in \cite{MR3128980}).

\begin{definition}\label{def:6.2.2jan}
	If $\lambda\in\kt_{U}^{*}$, for $\omega\in W(U,T_{U})$, put
	\begin{equation}
		\begin{split}
			a_{\lambda,\omega}=\langle \omega\cdot \lambda, \ii
			a_{1}\rangle\in\R.
		\end{split}
		\label{eq:4.2.5oct}
	\end{equation}
	Recall the real number $b_{\lambda,\omega}$ is already defined by 
	\eqref{eq:7.2.9vogel}, then 
	$b_{\lambda,\omega}=a_{\lambda,\omega}+a_{\rho_{\ku},\omega}$. In particular, we simply put $a_{\lambda}=a_{\lambda,1}$, 
	$b_{\lambda}=b_{\lambda,1}$.
	
\end{definition}

\begin{lemma}\label{lm:3.3.2asymb}
	If $\lambda\in P_{++}(U)$ is nondegenerate, then for $\omega\in 
	W(U,T_{U})$, $a_{\lambda,\omega}\neq 0$.
\end{lemma}

Now we fix two dominant weights $\lambda,\lambda_{0}\in P_{++}(U)$. Let $\{(E_{d},\rho^{E_{d}})\}|_{d\in\mathbb{N}}$ be 
the sequence of representations of $G$ given by the irreducible 
unitary representations of $U$ with the highest weights 
$d\lambda+\lambda_{0}$, 
$d\in\mathbb{N}$. 

Put $F_{d}=G\times_{K} E_{d}$. Let $\mathbf{D}^{X,F_{d},2}$ denote the 
associated de Rham-Hodge Laplacian. For $t>0$, let 
$\exp(-t\mathbf{D}^{X,F_{d},2}/2)$ denote the heat operator 
associated with $\mathbf{D}^{X,F_{d},2}/2$. By taking 
$\Lambda=d\lambda+\lambda_{0}$, we apply our results in previous 
subsection to the sequence $\cE_{X,\gamma}(F_d,t)$, $d\in \mathbb{N}$.

%%%%%%%%%%%%%%%%%%%%%%%%%%%%%%%%%%%%%%%%%%%%%%%%%%%%%%%%%%%%%%%%%%
\subsection{Asymptotics for identity orbital integrals}\label{section7.3}
In this subsection, we specialize our results in Subsection 
\ref{sec7.2bath} for $\gamma=1$ and $\Lambda=d\lambda+\lambda_{0}$. 
Now the set $W^{1}(\gamma)$ reduces to $\{1\}$, and $l(\gamma)=l$, 
$\varphi_{\gamma}(\sigma,\eta_{\omega}(\Lambda))=1$. We will drop the 
superscript $\gamma$ and subscript $\sigma$ in our notation.

Moreover, for $\omega\in 
W_{u}$, the representation 
$E^{\gamma=1}_{\omega,\sigma=1}$ is just $V_{\Lambda,\omega}$ 
introduced in \eqref{eq:5.3.11pap}, which is the irreducible unitary 
representation of $U_{M}$ with highest weight 
$\eta_{\omega}(\Lambda)$ given by \eqref{eq:5.3.11final}.

\begin{definition}
	By taking 
	$\Lambda=d\lambda+\lambda_{0}$ in \eqref{eq:6.4.46ss20}, we define the following functions 
	in $d$, for $j=0,1,\cdots, l$, $\omega\in 
	W_{u}$, set
	\begin{equation}
		\begin{split}
			Q^{\lambda,\lambda_{0}}_{j,\omega}(d)&=Q_{j,\omega}(d\lambda+\lambda_{0})\\
			&=\frac{(-1)^{j}\beta(a_{1})^{2j}}{j!(2l-2j)! 
			(8\pi^{2})^{j}}\dim V_{d\lambda+\lambda_{0},\omega} \\
			&\qquad\cdot\big[\omega^{Y_{\kb},2j}\langle 
			\omega(d\lambda+\lambda_{0}+\rho_{\ku}), 
			\Omega^{\ku_{\km}}\rangle^{2l-2j}\big]^{\mathrm{max}}.
		\end{split}
		\label{eq:qcoeffbath}
	\end{equation}
	By the Weyl dimension formula, $\dim 
	V_{d\lambda+\lambda_{0},\omega}$ is a polynomial in $d$. Then 
	$Q^{\lambda,\lambda_{0}}_{j,\omega}(d)$ is a polynomial in $d$ of degree $\leq \frac{\dim(\g/\kh)}{2} - 2j$.
\end{definition}

By Theorem \ref{thm:7.3.1bath} and \eqref{eq:qcoeffbath}, we get 
directly the following results.
\begin{theorem}\label{thm:7.3.2ss20}
	For $t>0$, we have the following identity
	\begin{equation}
		\cI_{X}(F_{d},t)=\frac{c_{G}}{\sqrt{2\pi t}}\sum_{j=0}^{l} 
		t^{-j}\sum_{\omega\in 
		W_{u}} \varepsilon(\omega) 
		e^{-2\pi^{2}t(da_{\lambda,\omega}+b_{\lambda_{0},\omega})^{2}}Q^{\lambda,\lambda_{0}}_{j,\omega}(d).
		\label{eq:6.4.48ss20}
	\end{equation}
\end{theorem}

\begin{theorem}\label{thm:7.3.3ss}
	Suppose that $\lambda$ is nondegenerate with respect to $\theta$. 
	For $d\in\bN$ large enough and for $s\in\C$ with $\Re(s)\gg 0$ , 
	$\mathcal{MI}_{X}(F_{d},s)$ is 
	well-defined and holomorphic, which admits a unique meromorphic extension to 
	$s\in\C$ and is holomorphic at $s=0$.
	
	Moreover, we have the following identities,
	\begin{equation}
		\begin{split}
			\mathcal{MI}_{X}(F_{d},s)=&-\frac{c_{G}}{\sqrt{2\pi}}\sum_{j=0}^{l}\frac{\Gamma(s-j-\frac{1}{2})}{\Gamma(s)}\\
			&\qquad\qquad\cdot\big[\sum_{\omega\in 
			W_u}\varepsilon(\omega)Q^{\lambda,\lambda_{0}}_{j,\omega}(d)(2\pi^2(da_{\lambda,\omega}+b_{\lambda_{0},\omega})^2)^{j+\frac{1}{2}-s}\big],
		\end{split}
		\label{eq:7.3.53kk20}
	\end{equation}
	and
	\begin{equation}
		\begin{split}
			\cP\cI_{X}(F_{d})=&-\frac{c_{G}}{\sqrt{2}}\sum_{j=0}^{l}\frac{(-4)^{j+1}(j+1)!}{(2j+2)!}\\
			&\cdot\big[\sum_{\omega\in W_u}\varepsilon(\omega)Q^{\lambda,\lambda_{0}}_{j,\omega}(d)(2\pi^2(da_{\lambda,\omega}+b_{\lambda_{0},\omega})^2)^{j+\frac{1}{2}}\big].
		\end{split}
		\label{eq:7.3.49kk20}
	\end{equation}
	In particular, the quantity $\cP\cI_{X}(F_{d})$ is a polynomial in $d$ for 
	$d$ large enough, whose coefficients 
	depend only on the given root system and $\lambda$, $\lambda_{0}$, and has 
	degree $\leq \frac{\dim(\g/\kh)}{2}+1$.
\end{theorem}
\begin{proof}
	Since $\lambda$ is nondegenerate, by Lemma 
	\ref{lm:3.3.2asymb}, $a_{\lambda,\omega}\neq 0$, 
	$\omega\in W_{u}$. Then there exists $d_{0}\in 
	\bN$ such that for $d\geq d_{0}$, 
	$(da_{\lambda,\omega}+b_{\lambda_{0},\omega})^{2}>0$. Then by Theorem 
	\ref{thm:7.3.6bath}, we get 
	first part of this theorem and \eqref{eq:7.3.53kk20}, 
	\eqref{eq:7.3.49kk20}.
	
	Note that 
	$[(da_{\lambda,\omega}+b_{\lambda_{0},\omega})^2]^{1/2}=|da_{\lambda,\omega}+b_{\lambda_{0},\omega}|$. For $d\gg d_{0}$, $|da_{\lambda,\omega}+b_{\lambda_{0},\omega}|=\mathrm{sign}(a_{\lambda,\omega})(da_{\lambda,\omega}+b_{\lambda_{0},\omega})$. Then we see that $\cP\cI_{X}(F_{d})$ is a polynomial in $d$ for $d$ large enough. This 
	completes the proof of our theorem.
\end{proof}

As explained in Remark \ref{rm:important}, when $G$ has noncompact 
center with $\delta(G)=1$ (but $U$ is still assumed to be compact), 
most of the above computations can be reduce into very simple ones. 
Recall that $a_{\lambda}, b_{\lambda_{0}}\in\R$ are defined in 
Definition \ref{def:6.2.2jan}. 
\begin{corollary}\label{cor:7.3.6kkss}
	Assume that $U$ is compact and that $G$ has noncompact center 
	with $\delta(G)=1$, and assume that $\lambda$ is nondegenerate. Then for $t>0$, $s\in\C$,
	\begin{equation}
		\begin{split}
			&\cI_{X}(F_{d},t)=\frac{c_{G}}{\sqrt{2\pi 
			t}}e^{-2\pi^{2}t(da_{\lambda}+b_{\lambda_{0}})^{2}}\dim 
			E_{d},\\
			& 
			\mathcal{MI}_{X}(F_{d},s)=-\frac{c_{G}}{\sqrt{2\pi}}\frac{\Gamma(s-\frac{1}{2})}{\Gamma(s)}\big(2\pi^{2}(da_{\lambda}+b_{\lambda_{0}})^{2}\big)^{1/2-s}\dim 
			E_{d}.
		\end{split}
		\label{eq:7.3.57std}
	\end{equation}
	Furthermore,
	\begin{equation}
		\cP\cI_{X}(F_{d})=2\pi 
		c_{G}|da_{\lambda}+b_{\lambda_{0}}|\dim E_{d}.
		\label{eq:7.3.58std}
	\end{equation}
\end{corollary}
\begin{proof}
	By the hypothesis, we get that $l=0$, $W_{u}=\{1\}$ 
	and $Q^{\lambda,\lambda_{0}}_{0,1}(d)=\dim E_{d}$. Then 
	\eqref{eq:7.3.57std}, \eqref{eq:7.3.58std} are just special cases 
	of \eqref{eq:6.4.48ss20}, \eqref{eq:7.3.53kk20} and 
	\eqref{eq:7.3.49kk20}.
	
	However, we can prove them more directly using a result of Proposition \ref{prop:3.6.5nov19}. It 
	is enough to prove the first identity in \eqref{eq:7.3.57std}. 
	Note that by \eqref{eq:5.3.14conf}, we have
	\begin{equation}
		X'=M/K,
		\label{eq:7.3.59std}
	\end{equation}
	with $\delta(X')=0$.
	
	By \cite[Proposition 5.2]{MR3128980} or 
	\cite[Proposition 4.1]{Shen_2016}, we have
	\begin{equation}
		\begin{split}
			&[e(TX',\nabla^{TX'})]^{\mathrm{max}}\\
			&=(-1)^{\frac{m-1}{2}}\frac{|W(U_{M},T)|/|W(K,T)|}{\mathrm{Vol}(U_{M}/K)}.
		\end{split}
		\label{eq:7.4.19pp20}
	\end{equation}
	Then by \eqref{eq:eq:6.4.47ss20}, we have
	\begin{equation}
		[e(TX',\nabla^{TX'})]^{\mathrm{max}}=-c_{G}.
		\label{eq:7.3.61std}
	\end{equation}
	By \eqref{eq:4.1.28std} and \eqref{eq:4.2.5oct}, we have
	\begin{equation}
		\alpha_{E_{d}}=-2\pi(da_{\lambda}+b_{\lambda_{0}}).
		\label{eq:7.3.62std}
	\end{equation}
	
	Combing \eqref{eq:3.6.21novend} and \eqref{eq:7.4.19pp20} - 
	\eqref{eq:7.3.62std}, we get the first identity in 
	\eqref{eq:7.3.57std}, and hence the other identities. This gives 
	a second proof to this corollary. 
\end{proof}

%%%%%%%%%%%%%%%%%%%%%%%%%%%%%%%%%%%%%%%%%%%%%%%%%%%%%%%%%%%%%%%%%%%%%%%%%%%
\subsection{Connection to M\"{u}ller-Pfaff's 
results}\label{section7.4kk}
In this subsection, we assume that $G$ has compact center with 
$\delta(G)=1$. We explain here how to connect our computations in previous 
subsection to M\"{u}ller-Pfaff's results in \cite{MR3128980}.

For $\gamma=1$, $\omega\in W_{u}$, the function 
$P^{\gamma}_{\omega,\sigma, \Lambda}$ defined in \eqref{eq:7.5.35ss20s} now reduces to 
\begin{equation}
	P_{\omega,\Lambda}(z)=\dim V_{\Lambda,\omega} \bigg[\exp\big(\langle 
	\eta_{\omega}(\Lambda)+\rho_{\ku_{\km}}+z 
	\sqrt{-1}a^{1}, \Omega^{\ku(\kb)}\rangle\big)\bigg]^{\mathrm{max}}.
	\label{eq:6.4.52ss20s}
\end{equation}

We can verify directly that
\begin{equation}
	P_{\omega,\Lambda}(z)=\frac{\mathrm{Vol}(U_{M}/T)}{\mathrm{Vol}(U/T_{U})}\Pi_{\alpha^{0}\in R^{+}(\ku,\kt_{U})}\frac{\langle \alpha^{0},\eta_{\omega}(\Lambda)+\rho_{\ku_{\km}}+z\ii a^{1}\rangle}{\langle \alpha^{0},\rho_{\ku}\rangle}.
	\label{eq:7.4.2kk20}
\end{equation}
The scalar product in \eqref{eq:7.4.2kk20} is taken with 
respect to $-B|_{\ku}$. Up to a universal constant, $P_{\omega,\Lambda}(z)$ is just the polynomial related to the 
Plancherel measure of representation $V_{\Lambda,\omega}$ as 
given in \cite[Eq. (6.10)]{MR3128980}. Note that there is no factor 
$(2\pi)^{2l}$ appeared in \eqref{eq:7.4.2kk20} 
because of our normalization for $[\cdot]^{\mathrm{max}}$.

By Theorem \ref{thm:7.3.4bath}, 
we have the following result for sufficiently large $d$.
\begin{corollary}\label{cor:7.4.3ss20}
	Suppose that $\lambda$ is nondegenerate with respect to $\theta$. 
	Then
	\begin{equation}
		\cP\cI_{X}(F_{d})=2\pi c_{G} \sum_{\omega\in W_{u}}\varepsilon(\omega) 
		\int_{0}^{|da_{\lambda,\omega}+b_{\lambda_{0},\omega}|} P_{\omega, d\lambda+\lambda_{0}}(t)dt.
		\label{eq:7.4.3ss20}
	\end{equation}
\end{corollary}

By \cite[Lemma 6.1]{MR3128980}, we can get the following identity, 
\begin{equation}
	|W(K,T)|=2|W(K_{M},T)|.
	\label{eq:6.4.54ss20}
\end{equation}
Combining \eqref{eq:eq:6.4.47ss20}, \eqref{eq:7.4.2kk20}, 
\eqref{eq:6.4.54ss20}, we see that the formula in Corollary \ref{cor:7.4.3ss20},
is exactly the same formula of M\"{u}ller-Pfaff 
\cite[Proposition 6.6]{MR3128980} for $\cP\cI_{X}(F_{d})$.

Recall that the $U$-representation $E_{d}$ has highest weight 
$d\lambda+\lambda_{0}\in P_{++}(U)$. Then by Weyl dimension formula,  $\dim E_{d}$ is a polynomial in $d$. If $\lambda$ 
is regular, then the degree (in $d$) of $\dim E_{d}$ is 
$\frac{\dim \g/\kh}{2}$.

For determining the leading term of 
$\cP\cI_{X}(F_{d})$, as mentioned in the introduction 
part, we can specialize the result of Bismut-Ma-Zhang 
\cite[Theorem 0.1]{BMZ2015toeplitz} as in \cite[Section 
8]{BMZ2015toeplitz} for the symmetric space $X$. 
Here to emphasize $\cP\cI_{X}(F_{d})$ being a polynomial in $d$, we state a result of M\"{u}ller-Pfaff 
\cite[Proposition 1.3]{MR3128980} as follows.
\begin{proposition}\label{prop:7.4.4}
	Suppose that $\lambda$ is nondegenerate and that $\lambda_{0}=0$. Then there 
	exists a constant $C_{X,\lambda}\neq 0$ such that
	\begin{equation}
		\cP\cI_{X}(F_{d})=C_{X,\lambda} d \dim 
		E_{d} + R(d),
		\label{eq:7.4.9kk20}
	\end{equation}
	where $R(d)$ is a polynomial whose degree is no greater 
	than the degree of $\dim E_{d}$.
\end{proposition}
\begin{remark}
	Note that in \cite[Proposition 1.3]{MR3128980}, M\"{u}ller and 
	Pfaff proved Proposition 
	\ref{prop:7.4.4} by reducing the problems to the cases $G=\mathrm{SL}_{3}(\R)$ and 
	$\mathrm{SO}^{0}(p,q)$($pq>1$ odd). In 
	particular, for certain examples of $\lambda$, they also worked out explicitly the constant 
	$C_{X,\lambda}$ \cite[Corollaries 1.4 \& 1.5]{MR3128980}. 
	
	Similarly, if we take a nonzero 
	$\lambda_{0}$, we can repeat their computations for 
	$G=\mathrm{SL}_{3}(\R)$ and $\mathrm{SO}^{0}(p,q)$ ($pq>1$ 
	odd) in order to get more explicit information on the leading 
	terms of $\cP\cI_{X}(F_{d})$.
\end{remark}

An important step in 
M\"{u}ller-Pfaff's proof to Proposition \ref{prop:7.4.4} is reducing 
the computation of $\cP\cI_{X}(F_{d})$ to 
the cases where $\g=\mathfrak{sl}_{3}(\R)$ or 
$\mathfrak{so}(p,q)$ with $pq>1$ odd. Such reduction is already 
explained in Subsection \ref{subs:4.2std}. More precisely, we have
\begin{equation}
	X=X_{1}\times X_{2},
	\label{eq:}
\end{equation}
where $X_{1}$ is one case listed in \eqref{eq:4.0.33kk}, and $X_{2}$ 
is a symmetric space rank $0$. 

We use the notation in Subsection \ref{subs:4.2std} and assume $G$ to 
be semisimple. Let 
$\lambda_{i}$, $\lambda_{0,i}$ be dominant weights of $U_{i}$, 
$i=1,2$ such that 
\begin{equation}
	\lambda=\lambda_{1}+\lambda_2, \; 
	\lambda_{0}=\lambda_{0,1}+\lambda_{0,2}.
	\label{eq:7.4.16kk20}
\end{equation}
Now we consider the sequence $d\lambda+\lambda_{0}$, $d\in\bN$. Then
\begin{equation}
	E_{d}=E_{d\lambda_{1}+\lambda_{0,1}}\otimes 
	E_{d\lambda_{2}+\lambda_{0,2}}.
	\label{eq:7.4.17kk20}
\end{equation}

Since $G_{2}$ is equal rank, the nondegeneracy of $\lambda$  with respect 
to $\theta$ is equivalent to the nondegeneracy of $\lambda_{1}$ with 
respect to $\theta_{1}$. Then by Proposition \ref{prop:4.2.2bio}, 
after taking the Mellin transform, we have
\begin{equation}
	\mathcal{M}\cI_{X}(F_{d},s)=[e(TX_{2},\nabla^{TX_{2}})]^{\mathrm{max}_{2}}\dim E_{d\lambda_2+\lambda_{0,2}}\mathcal{M}\cI_{X_{1}}(F_{d\lambda_1+\lambda_{0,1}},s).
	\label{eq:7.4.19kk20}
\end{equation}
Then
\begin{equation}
	\cP\cI_{X}(F_{d})=[e(TX_{2},\nabla^{TX_{2}})]^{\mathrm{max}_{2}}\dim 
	E_{d\lambda_2+\lambda_{0,2}}\cP\cI_{X_{1}}(F_{d\lambda_1+\lambda_{0,1}}).
	\label{eq:7.4.20kk20}
\end{equation}
Then we only need to evaluate 
$\cP\cI_{X_{1}}(F_{d\lambda_1+\lambda_{0,1}})$ explicitly, which 
has been dealt with in \cite[Section 6]{MR3128980}.

%%%%%%%%%%%%%%%%%%%%%%%%%%%%%%%%%%%%%%%%%%%%%%%%%%%%%%%%%%%%%%%%%%%%%
\subsection{Asymptotic elliptic orbital integrals}\label{subsection7.5}
\begin{definition}\label{def:7.5.1kk20}
	A function $f(d)$ in $d$ is called an exponential polynomial in $d$ if 
	it is a finite sum of the term $c_{j,s}e^{2\pi\sqrt{-1}sd}d^{j}$ 
	with $j\in\mathbb{N}$, $s\in\R$, $c_{j,s}\in \C$. The largest 
	$j\geq 0$ 
	such that $c_{j,s}\neq 0$ in $f(d)$ is called the degree of 
	$f(d)$.
	
	We say that the oscillating term $e^{2\pi\sqrt{-1}sd}$ is nice if 
	$s\in\mathbb{Q}$. We say that an exponential polynomial $f(d)$ in $d$ 
	is nice 
	if all its oscillating terms are nice. 
\end{definition}

\begin{remark}
	If $f(d)$ is a nice exponential polynomial in $d$, then there exists a 
	$N_{0}\in\mathbb{N}_{>0}$ such that the function $f(dN_{0})$ is a 
	polynomial in $d$.
\end{remark}

Note that by \eqref{eq:5.4.18kk20}, 
$\varphi_{\gamma}(\sigma,\eta_{\omega}(d\lambda+\lambda_{0}))$ is an 
oscillating term in $d$, which is nice when $\gamma\in T$ is of finite 
order. The following theorem is a direct 
consequence of Theorem \ref{thm:7.3.4bath}.

\begin{theorem}\label{thm:7.5.5kk20}
	Suppose that $\lambda$ is nondegenerate, 
	and that $\gamma=k\in T$. Then, for sufficiently large $d$, $\cP\cE_{X,\gamma}(F_{d})$ is an 
	exponential polynomial in $d$. Moreover, we have
	\begin{equation}
		\cP\cE_{X,\gamma}(F_{d})=2\pi c_{G}(\gamma)\cdot
		\sum_{\genfrac{}{}{0pt}{2}{\omega\in W_{u}}{\sigma \in 
		W^{1}(\gamma)}}\varepsilon(\omega) 
		\varphi_{\gamma}(\sigma,\eta_{\omega}(d\lambda+\lambda_{0}))\int_{0}^{|da_{\lambda,\omega}+b_{\lambda_{0},\omega}|} P^{\gamma}_{\omega,\sigma,d\lambda+\lambda_{0}}(t)dt.
		\label{eq:7.5.35ss20spa}
	\end{equation}
\end{theorem}

If we consider $G=\mathrm{Spin}(1, 2n+1)$, $n\geq 1$ as 
in \cite{Fedosova2015compact}, then up to a constant, 
the exponential polynomial
$\sum_{\sigma\in 
W^{1}(\gamma)}\varphi_{\gamma}(\sigma,\eta_{\omega}(d\lambda+\lambda_{0}))P^{\gamma}_{\omega,\sigma,d\lambda+\lambda_{0}}(t)$ is just the one defined by Fedosova in \cite[Proposition 5.1]{Fedosova2015compact}. This way, our results are compatible with her results in \cite[Theorem 1,1]{Fedosova2015compact} for hyperbolic orbifolds.

\begin{remark}
	Let $\mathrm{Char}(A)$ denote the character 
	ring of the complex representations of a compact Lie group $A$.
	One key ingredient in  \eqref{eq:7.5.35ss20} is an explicit decomposition of 
	characters of $U$ into characters of $U_{M}(\gamma)^{0}$. In the 
	diagram \eqref{eq:7.5.6kk20}, we give two different ways of this 
	decomposition. The 
	formula in \eqref{eq:7.5.35ss20} is obtained by the computations 
	along the lower path in \eqref{eq:7.5.6kk20}. We also have the 
	upper path, which is essentially the geometric localization 
	formula obtained in Theorem \ref{thm:6.2.1ss}.
	\begin{equation}
		\xymatrixcolsep{4pc}\xymatrix{
		& \mathrm{Char}(U(\gamma)^{0})\ar[dr]^{\otimes 
		\Lambda^{\bullet}\kn(\gamma)^{*}_{\C}} &
		\\
		\mathrm{Char}(U)\ar[dr]^{\otimes 
		\Lambda^{\bullet}\kn^{*}_{\C}} 
		\ar[ur]^{\mathrm{Kirillov\;for\;}\gamma\;\in\;U}& 
		&\mathrm{Char}(U_{M}(\gamma)^{0})\\
		& 
		\mathrm{Char}(U_{M})\ar[ur]^{\mathrm{Kirillov\;for\;}\gamma\;\in\; U_{M}} & }
		\label{eq:7.5.6kk20}
	\end{equation}
\end{remark}

We will use the same notation as in Section \ref{section6paris}. The 
following theorem is a consequence of the geometric localization 
formula obtained in Theorem \ref{thm:6.2.1ss}.

For $k\in T$, let 
$W^{1}_{U}(k)\subset W(U,T_{U})$ be defined as in 
\eqref{eq:5.4.13ss20} with respect to $R^{+}(\ku,\kt_{U})$. For $\sigma\in W^{1}_{U}(k)$, the term 
$\varphi^{U}_{k}(\sigma,d\lambda+\lambda_{0})$ defined 
as in \eqref{eq:6.2.7kk20} is an 
oscillating term, which is nice if $k$ is of finite 
order. 

\begin{theorem}\label{thm:7.5.3conf}
	Suppose that $\gamma=k\in T$ is elliptic and that $\lambda$ is 
	nondegenerate with respect to $\theta$. Then 
	for $\sigma\in W^{1}_{U}(k)$, $\sigma\lambda\in 
	P_{++}(\widetilde{U}(k))$ is nondegenerate with respect to 
	the Cartan involution $\theta$ on $\z(k)$. For $d\in\bN$, 
	let $E_{\sigma,d}^{k}$ be the irreducible unitary 
	representation of $\widetilde{U}(k)$ with highest weight 
	$d\sigma\lambda+\sigma(\lambda_{0}+\rho_{\ku})-\rho_{\ku(k)}$. 
	This way, we get a sequence of flat vector bundles 
	$\{F_{\sigma,d}^{k}\}_{d\in\bN}$ on $X(k)$. Then for sufficiently 
	large $d$, we have
	\begin{equation}
		\cP\cE_{X,\gamma}(F_{d})=\sum_{\sigma\in 
		W^{1}_{U}(k)}\varphi^{U}_{k}(\sigma,d\lambda+\lambda_{0}) 
		\cP\cI_{X(k)}(F_{\sigma,d}^{k}).
		\label{eq:7.5.2kk20}
	\end{equation}
\end{theorem}
\begin{proof}
	The nondegeneracy of $\sigma\lambda$ ($\sigma\in W^{1}_{U}(k)$) 
	follows easily from the nondegeneracy of $\lambda$ and the 
	definition of $W^{1}_{U}(k)$. For proving this theorem, we only need to prove 
	\eqref{eq:7.5.2kk20}. Actually, by Theorem \ref{thm:6.2.1ss}, for 
	$t>0$, we get
	\begin{equation}
		\cE_{X,\gamma}(F_{d},t)=\sum_{\sigma\in 
		W^{1}_{U}(k)}\varphi^{U}_{k}(\sigma,d\lambda+\lambda_{0}) 
		\cI_{X(k)}(F_{\sigma,d}^{k},t),
		\label{eq:7.5.3std}
	\end{equation}
	Then \eqref{eq:7.5.2kk20} follows from the linearity of Mellin 
	transform. This completes the proof of our theorem.
\end{proof}

%%%%%%%%%%%%%%%%%%%%%%%%%%%%%%%%%%%%%%%%%%
\section{A proof of Theorem \ref{thm:maintheorem}}\label{section8paris}
In this section, we complete the proof of Theorem 
\ref{thm:maintheorem}, then Theorem \ref{thm:0.00001s} (and Theorem 
\ref{thm:1.0.1bis}) follows as a 
consequence. We assume that $G$ is a connected 
linear real reductive Lie group with $\delta(G)=1$ and compact 
center, so that $U$ is a compact Lie group.

%%%%%%%%%%%%%%%%%%%%%%%%%%%%%%%%%%%%%%%%%%%%%%%%%%%%%%%%%%%%%%%%%%%%%%%%%%%%%%%
\subsection{A lower bound for the Hodge Laplacian on $X$}\label{section:lower}
We use the notation as in Subsection \ref{section3.6}. Recall that $e_1$, 
$\cdots$, $e_m$ is an orthogonal basis of $TX$ or $\pp$. 
Put
\begin{equation}
	C^{\g,H}=-\sum_{j=1}^{m} e_{j}^{2}\in U\g.
	\label{eq:8.1.1ss20}
\end{equation}
Let $C^{\g,H,E}$ be its action on $E$ via $\rho^{E}$. Then
\begin{equation}
	C^{\g,E}=C^{\g,H,E}+C^{\kk,E}.
	\label{eq:7.7.16pp}
\end{equation}

Let $\Delta^{H,X}$ be the Bochner-Laplace operator on bundle 
$\Lambda^{\bullet}(T^*X)\otimes F$ associated with the unitary 
connection $\nabla^{\Lambda^{\bullet}(T^*X)\otimes F,u}$. Put
\begin{equation}
	\begin{split}
		\Theta(F)=\frac{S^X}{4}-\frac{1}{8} \langle R^{TX}(e_i,e_j)e_k,e_\ell\rangle c(e_i)c(e_j)\widehat{c}(e_k)\widehat{c}(e_\ell)&\\
		-C^{\g,H,E}+\frac{1}{2}\big(c(e_i)c(e_j)-\widehat{c}(e_i)\widehat{c}(e_j)\big)R^{F}(e_i, e_j),
	\end{split}
	\label{eq:8.5.3lara}
\end{equation}
where $R^{F}$ is the curvature of the unitary connection $\nabla^{F}$ 
on $F$.

Then $\Theta(F)$ is a self-adjoint section of 
$\mathrm{End}(\Lambda^{\bullet}(T^*X)\otimes F)$, which is parallel with 
respect to $\nabla^{\Lambda^{\bullet}(T^*X)\otimes F,u}$. 
Equivalently, 
$\Theta(F)$ is an element in $\mathrm{End}(\Lambda^{\bullet}(\pp^{*})\otimes E)$ 
which commutes with $K$-action.
By 
\cite[Eq. (8.39)]{BMZ2015toeplitz}, we have
\begin{equation}
	\begin{split}
		\mathbf{D}^{X,F,2}=-\Delta^{H,X}+\Theta(F).
	\end{split}
	\label{eq:8.5.5lara}
\end{equation}
Then for $s\in \Omega_c^\cdot(X, F)$, we have
\begin{equation}
	\langle \mathbf{D}^{X,F,2}s,s\rangle_{L_2}\geq \langle \Theta(F)s,s\rangle_{L_2}.
	\label{eq:8.5.6lara}
\end{equation}

Let $\Delta^{H,X,i}$ denote the Bochner-Laplace operator acting on 
$\Omega^i(X,F)$, and let $p^{H,i}_t(x,x')$ be the kernel of 
$\exp(t\Delta^{H,X,i}/2)$ on $X$ with respect to $dx'$. We will 
denote by $p^{H,i}_t(g)\in\mathrm{End}(\Lambda^i(\pp^*)\otimes E)$ 
its lift to $G$ explained in Subsection \ref{section3.2}.
Let $\Delta^{X}_0$ be the scalar Laplacian on $X$ with the heat kernel $p^{X,0}_t$.

Let $||p^{H,i}_t(g)||$ be the operator norm of $p^{H,i}_t(g)$ in 
$\mathrm{End}(\Lambda^i(\pp^*)\otimes E)$. By \cite[Proposition 
3.1]{MP2013raysinger}, if $g\in G$, then
\begin{equation}
	||p^{H,i}_t(g)||\leq p^{X,0}_t(g).
	\label{eq:8.6.5qq}
\end{equation}

Let $p^H_t$ be the kernel of $\exp(t\Delta^{H,X}/2)$, then
\begin{equation}
	p^H_t=\oplus_{i=1}^p p^{H,i}_t.
	\label{eq:8.5.8sud}
\end{equation}
Let $q^{X,F}_t$ be the heat kernel associated with 
$\mathbf{D}^{X,F,2}/2$, by \eqref{eq:8.5.5lara}, for $g\in X$, 
\begin{equation}
	q^{X,F}_t(g)=\exp(-t\Theta(F)/2)p^{H}_t(g).
	\label{eq:8.5.9}
\end{equation}

Recall that $P_{++}(U)$ is the set of dominant weights of $U$ with 
respect to $R^{+}(\ku,\kt_{U})$ defined in Subsection 
\ref{section5.3pl}. As in Subsection \ref{section4.3}, we fix 
$\lambda, \lambda_{0}\in P_{++}(U)$ such that $\lambda$ is nondegenerate with respect to $\theta$.
Recall that for $\rd\in\bN$, $(E_\rd,\rho^{E_{d}})$ is the irreducible 
unitary representation of $U$ with highest weight $d\lambda+\lambda_{0}$, 
which extends uniquely to a representation of $G$. By \cite[Th\'{e}or\`{e}me 3.2]{MR2838248} \cite[Theorem 4.4 and Remark 
4.5]{BMZ2015toeplitz} and \cite[Proposition 7.5]{MR3128980}, there exist $c>0$, $C>0$ such that, for $\rd\in\mathbb{N}$,
\begin{equation}
	\Theta(F_\rd)\geq c\rd^2- C,
	\label{eq:8.5.4sudpps}
\end{equation}
where the estimate $d^{2}$ comes from the positive operator 
$C^{\g,H,E_{d}}$.
By \eqref{eq:8.5.5lara}, \eqref{eq:8.5.6lara}, \eqref{eq:8.5.4sudpps}, we get
\begin{equation}
	\mathbf{D}^{X,F_\rd,2}\geq c\rd^2- C.
	\label{eq:8.4.26pps}
\end{equation}

\begin{lemma}\label{lm:8.5.1sud}
	There exists $\rd_0\in \bN$ and $c_0>0$ such that if $\rd\geq 
	\rd_0$, $g\in G$
	\begin{equation}
		||q^{X,F_\rd}_t(g)||\leq e^{-c_0\rd^2 t} p^{X,0}_t(g).
		\label{eq:8.5.13sud}
	\end{equation}
\end{lemma}

\begin{proof}
	By \eqref{eq:8.5.4sudpps}, there exist $\rd_0\in\bN$, $c'>0$ such that if $\rd\geq \rd_0$, 
	\begin{equation}
		\Theta(F_\rd)\geq c' \rd^2.
		\label{eq:8.5.14sud}
	\end{equation}
	Then if $t>0$,
	\begin{equation}
		||\exp(-t\Theta(F_\rd)/2)||\leq e^{-c'\rd^2 t/2}.
		\label{eq:8.5.15sud}
	\end{equation}
	By \eqref{eq:8.6.5qq}, \eqref{eq:8.5.8sud}, \eqref{eq:8.5.9}, \eqref{eq:8.5.15sud}, we get \eqref{eq:8.5.13sud}. 
\end{proof}

The locally symmetric orbifold $Z$ is defined as $\Gamma\backslash 
X$, where $\Gamma$ is a cocompact discrete subgroup of $G$. For 
$\gamma\in\Gamma$, the number $m_{\gamma}\geq 0$ is given by 
\eqref{eq:3.3.2ff}, which only depends on the conjugacy class of 
$\gamma$ (in $G$ or $\Gamma$). Recall that $E[\Gamma]$ is the finite 
set of elliptic conjugacy classes in $\Gamma$.

For $t>0$, $x\in X$, $\gamma\in \Gamma$, set
\begin{equation}
	v_t(F_\rd,\gamma, x) = \mathrm{Tr_s}^{\Lambda^{\bullet}(T^*X)\otimes 
	F_\rd}\bigg[\big(N^{\Lambda^{\bullet}(T^*X)}-\frac{m}{2}\big)q^{X,F_\rd}_{t}(x,\gamma(x))\gamma\bigg].
	\label{eq:8.6.7mmk}
\end{equation}
Then by Lemma \ref{lm:8.5.1sud}, we have the 
following result.
\begin{lemma}\label{lm:8.6.3ugc}
	There exist $C_0>0$, $c_0>0$ such that if $\rd$ is large enough, for $t >0$, $x\in X$, $\gamma\in \Gamma$,
	\begin{equation}
		|v_t(F_\rd, \gamma, x)|\leq C_0 (\dim E_\rd) e^{-c_0 d^2 t}p^{X,0}_{t}(x,\gamma(x)).
		\label{eq:8.6.8mk2}
	\end{equation}
\end{lemma}

Set
\begin{equation}
	m_{\Gamma}=\inf_{[\gamma]\in 
	[\Gamma] - E[\Gamma]} m_{\gamma}.
	\label{eq:ss8.6.1copy}
\end{equation}
By \cite[Proposition 1.8.5]{liu:tel-01841334}, $m_{\Gamma}>0$.

\begin{proposition}\label{prop:8.6.5ugc}
	There exist 
	constants $C>0$, $c>0$ such that if $x\in X$, $t\in\, ]0,1]$, then
	\begin{equation}
		\sum_{\gamma\in\Gamma, \gamma \mathrm{\;nonelliptic}} p^{X,0}_t(x,\gamma(x))\leq C\exp(-c/t).
		\label{eq:8.6.11ugc}
	\end{equation}
\end{proposition}

\begin{proof}
	By \cite[Theorem 3.3]{Donnelly1979asymptotic}, there exists $C_0>0$ such that when $0<t\leq 1$,
	\begin{equation}
		p_t^{X,0}(x,x')\leq C_0 t^{-m/2} \exp(-\frac{d^2(x,x')}{4t}).
		\label{eq:8.6.6}
	\end{equation}
	
	By \cite[Lemma 1.8.6]{liu:tel-01841334}, there exist $c>0$, $C>0$ such that for $R>0$, $x\in X$,
	\begin{equation}
		\#\{\gamma\in \Gamma\;|\;\gamma 
		\mathrm{\;nonelliptic}, d_{\gamma}(x)\leq R\}\leq C\exp(cR).
		\label{eq:8.6.1wwx}
	\end{equation}
	By \eqref{eq:ss8.6.1copy}, \eqref{eq:8.6.6}, \eqref{eq:8.6.1wwx}, and using the same arguments as in the proof of \cite[Proposition 3.2]{MP2013raysinger}, we get \eqref{eq:8.6.11ugc}.
\end{proof}

%%%%%%%%%%%%%%%%%%%%%%%%%%%%%%%%%%%%%%%%%%%%%%%%%%%%%%%%%%%%%
\subsection{A proof of Theorem \ref{thm:maintheorem}}\label{section:torsionZ}
In this subsection, we complete our proof to Theorem 
\ref{thm:maintheorem}. Note that every elliptic element $\gamma\in 
\Gamma$ is of finite order, then the Part (2) in Theorem 
\ref{thm:maintheorem} is an easy consequence of Theorem \ref{thm:7.5.3conf}. We only need 
to prove the Part (1). We restate it as follows.

\begin{proposition}\label{prop:8.7.1ugc}
	Let $\Gamma\subset G$ be a cocompact discrete subgroup and set 
	$Z=\Gamma\backslash X$. There exists $c>0$ such that for $d$ large enough, 
	\begin{equation}
		\begin{split}
			\mathcal{T}(Z, 
			F_{d})=&\frac{\mathrm{Vol}(Z)}{|S|}\cP\cI_{X}(F_d)\\
			&+\sum_{[\gamma]\in 
			E^{+}[\Gamma]}\frac{\mathrm{Vol}(\Gamma\cap Z(\gamma)\backslash X(\gamma))}{|S(\gamma)|}  \cP\cE_{X,\gamma}(F_d)+\mathcal{O}(e^{-cd}),
		\end{split}
		\label{eq:8.7.14ugc}
	\end{equation}
	where $E^{+}[\Gamma]=E^{+}[\Gamma]\backslash \{[1]\}$ is the finite set of 
	nontrivial elliptic classes in $[\Gamma]$.
\end{proposition}

\begin{proof}
	By \eqref{eq:8.4.26pps}, we have
	\begin{equation}
		\mathbf{D}^{Z,F_\rd,2}\geq c\rd^2- C.
		\label{eq:8.AAAA}
	\end{equation}
	Then if $\rd$ is large enough, we have
	\begin{equation}
		H^{\cdot}(Z, F_\rd)=0.
		\label{eq:8.4.5ppap}
	\end{equation}
	Then $\mathcal{T}(Z, F_{d})$ can be computed using 
	\eqref{eq:7.8.11bonn}.

	As in \eqref{eq:7.8.7kkkk}, for $t>0$, set 
	\begin{equation}
		b(F_{d},t)=(1+2t\frac{\partial}{\partial 
		t})\mathrm{Tr_s}\bigg[\big(N^{\Lambda^{\bullet}(T^*Z)}-\frac{m}{2}\big)\exp(-t\mathbf{D}^{Z,F_{d},2}/2)\bigg].
		\label{eq:8.8.7kkkk}
	\end{equation}
	As in \cite[Subsection 7.2]{BMZ2015toeplitz}, by \eqref{eq:8.AAAA}, 
	there exist constants $\tilde{c}>0$, $\widetilde{C}>0$ such that for 
	$d$ large enough and for 
	$t>1/d$, 
	\begin{equation}
		|b(F_{d},t)|\leq \widetilde{C}\exp(-\tilde{c} d -\tilde{c} t).
		\label{eq:8.2.6parisconf}
	\end{equation}
	
	By \eqref{eq:7.8.11bonn}, we have
	\begin{equation}
		\mathcal{T}(Z,F_{d})=-\int_{0}^{+\infty} b(F_\rd, t)\frac{dt}{t}.
		\label{eq:8.6.9kkkk}
	\end{equation}
	We rewrite it as follows,
	\begin{equation}
		\mathcal{T}(Z,F_{d})=-\int_{1/d}^{+\infty} b(F_\rd, t)\frac{dt}{t}-\int_{0}^{d} b(F_\rd,t/\rd^2)\frac{dt}{t}.
		\label{eq:8.6.14kkkk}
	\end{equation}
	
	By \eqref{eq:8.2.6parisconf}, there exists $c>0$ such that for 
	$d$ large enough,
	\begin{equation}
		\int_{1/d}^{+\infty} b(F_\rd,t)\frac{dt}{t}=\mathcal{O}(e^{-cd}).
		\label{eq:8.6.11111}
	\end{equation}

	By \eqref{eq:3.5.1ksd}, \eqref{eq:8.6.7mmk}, \eqref{eq:8.8.7kkkk}, we get
	\begin{equation}
		\begin{split}
			b(F_\rd,t)=(1+2t\frac{\partial}{\partial t})\int_Z 
			\frac{1}{|S|}\sum_{\gamma\in\Gamma} v_t(F_\rd,\gamma,z)dz.
		\end{split}
		\label{eq:8.6.15kkkk}
	\end{equation}
	We split the sum in \eqref{eq:8.6.15kkkk} into two parts,
	\begin{equation}
		\sum_{\gamma\in \Gamma, \gamma \text{ elliptic}} + \sum_{\gamma\in \Gamma, \gamma \text{ nonelliptic}}
		\label{eq:8.7.17ugc}
	\end{equation}
	so that we write
	\begin{equation}
		b(F_\rd,t)=b_{\mathrm{ell}}(F_{d},t)+b_{\mathrm{nonell}}(F_{d},t).
		\label{eq:muge}
	\end{equation}
	Similar to the Selberg's trace formula in Subsection \ref{section3.5}, we get
	\begin{equation}
		b_{\mathrm{ell}}(F_{d},t)=\sum_{[\gamma]\in 
		E[\Gamma]}\frac{\mathrm{Vol}(\Gamma\cap Z(\gamma)\backslash 
		X(\gamma))}{|S(\gamma)|}(1+2t\frac{\partial}{\partial t})\cE_{X,\gamma}(F_{d},t).
		\label{eq:8.2.12std}
	\end{equation}
	
	By \eqref{eq:6.4.48ss20} and \eqref{eq:7.5.3std}, the terms in $\cE_{X,\gamma}(F_{d},t)$ 
	are of the form 
	\begin{equation}
		t^{-j+1/2}\exp(-2\pi^{2}t(da'+b')^{2})Q(d),
		\label{eq:8.2.13std}
	\end{equation}
	where $Q(d)$ is a nice exponential polynomial in $d$, and $a',b'\in\R$ with 
	$a'\neq 0$ due to the nondegeneracy of $\lambda$. By \eqref{eq:8.2.13std}, 
	there exists $c>0$ such that for $d$ large enough, 
	\begin{equation}
		\int_{0}^{d}b_{\mathrm{ell}}(F_{d},t/d^{2})\frac{dt}{t}=\int_{0}^{+\infty}b_{\mathrm{ell}}(F_{d},t)\frac{dt}{t}+\mathcal{O}(e^{-cd}).
		\label{eq:8.2.15std}
	\end{equation}
	
	Using Proposition \ref{prop:6.3.1est} and by 
	\eqref{eq:8.2.13std}, we get
	\begin{equation}
		\cP\cE_{X,\gamma}(F_{d})=-\int_{0}^{+\infty}(1+2t\frac{\partial}{\partial t})\cE_{X,\gamma}(F_{d},t)\frac{dt}{t}.
		\label{eq:8.2.16std}
	\end{equation}
	
	Now we consider the contribution from the nonelliptic elements. If $x\in X$, put
	\begin{equation}
		h_t(F_\rd,x)= \frac{1}{|S|}\sum_{\gamma\in \Gamma, \gamma \text{ 
		nonelliptic}} v_t(F_\rd,\gamma,x).
		\label{eq:8.7.18ugc}
	\end{equation}
	Then
	\begin{equation}
		b_{\mathrm{nonell}}(F_{d},t)=(1+2t\frac{\partial}{\partial t})\int_Z 
		h_{t}(F_\rd,z)dz.
		\label{eq:8.2.18std}
	\end{equation}
	Now we prove the following uniform estimates for $x\in X$,
	\begin{equation}
		\int_0^\rd (1+2t\frac{\partial}{\partial t})
		h_{t/\rd^2}(F_\rd,x)\frac{dt}{t}= \mathcal{O}(e^{-c\rd}).
		\label{eq:8.7.19ugc}
	\end{equation}
	
	Indeed, using Lemma \ref{lm:8.6.3ugc} and Proposition \ref{prop:8.6.5ugc}, there exists $C>0$, $c'>0$, $c''>0$ such that if $\rd$ is large enough, $0<t\leq \rd$, then
	\begin{equation}
		|h_{t/\rd^2}(F_\rd,x)|\leq  C \dim(E_\rd) e^{-c' t}\exp(-c''\rd^2/t).
		\label{eq:8.7.20ugc}
	\end{equation}

	Recall that $\dim E_{d}$ is a polynomial in $d$. Then by \eqref{eq:8.7.20ugc}, we have
	\begin{equation}
		\begin{split}
			&\big|\int_0^1 h_{t/\rd^2}(F_\rd,x)\frac{dt}{t}\big|\leq C e^{-c''\rd^2/2}\dim(E_\rd) \int_0^1 e^{-c''\rd^2/2t}\frac{dt}{t}=\mathcal{O}(e^{-c\rd}),\\
			&\big|\int_1^\rd h_{t/\rd^2}(F_\rd,x)\frac{dt}{t}\big|\leq C e^{-c''\rd}\dim(E_\rd) \int_1^\rd e^{-c't}\frac{dt}{t}=\mathcal{O}(e^{-c\rd}).
		\end{split}
		\label{eq:8.7.24ugc}
	\end{equation}
	By \eqref{eq:8.7.20ugc} - \eqref{eq:8.7.24ugc}, we get 
	\eqref{eq:8.7.19ugc}. 
	
	At last, we assembly together \eqref{eq:8.6.14kkkk}, 
	\eqref{eq:8.6.11111}, \eqref{eq:muge}, \eqref{eq:8.2.15std} - \eqref{eq:8.7.19ugc}, we get exactly \eqref{eq:8.7.14ugc}. This completes the proof of our proposition.
\end{proof}

Note that since $\mathcal{T}(Z,F_{d})$ is always real number, 
then \eqref{eq:8.7.14ugc} still holds if we take the real part of 
$\cP\cE_{X,\gamma}(F_{d})$ instead.

\bibliography{References2}
\bibliographystyle{abbrv}

\end{document}